\documentclass[11pt]{article}
\usepackage{epsfig, psfrag, color, amsmath, amscd, geometry}
\input epsf.sty
\epsfverbosetrue %messages will show height and width

\usepackage{latexsym}
\usepackage{amsfonts}
\usepackage{amssymb}

\newcommand{\R}{\mathbb{R}}

\newcommand{\Z}{\mathbb{Z}}
\newcommand{\C}{\mathbb{C}}
\newcommand{\IP}{\mathbb{P}}

\newcommand{\N}{\mathbb{N}}

\newcommand{\IS}{\mathbb{S}}

\newcommand{\rs}{\mbox{$\widehat{\C}$}}

\def\DDD{{\mathcal D}}

\def\HHH{{\mathcal H}}

\def\KKK{{\mathcal K}}

\def\MMM{{\mathcal M}}

\def\MMM{{\mathcal M}}

\def\YYY{{\mathcal Y}}
\def\NNN{{\mathcal N}}

\def\RRR{{\mathcal R}}

\def\UUU{{\mathcal U}}
\def\VVV{{\mathcal V}}

\newcommand{\wt}[1]{\widetilde{#1}}

\newtheorem{thm}{Theorem}[section]
\newtheorem{defn}[thm]{Definition}

\newtheorem{prop}[thm]{Proposition}

\newtheorem{lemma}[thm]{Lemma}
\newtheorem{cor}[thm]{Corollary}

\newcommand{\qed}{\nopagebreak \begin{flushright}%end-of-proof
       \rule{2mm}{2.5mm} \end{flushright}}
\newcommand{\qedspecial}[1]{\nopagebreak \begin{flushright}%end-of-proof
       \rule{2mm}{2.5mm}{\bf #1} \end{flushright}}
\newcommand{\bdry}{\partial}                     %boundary
\newcommand{\id}{\mbox{\rm id}}                  %identity
\newcommand{\cl}{\overline}                      %closure
\newcommand{\Sym}{\mbox{\rm Sym}}                    %short for symmetry
\newcommand{\Aut}{\mbox{\rm Aut}}                        %Aut
\newcommand{\intersect}{\cap}                    %intersection
\newcommand{\union}{\cup}                        %union
\newcommand{\diam}{\mbox{\rm diam}}					         %diameter

\newcommand{\interior}{\mbox{int}}     %interior
\newcommand{\mtwo}[4]                            %2x2 matrices--
{\mbox{$\left(\begin{array}{cc}                  %takes four arguments
#1 & #2 \\
#3 & #4 
\end{array}
\right)$}}

\newcommand{\dettwo}[4]                          %2x2 matrices--
{\mbox{$\left|\begin{array}{cc}                  %takes four arguments
#1 & #2 \\
#3 & #4 
\end{array}
\right|$}}

\newcommand{\pf}{\noindent {\bf Proof: }}

\newcommand{\be}{\begin{enumerate}}
\newcommand{\eb}{\end{enumerate}}
\newcommand{\bi}{\begin{itemize}}
\newcommand{\ib}{\end{itemize}}
\newcommand{\bl}{\begin{list}}
\newcommand{\lb}{\end{list}}

\newcommand{\gap}{\vspace{5pt}}                 %make a space of a blank line

\newcommand{\Stab}{\mbox{\rm Stab}}

\newcommand{\genby}[1]{\mbox{$\langle #1 \rangle$}}

\def\dis{\displaystyle}
\newcommand{\roundness}{\mbox{\rm Round}}

\newcommand{\wtU}{{\widetilde{U}}}
\newcommand{\wtV}{\widetilde{V}}
\newcommand{\wtW}{\widetilde{W}}

\newcommand{\hatU}{\widehat{U}}
\newcommand{\hatV}{\widehat{V}}
\newcommand{\wthatU}{\widetilde{\widehat{U}}}
\newcommand{\wthatV}{\widetilde{\widehat{V}}}
\newcommand{\wtR}{\widetilde{R}}
\newcommand{\wtS}{\widetilde{S}}
\newcommand{\wtG}{\widetilde{\Gamma}}

\newcommand{\wtmho}{\widetilde{\mho}}
\newcommand{\umb}{\mbox{\textsc{Umb}}}
\newcommand{\wtphi}{{\widetilde{\phi}}}
\newcommand{\wtpsi}{{\widetilde{\psi}}}
\newcommand{\roundup}[1]{\lceil #1 \rceil}
\newcommand{\txi}{\tilde{\xi}}
\newcommand{\tv}{\tilde{v}}

\newcommand{\tx}{\tilde{x}}

\newcommand{\tw}{\tilde{w}}

\newcommand{\tzeta}{\tilde{\zeta}}
\newcommand{\XX}{\mbox{$\mathfrak{X}$}}
\newcommand{\YY}{\mbox{$\mathfrak{Y}$}}
\newcommand{\level}[2]{[#1|#2]}
\newcommand{\gromprod}[2]{(#1|#2)}
\newcommand{\divlevel}[2]{\{#1|#2\}}

\def\hde{\hat{\delta}}
\def\hrho{\hat{\rho}}

\begin{document}

\title{Finite type coarse expanding conformal dynamics}

\date{\today}

\author{Peter Ha\"issinsky, Universit\'e de Provence \\ and \\ Kevin M. Pilgrim, Indiana University Bloomington}

\maketitle

\begin{abstract}
We continue the study of noninvertible topological dynamical
systems with expanding behavior.  We introduce the class of
{\em finite type} systems which are characterized by the condition that, 
up to rescaling
and uniformly bounded distortion, there are only finitely many
iterates.  We show that subhyperbolic rational maps and finite subdivision rules (in the sense of Cannon, Floyd, Kenyon, and Parry)
with bounded valence and mesh going to zero are of finite type.
In addition, we show that the limit dynamical system associated 
to a selfsimilar,
contracting, recurrent, level-transitive group action (in the sense of V. Nekrashevych) is of finite type.
The proof makes essential use of an analog of the finiteness of cone types
property enjoyed by hyperbolic groups.  
\end{abstract}

\tableofcontents

%\chapter{Introduction}

\section{Introduction}
\label{ch:intro}

Consider a classical expanding conformal dynamical system on the Riemann sphere $\rs$ equipped 
with its spherical metric--that is, a finitely generated convex cocompact Kleinian group $\Gamma$ 
of M\"obius transformations, or a hyperbolic rational function $f$.  
The chaotic set $X$ (the limit set, in case of a group, or Julia set, in case of a map) 
is {\em quasi-selfsimilar}: given any ball $B \subset X$, there is a group element (or iterate) 
$\psi: B \to X$ which is nearly a similarity and whose image has a definite size independent of $B$.  
This is sometimes known as the {\em principle of the conformal elevator}:  
the dynamics transports geometric features at small scales to large 
(and, by taking inverses, large scales to small) with uniformly bounded distortion.  

The expansive nature of such systems implies that they are {\em finitely generated} 
in the sense of Gromov \cite{gromov:hyperbolic_groups}.  
Roughly, this means that they are quotients of a subshift of finite type.  
Moreover, they are {\em finitely presented}, that is, the equivalence relation defining the quotient is 
again a subshift of finite type.    Other finitely presented systems in which 
the principle of the conformal elevator holds include the action 
of a Gromov hyperbolic group on its boundary equipped with a visual metric.   
For such groups, the finiteness comes from the {\em finiteness of cone types}, 
observed by Cannon in the classical case.  Not all finitely presented systems are conformal, 
since the former include e.g. Anosov maps on tori.  
For details, see  \cite{coornaert:papa:symbolic} and the references therein.

In one-dimensional complex dynamics, there are classes more general than hyperbolic 
for which the principle of the conformal elevator still holds.  For example,  it holds for  {\em subhyperbolic} maps --- 
those whose critical points are either in the Fatou set and converge to attracting cycles, 
or else in the Julia set and are eventually periodic.   
As topological dynamical systems, the set of conjugacy classes of such maps is countable.  
The principle holds as well for the more general {\em semi-hyperbolic} maps --- 
those with neither recurrent critical points nor indifferent cycles.    
The latter class is much larger, however, containing uncountably many distinct topological conjugacy 
classes even in the family of quadratic polynomials.

In earlier work \cite{kmp:ph:cxci}, we introduced a broad class of metric noninvertible discrete-time 
dynamical systems $f: X \to X$ which generalizes the class of semi-hyperbolic rational maps.   
A key role is played by a finite cover $\UUU_0$ of $X$ by connected open sets and the sequence of covers 
$\UUU_{n+1}=f^{-1}(\UUU_n), n=1, 2, 3, \ldots$ obtained by taking components of iterated preimages.  
Such systems, called {\em metric coarse expanding conformal (cxc)}, are defined so that the principle 
of the conformal elevator holds.  A metrization principle holds:  under reasonable hypotheses, a suitably expanding (precisely, 
a ``topologically cxc'')  dynamical system determines a natural quasisymmetry class of metrics 
in which the dynamics is metrically cxc.   
This allows us to adapt techniques from classical conformal dynamics to study
noninvertible topological dynamical systems.
Precise statements  are given in Chapter 2 below.   
A distinguishing feature of such systems is that, like for subhyperbolic rational maps, 
the dynamics need not be locally injective on the chaotic set.  In many respect, our results
in \cite{kmp:ph:cxci} suggest that our class of cxc maps share many properties with hyperbolic groups, thus
extending Sullivan's dictionary to this wider setting.  However, the fact that 
hyperbolic groups are automatic, which is a consequence of the finiteness of cone types for such groups, 
does not seem to have a counterpart for
cxc dynamics in general.

In the present work, we single out a subset of the metric cxc systems  comprising maps 
which satisfy finiteness features that we believe to be analogous to the finiteness of cone types.  
These dynamical systems are characterized by the existence of what we call a {\em dynatlas}---a finite set $\MMM$ of local 
model maps $g_m: \wtV_m \to V_m, m \in \MMM$ such that the restriction of any (suitable) iterate 
$f^k: \wtU \to U$, $\wtU \in \UUU_{n+k}, U \in \UUU_n$ is, after rescaling, nearly isometric to one of 
the model maps $g_m$.   We term such systems {\em cxc systems of finite type}.  The set of finite type rational
maps is exactly the set of so-called subhyperbolic rational maps. 
Perhaps over-optimistically, we suspect that such finite type cxc systems 
are finitely presented.  
The lack of available (to us) general techniques from dynamics and the few topological assumptions 
made on the underlying space $X$ makes verification of this suspicion difficult.   

Our main results identify two important natural sources of examples of cxc systems of finite type 
which do not necessarily arise from classical Riemannian conformal dynamics.  
\gap

\noindent{\bf  Finite subdivision rules.}  We study those  {\em finite subdivision rules} (fsr's) of Cannon, Floyd, 
and Parry \cite{cfp:fsr} which have {\em bounded valence} and {\em mesh going to zero}; they yield postcritically 
finite branched, or {\em Thurston}, maps $f: S^2 \to S^2$ of the two-sphere to itself without periodic critical points. 
Such fsr's arise naturally in geometrization questions such as the characterization of rational functions and 
Cannon's conjecture concerning hyperbolic groups with two-sphere boundary.   
These maps turn out to be topologically cxc, and the metrization principle in \cite{kmp:ph:cxci}, 
shows that such fsr's naturally yield metric cxc systems.  

We show (Corollary \ref{cor:fsr_are_ft}) that in fact, the corresponding metric cxc dynamical systems are of finite type.  
More is actually true:  in a natural metric, up to similarity (not quasisimilarity), there are only finitely many tiles.  
This is derived from a general metrization principle (Theorem \ref{thm:top_implies_metric_ft}) 
which asserts that a topological dynamical system with suitable finiteness properties admits a natural metric 
in which the dynamics is metrically of finite type.  
\gap

\noindent{\bf Selfsimilar groups.}  In \cite{nekrashevych:book:selfsimilar} a general theory of so-called 
{\em selfsimilar groups} is developed which connects group theory and dynamical systems in both directions.  
In one direction, to a dynamical system, $f: X \to X$ one may associate a selfsimilar group action, 
its {\em iterated monodromy group}.  Under reasonable expansion hypotheses, this action is {\em contracting} 
and {\em recurrent}.  In the other direction, to a selfsimilar contracting recurrent action is associated 
a topological dynamical system $\bdry F_\Sigma: \bdry \Sigma \to \bdry \Sigma$.  
The underlying space $\bdry \Sigma$ is the boundary of an infinite Gromov hyperbolic graph, 
$\Sigma$, called the {\em selfsimilarity complex} associated to the action.  
The map $\bdry F_\Sigma$ is induced by a graph endomorphism $F_\Sigma: \Sigma \to \Sigma$.   
Under appropriate topological regularity and expansion hypotheses, the circle of ideas can be completed.  
That is, given $f: X \to X$, the iterated monodromy group associated to $f$ yields a selfsimilar recurrent contracting action, 
and the associated topological dynamical system $\bdry F_\Sigma: \bdry \Sigma \to \bdry \Sigma$ is conjugate to $f: X \to X$.   

We prove a finiteness principle (Theorem \ref{thm:finiteness_principles}) for the map $F_\Sigma$ analogous 
to the finiteness of cone types for a hyperbolic group.   This is used to conclude 
(Theorem \ref{thm:boundary_dynamics_finite_type}) that, in a visual metric on the boundary, 
such dynamical systems are metric cxc of finite type.  As a corollary, we obtain that the quasi-isometry type 
of the selfsimilarity complex $\Sigma$ and, therefore, the quasisymmetry class of metrics on its boundary are 
invariants of the induced topological dynamical system, hence of the group action.  
In particular, the Ahlfors regular conformal dimension is a numerical invariant of this group action.  
It would be interesting to know how this invariant is related to other such quantities, e.g. contraction coefficients, 
growth functions, etc.  
\gap

\noindent{\bf Organization.} In \S 2, we define the class of topological and metric cxc systems, 
concluding with the formal definition of metric finite type.  
In \S 3, we give natural first classes of examples, starting with unbranched systems.  
Next we prove that subhyperbolic maps are finite type.  
The proof given motivates the more abstract argument used to prove the metrization theorem,  
Theorem \ref{thm:top_implies_metric_ft}.  \S\S 4, 5, and 6 present the connections to selfsimilar groups.  
\gap

\noindent{\bf Acknowledgements.}  The authors would like to thank Jim Cannon and Volodia Nekrashevych for useful conversations, and the anonymous referee for numerous detailed comments. 
Both authors were partially supported by the  project ANR ``Cannon'' (ANR-06-BLAN-0366).  
The second author was supported by NSF grant DMS-0400852.  

%\input sec1.tex

%\chapter{Definition}

\section{Definition and first properties}
\label{secn:finite_type}

We first recall some definitions and results from \cite{kmp:ph:cxci}.   We  then define
the class of metric finite type maps, and show they are coarse expanding conformal.

\subsection{Finite branched coverings} 
\label{secn:fbc}

Suppose $X, Y$ are locally compact Hausdorff spaces, and let $f: X \to Y$ be a finite-to-one continuous map.  
The {\em degree} of
$f$ is  
\[ \deg(f) = \sup\{\# f^{-1}(y): y \in Y\}.\] For $x \in X$, the {\em
local degree of $f$ at $x$} is
\[ \deg(f; x) = \inf_U \sup\{\#f^{-1}(\{z\}) \intersect U : z \in f(U)\}\] where $U$
ranges over all neighborhoods of $x$.  

\begin{defn}[finite branched covering]
The map $f$ is a {\em finite branched covering} (abbrev. fbc) provided $\deg(f)<\infty$ and 
\bi 
\item[(i)] $$
\sum_{x \in f^{-1}(y)} \deg(f;x) = \deg f$$ holds for each $y \in Y$;
\item[(ii)] for every $x_0 \in X$, there are compact neighborhoods
$U$ and $V$ of $x_0$ and $f(x_0)$ respectively such that
\[ \sum_{x \in U, f(x)=y} \deg(f; x) = \deg(f; x_0)\]
for all $y\in V$.\ib  
\end{defn}

When $X, Y$ are connected, locally connected, and compact and $f: X \to Y$ is finite-to-one, 
closed, open, and continuous, the second condition (ii) is implied by the first; see \cite[Lemma 2.5]{edmonds:fbc}. 

The composition of  fbc's is an fbc, and the degrees of fbc's multiply under compositions.  
In particular, local degrees of fbc's multiply under compositions.

Condition (ii) implies that if $x_n \to x_0$, then $\deg(f; x_n) \leq \deg(f; x_0)$.  It follows that  
the {\em branch set} $B_f = \{x \in X : \deg(f; x)>1\}$ is closed.    
The set of {\em branch values} is defined as $V_f = f(B_f)$.

\begin{lemma}Let $X,Y$ be Hausdorff locally compact topological spaces. An fbc $f:X\to Y$ of
degree  $d$ is onto, proper and open:  the inverse image of a compact subset is compact
and the image of an open set is open. Furthermore, $B_f$ and $V_f$ are nowhere dense.
\end{lemma}

Many arguments are done using pull-backs of sets and restricting to connected components. It is therefore necessary
to work with fbc's defined on sets $X$ and $Y$ enjoying more properties. 
When  $X$ and $Y$, in addition to being locally compact and Hausdorff, are  assumed locally connected,  the following fundamental facts are known (cf. \cite{edmonds:fbc}). 
\bi
\item If $V \subset Y$ is open and connected,
and $U \subset X$ is a connected component of $f^{-1}(V)$, then $f|U: U \to V$ 
is an fbc as well.  

\item If $y \in Y$, and $f^{-1}(y) = \{x_1, x_2, \ldots, x_k\}$, then there exist arbitrarily small connected open neighborhoods $V$ of $y$ such that 
\[ f^{-1}(V)=U_1 \sqcup U_2 \sqcup \ldots \sqcup U_k\]
is a disjoint union of connected open neighborhoods $U_i$ of $x_i$ such that $f|U_i: U_i \to V$ is an fbc of degree $\deg(f; x_i)$, $i=1, 2, \ldots, k$.  

\item if $f(x)=y$, $\{V_n\}$ is sequence of nested open connected sets with 
$\intersect_n V_n = \{y\}$, and if $\wtV_n$ is the component of 
$f^{-1}(V_n)$  containing $x$, then $\intersect_n \wtV_n = \{x\}$. 

\ib

\subsection{Topological cxc systems}  
\label{secn:def_top_cxc}   In this section, we state the topological axioms underlying the definition of a cxc system.

Let $\XX_0, \XX_1$ be Hausdorff, locally compact, locally connected topological spaces, each with finitely many connected components.  We further assume that  $\XX_1$ is an open subset of $\XX_0$ and that $\cl{\XX_1}$ is compact in $\XX_0$.    
Note that this latter condition implies that if $\XX_0 = \XX_1$, then $\XX_0$ is compact.  

Let $f: \XX_1 \to \XX_0$ be a finite branched covering map of degree $d \geq 2$, and for $n \geq 0$ put 
\[ \XX_{n+1} = f^{-1}(\XX_n).\]
Then $f: \XX_{n+1} \to \XX_n$ is again an fbc of degree $d$ and since $f$ is proper, $\cl{\XX_{n+1}}$ is compact in $\XX_n$, hence in $\XX_0$.  

The {\em nonescaping set}, or {\em repellor}, of $f: \XX_1 \to \XX_0$ is
\[ X = \{ x \in \XX_1 | f^n(x) \in \XX_1 \;\forall n > 0\} = \bigcap_n \cl{\XX_n}.\]
We make the technical assumption that the restriction $f|X: X \to X$ is also an fbc of degree equal to $d$.  
This implies that  $\#X \geq 2$.   Also,  $X$ is {\em totally invariant}:  $f^{-1}(X) = X = f(X)$.  

The following is the essential ingredient in this work.  Let $\UUU_0$ be a finite cover of $X$ by open, connected subsets of $\XX_1$ whose intersection with $X$ is nonempty.  A  {\em preimage} of a connected set $A$ is defined as a connected component of $f^{-1}(A)$. Inductively, set $\UUU_{n+1}$ to be the open cover whose elements $\wtU$ are preimages of elements of $\UUU_n$.   We denote by $\mathbf{U} = \union_{n \geq 0} \UUU_n$ the collection  of all such open sets thus obtained.  For $U \in \mathbf{U}$, the {\em level} of $U$, denoted $|U|$, is the value of $n$ for which $U \in \UUU_n$. 

We say $f: (\XX_1,X) \to (\XX_0, X)$ is {\em topologically coarse expanding conformal with repellor $X$} (topologically cxc, for short) provided there exists  a finite covering $\UUU_0$ as above, such that the following axioms hold.  

\be  

\item {\bf [Expans]}  The mesh of the coverings $\UUU_n$ tends to zero as $n \to \infty$.  That is, for any finite open cover $\YYY$ of $X$ by open sets of $\XX_0$, there exists $N$ such that for all $n \geq N$ and all $U \in \UUU_n$, there exists $Y \in \YYY$ with $U \subset Y$.  

\item {\bf [Irred]}  The map $f: \XX_1 \to \XX_0$ is {\em locally eventually onto near $X$}:  for any $x \in X$ and any neighborhood $W$ of $x$ in $\XX_0$, there is  some $n$ with $f^n(W) \supset X$

\item {\bf [Deg]} The set of degrees of maps of the form $f^k|\wtU: \wtU \to U$, where $U \in \UUU_n$, $ \wtU \in \UUU_{n+k}$, and $n$ and $k$ are arbitrary, has a finite maximum, denoted $p$.  
\eb

Axiom [Expans] is equivalent to saying that, when $\XX_0$ is a metric space, the diameters of the elements of $\UUU_n$ tend to zero as $n \to \infty$.  
\gap

The elements of $\UUU_0$ will be referred to as {\em level zero good open sets}.  While as subsets of $\XX_0$ they are assumed connected, their intersections with the repellor $X$ need not be.   Also, the elements of $\mathbf{U}$, while connected, might nonetheless be quite complicated topologically--in particular they need not be contractible.  

If $\XX_0 = \XX_1 = X$, then the elements of $\mathbf{U}$ are connected subsets of $X$. \gap

\gap

\noindent{\bf Conjugacy.}  
Suppose $f: \XX_1 \to \XX_0$ and $g: \YY_1 \to \YY_0$ are fbc's with repellors $X$, $Y$ as in the definition of topologically cxc.  A homeomorphism $h: \XX_0 \to \YY_0$ is  called a {\em conjugacy}  if it makes the diagram
%commutative diagram 
\[ 
\begin{array}{ccc}
(\XX_1, X) & \stackrel{h}{\longrightarrow} & (\YY_1, Y)
\\
f \downarrow & \; & \downarrow g \\
(\XX_0, X) & \stackrel{h}{\longrightarrow} &(\YY_0, Y) \\
\end{array}
\]
commute.  (Strictly speaking, we should require only that $h$ is defined near $X$; however, we 
will not need this more general point of view here.)  

It is clear that the property of being topologically cxc is closed under conjugation.

\gap

\subsection{Metric cxc systems} 
\label{subsecn:metric_cxc}
In this section, we state the definition of metric cxc systems; we will henceforth drop the adjective, metric.  
Given a metric space, $|a-b|$ denotes the distance between $a$ and $b$, 
and $B(a,r)$ denotes the open ball of radius $r$ about $a$.
\gap

\noindent{\bf Roundness.} Let $Z$ be a metric space and let $A$ be a bounded, proper  subset of $Z$ with nonempty interior.  Given $a \in \interior(A)$, define the {\em outradius} of $A$ about $a$ as  
\[ L(A,a)=\sup\{|a-b|: b \in A\}\]
and the {\em inradius} of $A$ about $a$ as 
\[ \ell(A,a)=\sup\{ r : r \leq L(A,a) \; \mbox{ and } \; B(a,r) \subset A\}.\]
The condition $r \leq L(A,a)$ is necessary to guarantee that the outradius is at least the inradius. 
For instance, if $A=\{a\}$ is an isolated point
in $Z$, then 
$$ L(A,a)=0 < \sup\{r: B(a,r)\subset A\}= d(a,Z\setminus A) \,.$$ 
The outradius is intrinsic--it depends only on the restriction of the metric to $A$.  
In contrast, the inradius depends on how $A$ sits in $Z$. 
The {\em roundness of $A$ about $a$} is defined as 
\[ \roundness(A,a) = L(A,a)/\ell(A,a) \in [1, \infty).\]
One says $A$ is {\em $K$-almost-round} if $\roundness(A,a) \leq K$ for some $a \in A$, and this implies that for some $s>0$, 
\[ B(a,s) \subset A \subset B(a,Ks).\]
\gap

Isometric open embeddings which are not surjective may distort roundness:\\
\gap 

\noindent{\bf Example.} Consider
in $\R^2$ the metric spaces $X=\R\times\{0\}$ and $Y=X\cup \left( \{0\}\times [c,\infty)\ \right)$
for some constant $c>0$.
Then the inclusion $X\subset Y$ is an isometric and open embedding,
and, for any interval $(-r,r)\times\{0\}$, $r>c$, centered at the origin,
its roundness at the origin in $X$ is $1$, but $r/c$ in $Y$.
\gap

\noindent{\bf Metric cxc systems.} Suppose we are given a topological cxc system $f:\XX_1 \to \XX_0$ with level zero good neighborhoods $\UUU_0$,  and that $\XX_0$ is now endowed with a metric  
compatible with its topology.  The resulting metric dynamical system equipped with the covering $\UUU_0$  is called {\em
coarse expanding conformal}, abbreviated cxc, provided there exist
\bi
\item continuous, increasing embeddings $\rho_{\pm}:[1,\infty) \to [1,\infty)$, the {\em forward and backward
roundness distortion functions}, and

\item increasing homeomorphisms $\delta_{\pm}:\R_+ \to\R_+$, the {\em forward and backward relative
diameter distortion functions}
\ib
satisfying the following axioms:
\be
\setcounter{enumi}{3}

\item {\bf [Round]} $(\forall n, k)$ and for all
\[ U \in \UUU_n, \;\;\wtU \in \UUU_{n+k}, \;\; \tilde{y} \in \wtU, \;\;  
y
\in U\]
if
\[ f^{\circ k}(\wtU) = U, \;\;f^{\circ k}(\tilde{y}) = y \]
then the {\em backward roundness bound}
\begin{equation}
\label{eqn:backward_roundness_bound}
 \roundness(\wtU, \tilde{y}) <
\rho_-(\roundness(U,y))
\end{equation}
and the {\em forward roundness bound}
\begin{equation}
\label{eqn:forward_roundness_bound}
\roundness(U, y) <
\rho_+(\roundness(\wtU,
\tilde{y})).
\end{equation}
hold.

\item {\bf [Diam]} $(\forall n_0, n_1, k)$ and for all
\[ U \in \UUU_{n_0}, \;\;U' \in \UUU_{n_1}, \;\;\wtU \in \UUU_{n_0+k},
\;\;\wtU'
\in
\UUU_{n_1+k}, \;\; \wtU' \subset \wtU, \;\; U' \subset U\]
if
\[ f^k(\wtU) = U, \;\;f^k(\wtU') = U'\]
then
\[ \frac{\diam\wtU'}{\diam\wtU} < \delta_-\left(\frac{\diam U'}{\diam
U}\right)\]
and
\[ \frac{\diam U'}{\diam U} < \delta_+\left(\frac{\diam \wtU'}{\diam
\wtU}\right)\]

\eb

For the Axiom [Diam], note that only the values smaller than $1$ are relevant
for $\delta_{\pm}$.
The Axiom [Expans] implies that the maximum diameters of the elements of $\UUU_n$ tend to zero uniformly in $n$.   Since $\UUU_0$ is assumed finite, each covering $\UUU_n$ is finite, so for each $n$ there is a  minimum diameter of an element of $\UUU_n$.  By \cite[Prop. 2.4.1]{kmp:ph:cxci}, the space $X$ is perfect, that is, contains no isolated points.  By assumption, each $U \in \mathbf{U}$ contains a point of $X$.  So, it contains many points of $X$, and so has positive diameter.   Hence there exist decreasing positive sequences $c_n, d_n \to 0$ such that the {\em diameter bounds} hold:
\begin{equation}
\label{eqn:diameter_bounds}
0 < c_n \leq  \inf_{U \in \UUU_n} \diam U \leq \sup_{U \in \UUU_n} \diam U \leq  d_n.
\end{equation}

{\noindent\bf Canonical gauge.} A homeomorphism $h$ between metric spaces $(X,d_X)$ and $(Y,d_Y)$ is 
called {\em quasisymmetric} provided 
there exists a homeomorphism $\eta: [0,\infty) \to [0,\infty)$ such that 
$$d_X(x,a) \leq td_X(x,b) \implies d_Y(h(x),h(a)) \leq \eta(t)d_Y(h(x),h(b))$$ 
for all triples of points $x, a, b \in X$ and all $t \geq 0$.    

\gap

In \cite{kmp:ph:cxci} the following results are proved.

\begin{thm}[Invariance of cxc]
\label{thm:invariance_of_cxc}
Suppose $f: (\XX_1, X) \to (\XX_0, X)$ and $g:(\YY_1, Y) \to (\YY_0, Y)$ are two topological cxc systems which are conjugate via a homeomorphism $h: \XX_0 \to \YY_0$, where $\XX_0$ and $\YY_0$ are metric spaces. 
\be
\item If $f$ is metrically cxc and $h$ is quasisymmetric, then $g$ is metrically cxc, quantitatively.
\item If $f, g$ are both metrically cxc, then $h|_X: X \to Y$ is quasisymmetric, quantitatively. 
\eb

\end{thm}

The {\em conformal gauge} of a metric space $(X,d_X)$ is the set of metric spaces quasisymmetric to $X$.  The previous theorem shows that the gauge of $X$ depends only on the conjugacy class of $f: \XX_1 \to \XX_0$.  This is not quite intrinsic to the dynamics on $X$.  However, if $\XX_0=\XX_1=X$ then one has the following metrization theorem.

\begin{thm}[Canonical gauge]\label{topcxc_cxc1} 
If $f:X\to X$ is a topological cxc map, where $\XX_1=\XX_0=X$, then
there exists a unique conformal gauge on $X$ defined by a metric $d$ such that 
$f:(X,d)\to (X,d)$ is metric cxc.\end{thm}

The metric may be defined as follows (see \cite[\S 3.1]{kmp:ph:cxci} for details).

Suppose $f: \XX_1 \to \XX_0$ is topologically cxc with respect to an open covering 
$\UUU_0$ as in Theorem \ref{topcxc_cxc1}.  
Let $\UUU_{-1}$ be the open covering consisting of a single element, 
denoted $o$, defined as the union of the components of the interior of $\XX_1$.
Let $\Gamma$ be the graph whose edges are the elements of $\mathbf{U}$, together with the distinguished root vertex $o$.  
The set of edges is defined as a disjoint union of two types:  
{\em horizontal} edges join elements $U_1, U_2 \in \UUU_n$ if and only if 
$X\cap U_1 \intersect U_2 \neq \emptyset$, 
while {\em vertical} edges join elements $U \in \UUU_n, V\in \UUU_{n\pm 1}$ at consecutive levels if and only 
if $X\cap U \intersect V \neq \emptyset$.  
There is a natural map $F: \Gamma \to \Gamma$. 
Near the root $o$, it sends $o$ to $o$, collapses vertical edges meeting $o$ to $o$, 
sends each vertex $U \in \UUU_0$ to $o$,  collapses horizontal edges joining points of $\UUU_0$ also to $o$, and
send the edge between vertices $(U_0,U_1)\in \UUU_0\times \UUU_1$ to the edge $(o,f(U_1))$.
At every other vertex, the map $F$ is defined so that $F(U)=f(U)$; it extends to a graph homomorphism 
on the subgraph spanned by vertices at level one and higher.  

\gap

The usual graph-theoretic combinatorial distance $d( , )$ gives rise to a path metric on $\Gamma$ such 
that edges are isometric to unit intervals.  

One may define its compactification in the following way.
Fix $\varepsilon > 0$.  For $x \in \Gamma$ 
let $\varrho_\varepsilon (x)=\exp(-\varepsilon d(o,x))$.  Define  
a new metric $d_\varepsilon$ on $\Gamma$ by
\[ d_\varepsilon(x,y) = \inf_\gamma \int_\gamma \varrho_\varepsilon \ ds\]
where as usual the infimum is over all rectifiable curves in $\Gamma$ joining $x$ to $y$. 
The resulting metric space $\Gamma_\varepsilon$ is incomplete. Its complement in its
completion defines the boundary $\bdry_\varepsilon \Gamma$.  

\gap

Since $\Gamma$ is Gromov-hyperbolic \cite[Thm 3.3.1]{kmp:ph:cxci}, 
the metric $d_\varepsilon$ is a visual metric and the boundary $\bdry_\varepsilon\Gamma$ is the hyperbolic boundary
as soon as $\varepsilon$ is sufficiently small.  

In that case,   the boundary $\bdry_\varepsilon \Gamma$ 
coincides with the set of the classes of asymptotic geodesic rays (in the metric $d$) emanating from $o$,
and is homeomorphic to $X$. More precisely,  
the natural map $X \to \bdry_\varepsilon \Gamma$ given by $\phi(z) = \lim U_n(z)$ 
where $z \in U_n \in \UUU_n \subset \Gamma_\varepsilon$ is well-defined and a homeomorphism  
conjugating $f$ on $X$ to the map on $\bdry_\varepsilon \Gamma$ induced by the cellular map $F$.  
The conformal gauge of the metric $d_\varepsilon$ defines the canonical gauge of $f$.

We denote by $B_\varepsilon(x,r)$ the ball of radius $r$ in the metric $d_\varepsilon$.  
The metric $d_\varepsilon$ has the following properties (\cite[Prop. 3.2.2, 3.2.3]{kmp:ph:cxci}):
\bi
\item $F^k(B_\varepsilon(x,r))=B_\varepsilon(F^k(x), \exp(k\varepsilon)r)$ if $r$ is sufficiently small, and
\item if $F^k|_{B_\varepsilon(x,4r)}$ is injective, then $f^k|_{B_\varepsilon(x,r)}$ is a similarity with factor $\exp(k\varepsilon)$.
\ib

\gap

\subsection{Maps of metric finite type}

A random cxc system may appear rather inhomogeneous:  for example, one may 
conjugate $z \mapsto z^2$ on the standard Euclidean circle $\IS^1$ with a horrible quasisymmetric map which is  the identity off a neighborhood of the preperiodic point $-1$.  
In many cases, however, one finds an extra degree of homogeneity present in a cxc 
system.  

\gap

\noindent{\bf Quasisimilarities.}  We will find the concept of quasisimilarity
useful 
for capturing the notion that a family of maps is nearly a collection of
similarities.  

\begin{defn}[Quasisimilarity]
\label{defn:quasisimilarity}
Let $h:X \to Y$ be a homeomorphism between metric spaces.  
We say that
$h$ is a {\em $C$-quasisimilarity} if there is some constant $\lambda>0$ such that
\[ \frac{1}{C} \leq \frac{|h(a)-h(b)|}{\lambda |a-b|} \leq C \]
for all $a,b \in X$.  
A family $\HHH$ of homeomorphisms (perhaps defined on
different spaces) consists of {\em uniform quasisimilarities} 
if there exists a
constant $C$ (independent of $h$) such that each
$h \in \HHH$  
is a $C$-quasisimilarity.  
\end{defn}

We will speak of $(C,\lambda)$-quasisimilarity if we want to emphasise
the constant $\lambda$.

\gap

\noindent{\bf Example.}  Fix $0 < r < 1$.  If $f: \Delta \to \C$ is an analytic function which is injective on the unit disk $\Delta$, then the Koebe distortion principle implies that the restriction of $f$ to any smaller disk $\{z : |z| < r\}$ is a $(C_r, \lambda)$-quasisimilarity, where $\lambda = |f'(0)|$ and $C_r$ is independent of $f$.  
\gap

One establishes easily:
\be
\item If $\lambda = (1/C')\lambda'$ then a
$(C,\lambda)$-quasisimilarity is also a $(CC', \lambda')$-quasisimilarity.  
\item The inverse of a
$(C,\lambda)$-quasisimilarity is a $(C, 1/\lambda)$-quasisimilarity.
\item A $(C,\lambda)$-quasisimilarity distorts ratios of diameters by at most 
the factor $C^2$.
\eb

\begin{defn}[Dynatlas] 
A {\em dynatlas} is a finite set $\MMM$ whose elements, called {\em local model maps}, 
are fbc's $g: \wtV \to V$, where $\wtV, V$ are locally compact, connected metric spaces of diameter $1$.  
Given $m \in \MMM$, the corresponding local model map is denoted $g_m: \wtV_m \to V_m$.  
The sets $\wtV_m, V_m$, $m \in \MMM$, are called the {\em model open sets}, 
and the collection of model open sets is denoted $\VVV$, i.e. 
$\VVV = \{\wtV_m : m \in \MMM\} \union \{V_m : m \in \MMM\}$.  \end{defn}

\begin{defn}[metric finite type]
\label{defn:finite_type}
Let $f:(\XX_1,X) \to (\XX_0, X)$ be an fbc, and suppose $\XX_0$ is equipped with a metric compatible with its topology.  
Let $\UUU_0$ be a finite covering of $X$ by open connected subsets of $\XX_1$, and let $\mathbf{U} = \{\UUU_n\}$ be the sequence of coverings obtained by pulling back $\UUU_0$ under iterates of $f$.  

We say $f:(\XX_1, X) \to (\XX_0, X)$ is {\em metric finite type with respect to $\UUU_0$}  if the axioms [Expans]
and [Irred] hold, and if there exists a dynatlas $(\VVV, \MMM)$ and  
a  constant $C \geq 1$ with the following properties:
\be
\item $\forall U \in \mathbf{U}$, there exists $V \in \VVV$ and a 
$C$-quasisimilarity $\psi_U: U \to V$,
\item  every $V \in \VVV$ arises in this way---that is, $\forall V \in \VVV$, there exists $U \in \mathbf{U}$ and a $C$-quasisimilarity $\psi_U: U \to V$.  
\item Whenever $\wtU, U \in \mathbf{U}$ and $f^k: \wtU \to U$, the map   
\[  g_{\wtU, U} := \psi_U \circ f^k|\wtU \circ \psi_\wtU^{-1} \in \mathcal{M}.\]  
\eb 
\end{defn}

We may think of the set of maps $\psi_U$ as a set of ``local coordinate charts'' 
which comprise a family of uniform $C$-quasisimilarities.  
The property of being metric finite type may then be characterized as follows:  
up to quasisimilarity, there are only finitely many local models for the dynamics over 
the elements of the finite good cover $\UUU_0$.   
The second condition is a convenient minimality assumption: 
if a dynatlas satisfies the first and third, then by removing extraneous elements, it will satisfy the second.   

The property of being metric finite type is not invariant under quasisymmetric
conjugacies.  

\begin{thm}[metric finite type implies cxc]
\label{thm:finite_type_implies_cxc}
If $f: \XX_1 \to \XX_0$ is of metric finite type with respect to $\UUU_0$, 
then it is metric coarse expanding conformal with respect to $\UUU_{n_1}$ 
for some $n_1 \ge 0$.  
\end{thm}

The proof, which occupies the remainder of this subsection, is essentially straightforward except for one subtlety.  Since roundness is not an intrinsic quantity, care must be taken to show that the non-surjective embeddings $\psi^{-1}$ do not distort roundness too  much.   

\gap

We first establish some properties of quasisimilarity embeddings, that is, maps which are quasisimilarities onto their images.

\gap

{\noindent\bf Properties of quasisimilarities.}
We assume here that $V$, $Z$ are connected metric spaces,
$Z$ is complete and
$\psi : V\to Z$ is a $(C,\lambda)$-quasisimilarity embedding
with $\psi(V)=U$ open, and $V$ bounded. 
\gap

The following is easily verified:

\begin{prop}\label{prop:easy_properties}
\be
\item Both $\psi$ and $\psi^{-1}$ extend as a quasisimilarity
between the completion of $V$ and the closure of $U$.
\item For any open subsets $W_1, W_2\subset V$,
$$\frac{1}{C^2}\frac{\diam W_1}{\diam W_2}\le \frac{\diam\psi(W_1)}{\diam
\psi (W_2)}\le C^2\cdot \frac{\diam W_1}{\diam W_2}\,.$$
\item If $\psi$ is onto, then for $W\subset V$ and $x\in W$,
$$\roundness (\psi(W),\psi(x))\le C^2\roundness (W,x)\,.$$
\eb
\end{prop}

Next, we establish roundness distortion bounds for open embeddings which need not be onto.

\begin{prop}\label{prop:qsim_control}
Let $W\subset V$ and $x\in W$. The following hold 
\be
\item If $\diam W\le (1/2C^2)\diam V$, then
$$\roundness (\psi(W),\psi(x))\le\max\{C^2\roundness(W,x),\roundness(U,\psi(x))\}\,.$$
\item If $\roundness (U,\psi(x))\le R$ and  $\diam W\le (1/2RC^2)\diam V$, then
$$\roundness (\psi(W),\psi(x))\le C^2\roundness(W,x)\,.$$
\eb\end{prop}

\pf\be
\item The definition of roundness implies that
$$\roundness (\psi(W),\psi(x))\le \frac{\diam \psi(W)}{\ell_Z(\psi(W),\psi(x))}\,.$$
Since $\psi$ is a $(C,\lambda)$-quasisimilarity and $B(x, \ell_V(W,x)) \subset W$, we have (recalling $U=\psi(V)$) that 
$$B\left(\psi(x),\frac{\lambda \ell_V(W,x)}{C}\right)\cap U\subset \psi(W)\,.$$
Recall that by definition, $B(\psi(x), \ell_Z(U,\psi(x))) \subset U$.
\gap

We now consider two cases.  
If on the one hand $\ell_Z(U,\psi(x))\le \lambda \ell_V(W,x)/C$, then
$$B(\psi(x),\ell_Z(U,\psi(x)))=B(\psi(x),\ell_Z(U,\psi(x)))\cap U\subset\psi(W)$$
which implies
\[ \ell_Z(\psi(W),\psi(x)) = \ell_Z(U,\psi(x))\,.\]
But $L_Z(\psi(W),\psi(x)) \leq \diam\psi(W)$ and by assumption $\diam \psi(W)\le (1/2)\diam U$, so 
$$\roundness (\psi(W),\psi(x))\le \frac{\diam U}{2\ell_Z(U,\psi(x))}
\le\roundness(U,\psi(x))\,.$$

\gap

If on the other hand $\ell_Z(U,\psi(x))\ge \lambda \ell_V(W,x)/C$, then $$B\left(\psi(x), \frac{\lambda \ell_V(W,x)}{C}\right) \subset U$$ and so 
$$\roundness (\psi(W),\psi(x))\le C^2 \roundness(W,x)\,.$$

\item Since $\roundness(U,\psi(x))\le R$, it follows that
$$\ell_Z(U,\psi(x))\ge \frac{L_Z(U,\psi(x))}{R}\ge \frac{\diam U}{2R}\ge
\frac{\lambda\diam V}{2 RC}\,.$$
But $\diam W\le (1/2RC^2)\diam V$ so that
$$\diam \psi(W)\le \lambda C \diam W \le \frac{\lambda \diam V}{2RC}\le\ell_Z(U,\psi(x))\,.$$
Therefore, $$\psi(W)\subset B(\psi(x),\ell_Z(U,\psi(x)))\subset U$$
and (1) above implies $$\roundness(\psi(W),\psi(x))\le C^2\roundness (W,x)\,.$$
\qed\eb

\gap

We now assume that we are given dynamical system $f:(\XX_1, X) \to (\XX_0, X)$ of metric finite type.

\gap

{\noindent\bf Diameter bounds.}  
Since all model maps appear in the dynamics of $f$, it follows
that they are all uniformly continuous. Since the set $\MMM$ is finite,
there exists a common modulus of continuity function $\hde_+$:  
this is a self-homeomorphism of $\R_+$ such that, given $m\in \MMM$, any $x,y\in\wtV_m$,
$$|g_m(x)-g_m(y)|\le \hde_+(|x-y|)\,.$$

Define the function 
$\hde_-:[0,1]\to \R_+$ as follows. Given $m\in\MMM$, $W\subset V_m$,
we let $\hbox{presize}(W)=\max \diam \wtW$ where
the supremum ranges over all connected components of $g_m^{-1}(W)$.
For $r\in [0,1]$, set 
$$\hde_-(r)= \max_{m\in\MMM}\sup\{\hbox{presize}(W),\ W\subset V_m,\ \diam W\le r\}\,.$$
By taking its upper hull, we may assume that $\hde_-$ is a homeomorphism on its image, and
we may extend it to be a homeomorphism of $\R_+$ as well.

For future use, we note that, since our maps are open, 
there exists a positive increasing function $\hrho$ such that,
for any $m\in\MMM$, for any $x\in\wtV_m$, and any $r<1$, $g_m(B(x,r))$
contains the ball $B(g_m(x),\hrho(r))$.  Indeed, we may use representatives
given by appropriate restrictions of iterates of $f$: there is some $n_0$ such
that every model map appears, up to a $C$-quasisimilarity,
as $f^k:\wtU\to U$, with $\wtU,U\in \cup_{0\le j\le n_0}\UUU_j$, $1\le k\le n_0$.
Since $f$ and its iterates
are open maps, for each $x\in X$, $r>0$ and $n\ge 1$, there is some $\rho=\rho(x,r,n)$ such 
that $f^n(B(x,r))\supset B(f^n(x),\rho)$. We then get a uniform bound $\rho(r)$ independent
of $x$ and $n$
by compactness of $X$ and the finiteness of iterates to be considered.

\begin{prop}\label{prop:diam_bounds} The map $f$ satisfies the axiom
[Diam] with $$\delta_\pm(r)= C^2\hde_\pm(C^2 r)\,.$$
\end{prop}

\pf Let $n\ge 0$, $k>0$, $\wtU\in\UUU_{n+k}$ and $U=f^k(\wtU)\in\UUU_n$.
By definition of metric finite type, one can find $m\in\MMM$ such that

$$\begin{array}{rcl} \wtU & \stackrel{f^k}{\longrightarrow} & U\\
\wt{\psi} \downarrow & & \downarrow \psi\\
\wtV & \stackrel{g}{\longrightarrow} & V\end{array}$$ where
we have dropped the indices.

\gap

Now suppose $\wtW\subset \wtU$ and $W=f^k(\wtW)$. Then since $\diam V = \diam \wtV = 1$ and $\psi$ and $\wt{\psi}$ distort ratios of diameters by at most a factor of $C^2$,  
$$\frac{\diam W}{\diam U}\le C^2 \diam\psi(W)\le C^2\hde_+(\diam \wt{\psi}(\wtW))\le C^2\hde_+\left(C^2\frac{\diam\wtW}{\diam\wtU}\right)\,.$$
Similarly, if $\wtW$ is a connected component of a connected set $W\subset U$, then

$$\frac{\diam\wtW}{\diam\wtU}\le C^2\diam\wt{\psi}(\wtW)\le C^2
\hde_-(\diam \psi(W))\le C^2\hde_-\left(C^2\frac{\diam W}{\diam U}\right)\,.$$
By definition of $\hde_\pm$, $\delta_\pm$ are homeomorphisms of $\R_+$.
\qed

\gap

{\noindent\bf Roundness bounds for model maps.}
We assume that $g:\wtV\to V$ is a model map, $W$ is a connected subset
of $V$, $\wtW$ is a component of $g^{-1}(W)$, and $\tx\in\wtW$.
We write $x=g(\tx)$.

\begin{lemma}\label{lemma:round_model} Under these notations, the following hold.
\be
\item If $\roundness(\wtW,\tx)\le K$ and $\diam \wtW\ge c>0$, then
$\ell_V(W,x)\ge \hrho(c/2K)\,.$

\item  If $\roundness(W,x)\le K$ and $\diam W\ge c>0$, then
$\ell_{\wtV}(\wtW,\tx)\ge \hde_+^{-1}(c/2K)\,.$
\eb\end{lemma}
Since model sets have diameter $1$, roundness bounds follow at once.
\gap

\pf 1. The definition of roundness and the hypothesis imply
$B\left(\tx,\frac{\diam\wtW}{2K}\right)\subset \wtW\,.$
Applying the model map $g$, we find 
$B\left(x,\hrho\left(\frac{\diam\wtW}{2K}\right)\right)\subset W\,.$
Thus, $\ell_V(W,x)\ge \hrho(c/2K)\,.$

2. Similarly,
$B\left(x,\frac{\diam W}{2K}\right)\subset W\,.$
Hence 
$B\left(\tx,\hde_+^{-1}\left(\frac{\diam W}{2K}\right)\right)\subset \wtW$
and so 
$$\ell_{\wtV}(\wtW,\tx)\ge \hde_+^{-1}((1/2K)\diam W)\ge \hde_+^{-1}(c/2K)\,.$$
\qed

\gap

\noindent{\bf Roundness distortion.} 
Recall that by the diameter bounds (\ref{eqn:diameter_bounds}), there are decreasing positive sequences $c_n, d_n$ such that 
\[0 < c_n \leq  \inf_{U \in \UUU_n} \diam U \leq \sup_{U \in \UUU_n} \diam U \leq  d_n.\]
Recall that the {\em Lebesgue number} of an open cover $\UUU$ of a compact metric space is the largest $\delta>0$ such that 
every subset of diameter $\delta$ is contained in an element of $\UUU$.  
 Let $\delta_0$ denote the Lebesgue
number of the covering $\UUU_0$. 
Since $X$ is perfect, each ball $B(x,r), x \in X, r>0$ has positive diameter.   
Fix $r_0<\delta_0/2$ and cover $X$ by balls of radius $r_0$.  
Extracting a finite subcover we conclude that for some positive constant $q_0$, 
every ball of radius $r_0$ has diameter at least $2q_0$.  
Set $K_0=d_0/q_0$.  It follows that for any $x\in X$, one can find $U_0(x)\in\UUU_0$ such that
$\roundness (U_0(x),x)\le K_0$.  Define $U_n(x) \in \UUU_n$ to
be the component of $f^{-n}(U_0(f^n(x)))$ which contains $x$.

Our first result says that by dropping down some uniform number of levels, every pair $(x,U)$, $x \in U, U \in \mathbf{U}$, is contained ``deep inside'' a larger set in $\mathbf{U}$.

\begin{prop}\label{prop:uniform_round}\be
\item There is some $n_0$ such that, for any $n,k\ge 0$,
if $x\in U\in\UUU_{n_0+n+k}$, then 
$U$ has compact closure in $U_n(x)$.
\item There is a constant $R\geq 1$ such that, for any $x\in X$
and any $n\ge 0$, $\roundness(U_n(x),x)\le R$.\eb\end{prop}

\pf We fix $n_0$ large enough so that, for $n\ge n_0$,
$$d_n\le c_0\min\left\{ \frac{1}{3K_0},\frac{\hde_-^{-1}(1/2C^2)}{C^2}\right\}\,.$$ 

\be
\item
It follows from $$d_n\le\frac{c_0}{3K_0}$$ that if $n\ge n_0$, $W\in\UUU_n$ and $x\in W$ then $\cl{W}\subset U_0(x)$.  To see this, notice that
on the one hand, for any $x\in X$,
$$\ell(U_0(x),x)\ge \frac{\diam U_0(x)}{2K_0}\ge  \frac{c_0}{2K_0} 
\hbox{ which implies } B\left(x, \frac{c_0}{2K_0}\right) \subset U_0(x), $$ and on the other hand
$$W\subset B\left(x,\frac{c_0}{3K_0}\right) \subset \overline{B(x, \frac{c_0}{3K_0})} \subset B(x, \frac{c_0}{2K_0}) \subset U_0(x)\,.$$
Therefore, for $n,k\ge 0$,
if $x\in U\in\UUU_{n_0+n+k}$, then $f^n(U)\in\UUU_{n_0+k}$ so that $\cl{f^n(U)}\subset U_0(f^n(x))$ . It follows
that $ \cl{U}\subset U_{n}(x)$ since $f^n:U_{n}(x)\to U_0(f^n(x))$ is proper.

\item
By construction, for each $x \in X$, we have $B(x, c_0/K_0) \subset U_0(x)$.  
The finite set of maps $\{f, f^{\circ 2}, \ldots, f^{\circ (2n_0-1)}\}$ 
is uniformly equicontinuous, so there exists $r_0>0$ such that for each 
$x \in X$ and each $n \in \{1, \ldots, 2n_0-1\}$, we have $f^n(B(x, r_0)) \subset B(f^n(x), c_0/K_0)$.  
It follows that $\roundness(U_n(x), x) \leq K_1$ where $K_1=r_0/c_{2n_0-1}$.

Let $0\le\ell\le  n_0-1$ and $j\ge 1$. We consider the model map given in the following diagram
$$\begin{array}{rcl} U_{jn_0}(x) & \stackrel{f^{jn_0}}{\longrightarrow} & U_0(f^{jn_0}(x))\\
\wt{\psi} \downarrow & & \downarrow \psi\\
\wtV & \stackrel{g}{\longrightarrow} & V\end{array}$$
It follows from the point above that $U_{n_0+\ell}(f^{jn_0}(x))\subset U_0(f^{jn_0}(x))$ and $U_{(j+1)n_0+\ell}(x)\subset U_{jn_0}(x)$.

Since the maps $\psi$ are uniform quasisimilarities which are onto, $$\roundness(\psi(U_{n_0+\ell}(f^{jn_0}(x))),\psi(f^{jn_0}(x)))\le C^2 K_1$$
and 
$$\diam \psi(U_{n_0+\ell}(f^{jn_0}(x)))\ge \frac{1}{C^2}\frac{\diam U_{n_0+\ell}(f^{jn_0}(x))}{\diam U_{0}(f^{jn_0}(x))}\ge\frac{1}{C^2}\frac{c_{2n_0}}{d_0}\,.$$
From Lemma \ref{lemma:round_model} and the fact that $\diam \wtV=1$, it follows that
$$\roundness(\wt{\psi}(U_{(j+1)n_0+\ell}(x)),\wt{\psi}(x))\le  K_2:=\frac{1}{\hde_+^{-1}\left(\dis\frac{c_{2n_0}}{2C^2 d_0 K_1}\right)}\,.$$

But by the definition of $\hde_-$ and $n_0$, $$\diam \wt{\psi}(U_{(j+1)n_0+\ell}(x))\le \hde_-\left(C^2\frac{d_{n_0}}{c_0}\right)\le \frac{1}{2C^2}\,,$$
so that Proposition \ref{prop:qsim_control} (1) (applied to $\wt{\psi}^{-1}$ and remembering $\diam V  = 1$) implies
$$\roundness(U_{(j+1)n_0+\ell}(x),x)\le \max\{C^2  K_2,\roundness(U_{jn_0}(x),x)\}\,.$$

But since $\roundness(U_{n_0}(x),x)\le K_1$, it follows by induction that for any $n\ge n_0$,
$$\roundness(U_n(x),x)\le \max\{ C^2 K_2, K_1\}\,.$$
Letting $R= \max\{ C^2 K_2, K_1,K_0\}$, it follows that, for any $x\in X$ and $n\ge 0$,
$$\roundness(U_n(x),x)\le R\,.$$\qed\eb

We may now deduce the roundness bounds for $n$ large enough:

\begin{prop}\label{prop:roundness}  Let $n_0$ and $R$ be the constants provided by Proposition \ref{prop:uniform_round}.  There is some $n_1\ge n_0$ with the following properties. Let  $n\ge n_1$, $k\ge 1$, $\wtU\in\UUU_{n+k}$,
$\tx\in\wtU$, $x=f^k(\tx)$ and  $U=f^k(\wtU)$. Let 
$$c=\frac{1}{C^2}\delta_+^{-1}\left(\frac{c_{n_1}}{d_{0}}\right)\,. $$
\be
\item If $\roundness(\wtU,\tx)\le K$ then
$$\roundness(U,x)\le\frac{1}{2R\hrho\left(\dis\frac{c}{2K}\right)}\,.$$
\item If $\roundness(U,x)\le K$ then
$$\roundness(\wtU,\tx)\le\frac{1}{2R\hde_+^{-1}\left(\dis\frac{c}{2K}\right)}\,.$$\eb
\end{prop}

\pf Choose $n_1 \ge n_0$ so that 
$$d_{n_1}\le\frac{c_0}{C^2}\hde_-^{-1}\left(\frac{1}{2C^2 R}\right)\,.$$

Set $m=n-n_1$ and let $U_{m+k}(\tx)$ and $U_m(x)$ be the neighborhoods provided by Proposition \ref{prop:uniform_round}(1), so that $\overline{U}\subset U_m(x)$ and $\overline{\wtU} \subset U_{m+k}(\tx)$.  Let us consider the model map $$\begin{array}{rcl} U_{m+k}(\tx) & \stackrel{f^{k}}{\longrightarrow} & U_m(x)\\
\wt{\psi} \downarrow & & \downarrow \psi\\
\wtV & \stackrel{g}{\longrightarrow} & V\end{array}$$

Note that $n_1$ is the difference in levels between $\wtU, U$ and their corresponding supersets  $U_{m+k}(\tx), U_m(x)$.  

\gap

Let us for the time being consider the dynchart 
$$\begin{array}{rcl} U \subset  U_{m}(x) & \stackrel{f^{m}}{\longrightarrow} & f^m U_m(x)=U_0 \supset f^m(U)\\
\psi \downarrow & & \downarrow \psi_0\\
V & \stackrel{g}{\longrightarrow} & V_0\end{array}$$
Since the quasisimilarity $\psi_0$ distorts ratios of diameters by at most the factor $C^2$ and $\diam V_0 = \diam V = 1$,
$$\diam \psi_0(f^m(U))\le C^2\frac{\diam f^m(U)}{\diam U_0}\le C^2\frac{d_{n_1}}{c_0}$$
and so, by the definition of $\hde_-$, we have 
\begin{equation}
\label{eqn:star}
\diam \psi(U)\le \hde_-\left(C^2\frac{d_{n_1}}{c_0}\right)\le \frac{1}{2C^2 R}
\end{equation}
where the last inequality follows by our choice of $n_1$.

By Proposition \ref{prop:diam_bounds},   
one also finds
$$\frac{\diam U}{\diam U_m(x)}\ge \delta_+^{-1}\left(\frac{c_{n_1}}{d_0}\right)\,.$$

The same argument applied to $(\wtU,U_{m+k}(\tx))$ yields

$$\diam \wt{\psi}(\wtU)\le \frac{1}{2C^2 R}$$
and 
$$\frac{\diam\wtU}{\diam U_{m+k}(\tx)}\ge \delta_+^{-1}\left(\frac{c_{n_1}}{d_0}\right)\,.$$

\be\item If  $\roundness(\wtU,\tx)\le K$ then since $\wt{\psi}$ is a surjective quasisimilarity,
$$\roundness(\wt{\psi}(\wtU),\wt{\psi}(\tx))\le C^2 K$$ 
and, by Lemma \ref{lemma:round_model} (with $K$ replaced with $C^2K$),
$$\ell_V(\psi(U),\psi(x))\ge \hrho\left(\frac{1}{2KC^2}\delta_+^{-1}\left(\frac{c_{n_1}}{d_0}\right)\right)\ge \hrho\left(\frac{c}{2 K}\right)\,.$$

But since  $$\diam \psi(U)\le \frac{1}{2C^2R}\,,$$
Proposition \ref{prop:qsim_control}, (2) and inequality(\ref{eqn:star}) imply that 
$$\roundness(U,x)\le C^2\frac{\diam \psi(U)}{\ell_V(\psi(U),\psi(x))}\le\frac{1}{2R\hrho\left(\dis\frac{c}{2 K}\right)}\,.$$

\item Similarly, if  $\roundness(U,x)\le K$ then
$\roundness(\psi(U),\psi(x))\le C^2 K$, and, by Lemma \ref{lemma:round_model},
$$\ell_{\wtV}(\wt{\psi}(\wtU),\wt{\psi}(\tx))\ge \hde_+^{-1}\left(\frac{1}{2KC^2}\delta_+^{-1}\left(\frac{c_{n_1}}{d_0}\right)\right)
\ge \hde_+^{-1}\left(\frac{c}{2 K}\right)\,$$

But since  $$\diam \wt{\psi}(\wtU)\le \frac{1}{2C^2R}\,,$$
it follows that 
$$\roundness(\wtU,\tx)\le \frac{1}{2R\hde_+^{-1}\left(\dis\frac{c}{2 K}\right)}\,.$$\qed\eb

{\noindent\bf Proof of Theorem \ref{thm:finite_type_implies_cxc}:}
Let us first note that a map of metric finite type is topological cxc if the axiom [Deg] holds. But this axiom 
follows from the fact that $\MMM$ is a finite set.

\gap

The axiom [Diam] is given by Proposition \ref{prop:diam_bounds}, and the roundness control holds 
as soon as sets of level at least $n_1$ are considered, by Proposition \ref{prop:roundness}.\qed

%\input sec2.tex

%\chapter{Examples}

\section{Examples of metric finite type systems}

\subsection{Expanding maps on manifolds}

If $X$ is metric space, and $f:X\to X$ is continuous, we say that $f$ is 
{\em expanding} if, for any $x\in X$, there is a neighborhood $U$ of $x$ such that, 
for any distinct $y,z\in U$, one has $|f(y)-f(z)|>|y-z|$; cf. \cite[\S\,1]{gromov:expanding}.
We first recall the statement of \cite[Thm 4.5.1]{kmp:ph:cxci} from
which we will deduce that expanding maps are metric finite type.

\begin{thm}[From expanding to homothety]
\label{thm:from_expanding_to_homothety}
Let $f: M \to M$ be an expanding map of a compact connected
Riemannian manifold to itself.
Then there exists a distance function  $d$ on $M$ and constants
$\delta>0$ and $\rho>1$ such that for all $x,y \in M$, 
\[ d(x,y)<\delta \implies d(f(x),f(y)) = \rho \cdot d(x,y)\]
and such that balls of radius $\leq \delta$ are connected and
contractible.    
\end{thm}

\begin{cor}[Expanding implies metric finite type]
\label{cor:from_expanding_to_homothety}
The dynamical system $((M,d),f)$ is metric finite type, hence cxc.  
\end{cor}

This refines Corollary 4.5.2 of \cite{kmp:ph:cxci}, which asserts that $((M,d), f)$ is merely cxc.  \\

{\noindent\bf Proof of Corollary \ref{cor:from_expanding_to_homothety}:} 
We remark that $f: M \to M$ is necessarily a
covering map of degree $D=\deg f$.  Let $\UUU_0$ be a finite open cover of
$M$ by open balls of radius $\delta$.   If $U \in
\UUU$ then since $U$ is contractible we have 
\[ f^{-n}(U) = \bigcup_{1}^{D^n}  \wtU_i\]
where the union is disjoint and where each $f^n | \wtU_i : \wtU_i \to U$ is a
homeomorphism which multiplies distances by exactly the factor $\rho^n$. 
In the definition of metric finite type, let $\VVV = \UUU_0$, and take the model maps all to be the identity maps.   Given $\wtU \in \UUU_n$ we let $\psi_\wtU: \wtU \to U=f^n(\wtU)$.  Then if $f^k: \wtU \to U$ we see that $g_{\wtU, U} = \id_U$ by construction.  Since each chart $\psi_U$ is a similarity, it follows that conditions (1)-(3) in the definition of metric finite type hold.   
Verification of axioms [Expans] and [Irred] are straightforward (details are in \cite{kmp:ph:cxci}).  
\qed

\subsection{Subhyperbolic rational maps} 

A rational map $f: \rs \to \rs$ is {\em subhyperbolic} if it has neither critical points in 
the Julia set with infinite forward orbit, nor parabolic cycles.  
Equivalently, under iteration, each critical point either converges to or lands in an attracting cycle, 
or lands in a repelling periodic cycle.  We begin by proving 

\begin{thm}[Subhyperbolic implies metric finite type]
\label{thm:subhyperbolic_implies_ft}
Let $f$ be a subhyperbolic rational map with Julia set $J$.  
Then there are closed neighborhoods $\XX_0, \XX_1$ of $J$ in the sphere such that $f: \XX_1 \to \XX_0$ 
is metric finite type with repellor $J$, with good open sets given by a finite collection $\UUU_0$ of open spherical balls.  
\end{thm}

Before proceeding, we collect a few facts and concepts from the theory of univalent functions.

\gap

\noindent{\bf Distortion principles.}  Let $\Delta_s = \{z \in \C : |z| < s\}$ and $\Delta = \Delta_1$. We will need 
\begin{lemma}
\label{lemma:koebe}
For all  $s>1$ and small $\rho>0$, there exists a constant 
$C(s,\rho)>1$ such that if 
\[ \psi: \Delta_s \to \rs \]
is any holomorphic embedding whose image omits a spherical 
disk of radius $\rho$, then $\psi|\Delta$ is a $(C(s,\rho), \lambda)$-
quasisimilarity (with respect to the Euclidean metric on $\Delta$ 
and the spherical metric on $\psi(\Delta)$), where $\lambda = |\psi'(0)|$.
\end{lemma}

\pf By composing with a spherical isometry  we may assume the image lies in the Euclidean disk about the origin 
of radius $R=R(\rho)$.  On such a disk the Euclidean and  spherical 
metrics are comparable, with constant depending only on $R$.  The lemma 
then follows from the usual Koebe principle.
\qed

\noindent{\bf Concentric disks.}  
We will also need the following concept.  Suppose 
\[ x \in U \subset W\]
where $U,W$ are conformally isomorphic to disks.  We say that $(W,U,x)$ 
is {\em concentric} if the triple 
$(W,U,x)$ is holomorphically isomorphic to the triple $(\Delta_s, \Delta, 0)$ 
for some $s>1$, i.e. there exists a conformal isomorphism sending $W \to \Delta_s$, $U \to \Delta$, and $x \mapsto 0$.  We denote by 
\[ f: (\wtW, \wtU, \tx) \to (W,U,x)\]
a proper holomorphic map $f: \wtW \to W$ with the property that $\wtU=f^{-1}(U)$ where $U$ and $\wtU$ are disks,  $x \in U$, $\tx \in \wtU$, $f$, if ramified, is branched only at $\tx$, and $f(\tx)=x$.
If $(W, U, x)$ is concentric and 
\[ f: (\wtW, \wtU, \tx) \to (W,U,x)\]
then $(\wtW, \wtU, \tx)$ is also concentric.  

\gap

{\noindent\bf Proof of Theorem \ref{thm:subhyperbolic_implies_ft}:}  
We begin by choosing carefully the covering $\UUU_0$.

Since $f$ is subhyperbolic, there exists a neighborhood $\XX_0$ of $J$ such that 
$\XX_0 \intersect P_f = \{x_i\}_{i=1}^p$ is a subset of $J$ and $\XX_1 = f^{-1}(\XX_0)$ is relatively compact in $\XX_0$.  
Choose $r>0$ sufficiently small such that $B(x,3r/2) \subset \XX_1$ for all $x \in J$, and 
$B(x_i, 3r/2) \intersect P_f = \{x_i\}$, $1 \leq i \leq p$.  Let $W_i = B(x_i, 3r/2)$ and $U_i=B(x_i,r)$.  
Since $J'=J-\union_i B(x_i, r)$ is compact, there exist finitely many points $x_i \in J'$, $p+1 \leq i \leq p'$ 
such that the collection $U_i = B(x_i, r/2), i=p+1, \ldots , p'$ covers $J'$. A simple comparison 
of the spherical and Euclidean metrics yields the existence of some radius $r'\in (r/2,(3/4)r)$ such that the triples 
$(B(0,r'),B(0,r/2),0)$ and $(B(0,(3/2)r),B(0,r),0)$ are isomorphic.
Then by construction, if $W_i = B(x_i, r')$, $i=p+1, \ldots, p'$, we have that the collection 
$\{(W_i, U_i, x_i)\}_{1\le i\le p'}$ consists of isomorphic concentric triples such that $\UUU_0 = \{U_i\}$ 
is a covering of $J$.  
Moreover, $P_f \intersect W_i =\{x_i\}$ for $1 \leq i \leq p$ and is otherwise empty.

We call $(W_i, U_i, x_i)$ a {\em triple at level zero} and drop the subscripts in what follows.  
A {\em triple at level $n$} consists of a component $\wtW$ of $f^{-n}(W)$, a component $\wtU$ of $f^{-n}$ contained in $\wtW$, 
and a preimage $\tx \in f^{-n}(x) \intersect \wtU$, where $(W, U, x)$ is a triple at level zero.   
By construction, triples $(\wtW, \wtU, \tx)$ at level $n$ are concentric,  
since $f^n|\wtW$ is possibly ramified only at $\tx$.  
We set $\UUU_n$ to be the elements $\wtU$ occuring in triples $(\wtW, \wtU, \tx)$ at level $n$.  
The {\em degree} $D$ of a triple $(\wtW, \wtU, \tx)$ at level $n$ is defined to be the degree of $f^n|\wtW$, 
and we let $\DDD \subset \N$ denote the 
set of such degrees.  Since $f$ is subhyperbolic, $\DDD$ is finite.

The triples at level zero are all conformally isomorphic 
and concentric.  Hence, there is $s>1$ such that  for any triple $(W, U, x)$ at level zero, there exists a  
Riemann map to a fixed triple 
\[ \varphi_W: (W, U, x) \to (\Delta_s, \Delta, 0)\]
which is unique up to postcomposition by a rotation about the origin.  For each 
triple at level zero, we make such a choice arbitrarily.   Given $(\wtW, \wtU, \tx)$ a triple at level $n$ of degree $D$, let 
\[ \psi_{\wtW}: (\wtW, \wtU, \tx) \to (\Delta_{s^{1/D}}, \Delta, 0)\]
be the conformal isomorphism given by 
\[ \psi_{\wtW} = \left( \varphi_W \circ f^n|\wtW \right) ^{1/D}\]
where the principal branch of root is used, and set $\psi_{\wtU} = \psi_{\wtW}|\wtU$; 
these will be the maps as in the definition of metric finite type.  

We now verify the conditions in the definition.  If 
\[ f^k: (\wtW, \wtU, \tx) \to (W, U, x)\]
sends a triple at level $n+k$ of degree $\tilde{D}$ to one at level $n$ of degree $D$, 
then 
\[ \wtpsi^{-1}_{\wtU}  \circ f^k \circ \psi_U : (\Delta, 0) \to (\Delta, 0) \]
is just $z \mapsto z^m$, where $m=\tilde{D}/D$.  Thus, the set of 
model maps is just the set of maps of the unit disk to itself given by $g_m(z) = z^m$, 
where $m$ ranges over the set $\MMM$ of all local degrees of iterates of $f$ 
at points in the Julia set. Since $\DDD$ is finite, $\MMM$ is finite as well.  

The spherical diameters of the sets $\wtW$ 
arising in triples at level $n$ tend to zero uniformly, in fact, exponentially, in $n$ (cf. \cite{steinmetz:book:iteration}, 
Lemma 5.1.4.)  
Thus, for any triple $(W, U, x)$ at any level, the set $W$ omits a disk of some 
definite spherical radius $\rho$.  Lemma \ref{lemma:koebe} and the 
finiteness of the set of radii $s^{1/D}$, $D \in \DDD$, then implies 
that 
\[ \left\{ \psi_U : \Delta \to \rs, \;\; U \in \bigcup_n \UUU_n \right\}\]
is a family of uniform quasisimilarities, and the proof is complete.
\qed

We now provide a converse statement.  

\begin{thm}
Suppose $f: \IS^2 \to \IS^2$ is orientation-preserving and metric finite type with respect to the standard spherical 
metric.  Then $f$ is quasisymmetrically, hence quasiconformally equivalent 
to a postcritically finite rational map whose Julia set is the whole sphere.   
\end{thm}

\pf Theorem \ref{thm:finite_type_implies_cxc} implies that $f$ is cxc with respect to the standard 
spherical metric, and Theorem 4.2.7 in \cite{kmp:ph:cxci} implies that $f$ is qs, hence qc 
conjugate to a semihyperbolic rational map.  Theorem \ref{thm:when_X_is_a_surface} (see below) implies that $P_f$ is finite.
\qed

\subsection{Topologically finite type}
\label{secn:top_finite_type}

The constructions in the proof of Theorem \ref{topcxc_cxc1},  combined with those of the proof of Theorem \ref{thm:subhyperbolic_implies_ft}, yield a metrization result for a certain class of topological dynamical systems, which we now define precisely.  Since we wish the gauge of the constructed metric to depend only on the topological dynamics, we assume $\XX_0 = \XX_1 = X$.  
\gap

\noindent{\bf Topologically finite type dynamics.}  Let $f: X \to X$ be a finite branched covering with repellor $X$ as in \S \ref{secn:def_top_cxc}.  
Let $\UUU_0$ be an open cover  of $X$ by connected subsets of $X$.  An {\em iterate over $\UUU_0$} is a map of the form $f^k: \wtU \to U$ where $\wtU \in \UUU_{n+k}$ and $U \in \UUU_n$ for some $n$ and $k$.   

\begin{defn}
\label{defn:orbit_equivalent}
Elements  $U_1, U_2 \in \mathbf{U}$ are said to be {\em orbit isomorphic} provided there exists $U \in \mathbf{U}$, $k_1, k_2 \in \N$, and a homeomorphism $\psi: U_1 \to U_2$ such that $f^{k_i}: U_i \to U, i=1,2$ and $f^{ k_2} \circ \psi = f^{k_1}$.  
The map $\psi$ is called an {\em orbit isomorphism}.  
\end{defn}

For example, if $f^{\ell}: U_1 \to U_2$ is a homeomorphism, then $U_1$ and $U_2$ are orbit isomorphic.  Similarly,  if $f^{k_i}: U_i \to U, i=1,2$ are homeomorphisms, then $U_1$ and $U_2$ are orbit isomorphic.  
Recall that $|U|=n$ means that $U\in\UUU_n$.
It is easily verified that in the definition of orbit isomorphism, one may replace $U$ with $f^{|U|}(U) \in \UUU_0$ i.e., one may assume that $U \in \UUU_0$.  It follows easily that the relation of being orbit isomorphic is an equivalence relation.  By definition, distinct elements of $\UUU_0$ are never orbit isomorphic.  

The proof of the following lemma is straightforward.

\begin{lemma}
\label{lemma:orbit_eq_preserves_sets}
If $\psi: U_1 \to U_2$ is an orbit isomorphism and $U_1' \in \mathbf{U}$ is a subset of $U_1$, then $U_2' = \psi(U_1') \in \mathbf{U}$ and $|U_2'| - |U_2| = |U_1'|-|U_1|$.
\end{lemma}

\begin{defn}[Orbit isomorphism of iterates]
Two iterates $f^{k_1}: \wtU_1 \to U_1$ and $f^{k_2}:\wtU_2 \to U_2$ over $\UUU_0$ are {\em orbit isomorphic} if there exist orbit isomorphisms $\wtpsi: \wtU_1 \to \wtU_2$ and $\psi: U_1 \to U_2$ such that the diagram 
\[
\begin{array}{rcl}
\wtU_1 & \stackrel{\wtpsi}{\longrightarrow} & \wtU_2 \\
f^{k_1} \downarrow & \; & \downarrow f^{k_2}\\
U_1 & \stackrel{\psi}{\longrightarrow} & U_2
\end{array}
\]
commutes.
\end{defn}

\begin{defn}[Topologically finite type]
The map $f$ is said to be {\em topologically finite type with respect to $\UUU_0$} 
if it satisfies Axioms [Expans] and [Irred], and if in addition  there are only finitely many orbit isomorphism types of iterates over $\UUU_0$.   
\end{defn}

For example, if $f$ is a hyperbolic rational map with connected Julia set $X$ and $\UUU_0$ consists of open sets not separating $P_f$, then $f$ has exactly $\#\UUU_0$ isomorphism types of iterates and so is topologically finite type with respect to $\UUU_0$.    If $f$ is merely subhyperbolic, then it is topologically finite type with respect to the open covering $\UUU_0$ defined in the proof of Theorem \ref{thm:subhyperbolic_implies_ft}.

The following is the main result of this section.

\begin{thm}[Topological implies metric finite type]
\label{thm:top_implies_metric_ft}
If $f: X \to X$ is topologically finite type with respect to $\UUU_0$, 
then when $X$ is equipped with metric $d_\varepsilon$ constructed in \S \ref{subsecn:metric_cxc}, 
$f$ is metric finite type,  and the constant $C$ in the definition of metric finite type can be taken to be $1$.  
\end{thm}

\begin{cor}
\label{cor:literally_selfsimilar}
Let $f: S^2 \to S^2$ be a postcritically finite Thurston map 
which is topologically cxc with respect to an open covering $\UUU_0$. 
Then in the metric $d_\varepsilon$, $f$ is metric finite type, and up to {\em similarity}, 
there are only finitely many possibilities for the sets $U \in \mathbf{U}$.  
\end{cor}

In \S \ref{secn:fsr} we will discuss the implications of this corollary for finite subdivision rules.

The proof will rest on the following observation.  Recall that $|x-y|_\varepsilon$ denotes the distance between $x$ and $y$ in the metric $d_\varepsilon$.

\begin{prop}
\label{prop:orbit_equivalence_is_similarity}
Suppose $\psi: U_1 \to U_2$ is an orbit isomorphism and $B_\varepsilon(u, r) \subset U_1$. 
Let $\lambda = \exp(\varepsilon(|U_1|-|U_2|))$.  
Then for all $x, y \in B_\varepsilon(u, r)$, 
\[ |\psi(x)-\psi(y)|_\varepsilon = \lambda|x-y|_\varepsilon.\]
In particular, 
\[ \psi(B_\varepsilon(u, r)) = B_\varepsilon(\psi(u), \lambda r)\]
and $\psi|_{B_\varepsilon(u,r)}$ is a $(1,\lambda)$-quasisimilarity. 
\end{prop}

\pf Let $x_1=x, y_1=y, x_2=\psi(x_1), y_2=\psi(y_1)$.  
The points $x_1, y_1$ lie in $\cl{\Gamma}_\varepsilon$, which is a complete length space.   So there is a geodesic 
$\gamma_1: [0, L] \to \cl{\Gamma}_\varepsilon$ joining $x_1$ and $y_1$.  By reparameterizing $\gamma_1$ as $\gamma_1: \R \to \cl{\Gamma}$, 
we may assume $\gamma_1(\Z) \subset \mathbf{U}$ (the vertex set of $\Gamma$), that $\lim_{t \to -\infty}=x_1$, and that $\lim_{t \to +\infty}\gamma_1(t)=y_1$.
By \cite[Lemma 3.3.4]{kmp:ph:cxci}, 
\[ \cl{\bigcup_{n \in \Z}\gamma(n)} \subset B_\varepsilon(u, r) \subset U_1.\]
Since $\psi$ is an orbit isomorphism, $\psi(\gamma_1(n)) \in \mathbf{U}$ for all 
$n \in \Z$.  Moreover, $\gamma_1(n) \intersect \gamma_1(n+1) \neq \emptyset \implies \psi(\gamma_1(n)) \intersect \psi(\gamma_1(n+1) ) \neq \emptyset$.    
Lemma \ref{lemma:orbit_eq_preserves_sets} implies that there exists $\gamma_2: \Z \to \Gamma$ given by $\gamma_2(n)=\psi(\gamma_1(n))$ and such that 
$|\gamma_1(n)|-|\gamma_2(n)| = |U_1|-|U_2|$.    Moreover, $\gamma_2$ extends to a curve $\gamma_2: \R \to \Gamma$ such that for all $n \in \Z$, 
$\gamma_2|_{[n, n+1]}$ traverses a closed edge exactly once.  
Thus 
\[ |x_2-y_2|_\varepsilon \leq \int_{\gamma_2}\varrho_\varepsilon \ ds = \lambda \int_{\gamma_1}\varrho_\varepsilon\ ds =  \lambda |x_1-y_1|_\varepsilon.\]

Hence $B_\varepsilon(\psi(u), \lambda r) \subset U_2$.  By considering $\psi^{-1}$ and applying the same argument we conclude 
\[ |x_1-y_1|_\varepsilon \leq \lambda^{-1}|x_2-y_2|_\varepsilon.\]
Hence 
\[ \lambda|x_1-y_1|_\varepsilon \leq \lambda \cdot \lambda^{-1}|x_2-y_2|_\varepsilon \leq \lambda |x_1-y_1|_\varepsilon\]
so equality holds throughout. 
\qed

{\noindent\bf Proof of Theorem \ref{thm:top_implies_metric_ft}:}
By \cite{kmp:ph:cxci}, Prop. 3.3.2(1), there exists an integer $l>0$ with the following property.  
For any $U \in \mathbf{U}$ with $|U|$ sufficiently large, there exists a ball $B$ and 
$\hatU \in \mathbf{U}$ such that $U \subset B \subset \hatU$ and $|U|-|\hatU| = l$.   
The set $\hatU$ will play a role in this proof similar to that played by the set denoted $W$ in the proof of Theorem \ref{thm:subhyperbolic_implies_ft}:  
it will provide some "Koebe space"; the control of distortion will be provided here by Proposition \ref{prop:orbit_equivalence_is_similarity}.

In what follows, we assume that for each $U \in \mathbf{U}$, a choice of such larger $\hatU \supset U$ has been made;  
we refer to the couple $(\hatU, U)$ as a {\em pair}.  An {\em orbit isomorphism} of pairs is an orbit isomorphism 
$\psi: \hatU_1 \to \hatU_2$ such that the restriction $\psi|_{U_1}: U_1 \to U_2$ is also an orbit isomorphism; in this case  
$\psi: (\hatU_1, U_1) \to (\hatU_2, U_2)$ is a map of pairs.  An {\em iterate of pairs} is a map of pairs $f^k: (\wthatU, \wtU) \to (\hatU, U)$.  

Two iterates of pairs $f^{k_1}: (\wthatU_1, \wtU_1) \to (\hatU_1, U_1)$ and $f^{k_2}: (\wthatU_2, \wtU_2) \to (\hatU_2, U_2)$ are 
{\em orbit isomorphic} provided there are orbit isomorphisms $\wtpsi: (\wthatU_1, \wtU_1) \to (\wthatU_2, \wtU_2)$ and 
$\psi: (\hatU_1, U_1) \to (\hatU_2, U_2)$ such that the diagram 
\[
\begin{array}{rcl}
(\wthatU_1, \wtU_1) & \stackrel{\wtpsi}{\longrightarrow} & (\wthatU_2, \wtU_2) \\
f^{k_1} \downarrow & \; & \downarrow f^{k_2}\\
(\hatU_1, U_1) & \stackrel{\psi}{\longrightarrow} & (\hatU_2, U_2)
\end{array}
\]
commutes.

Now suppose $f: X \to X$ is topologically finite type with respect to $\UUU_0$.  
It follows immediately that Axiom [Deg] holds, and consequently, that the graph $\Gamma$ defined with respect to the open covering $\UUU_0$ 
is uniformly locally finite.  It follows that in $\Gamma$ defined with the standard graph metric which makes each edge isometric to $[0,1]$, 
the ball of radius $l$ about a given vertex corresponding to a set $\hatU$ contains at most finitely many, say $T$, vertices corresponding to sets $U$.  
Since the number of isomorphism classes of iterates $f^k: \wthatU \to \hatU$ is finite, the number of isomorphism classes of iterates of pairs 
$f^k: (\wthatU, \wtU) \to (\hatU, U)$ is at most $T$ times this number, hence is also finite.

Let $\MMM$ be an index set enumerating the isomorphism classes of iterates of pairs, so that each iterate of pairs is orbit isomorphic via isomorphisms 
$\wtpsi, \psi$ to a map of the form 
\[ g_m: (\wthatV_m, \wtV_m) \to (\hatV_m, V_m), \;\;\; m \in \MMM.\]
By Proposition \ref{prop:orbit_equivalence_is_similarity} and the construction of the neighborhoods $\hatU$, in the metric $d_\varepsilon$, the maps 
$\wtpsi, \psi$ are similarities when restricted to $\wtU$, $U$, respectively.  Hence, in the metric $d_\varepsilon$, 
the set of maps $\{g_m|_{\wtV}: \wtV \to V\}_{m \in \MMM}$ form a dynatlas, and so $f: X \to X$ is metric finite type.  
\qed

{\noindent\bf Proof of Corollary \ref{cor:literally_selfsimilar}:}  Suppose $f: S^2 \to S^2$ is a topologically cxc Thurston map.  
Axiom [Deg] implies in particular that $f$ does not have periodic critical points.  Choose $\UUU_0$ to be any collection of open disks $U$ such that, 
if the closure of $U$ meets $P_f$, then it does so in at most one point, and this point lies in $U$.  It follows easily as 
in the proof of Theorem \ref{thm:subhyperbolic_implies_ft} that there are only finitely many orbit isomorphism types of iterates over $\UUU_0$.  
The Corollary then follows by Theorem \ref{thm:top_implies_metric_ft}.
\qed

\subsection{Maps of finite type on a surface}

When $\XX_0=\XX_1=X$ is a compact surface without boundary, the possibilities for a finite type map are greatly restricted.

\begin{thm}  
\label{thm:when_X_is_a_surface}
Suppose $f: X \to X$ is topologically finite type, where $X$ is a surface without boundary.  Then either $X$ is a torus or Klein bottle and $f$ is 
an unramified covering map, or $X=S^2$ or $\R\IP^2$ and the postcritical set $P_f$ is finite.  
\end{thm}

For a surface with boundary, the branch set of a topologically finite type map can be infinite:  if $f: [-2, 2] \to [-2, 2]$ is the map $f(x)=x^2-2$, 
then $f\times f$ is topologically finite type on 
$[-2,2] \times [-2, 2]$; as $\UUU_0$ one may take slightly thickened neighborhoods of the four corner squares of side length two.  
\gap

\pf The Riemann-Hurwitz formula implies that the possibilities for $X$ are those given in the statement, and that $f$ is unramified in the genus one case.   

Let $\MMM$ be a finite set indexing representatives of orbit isomorphism classes.    
Thus, $\MMM=\{g_m: \wtV_m \to V_m\}$ is a finite set of fbc's with the property that, given any $U \in \UUU_0$ and 
$f^k: \wtU \to U$, there exist a model map $g_m: \wtV_m \to V_m$ and maps $\psi_\wtU: \wtU \to \wtV_m$ and 
$\psi_U: U \to V_m$ such that 
$\psi_U \circ f^k \circ \psi_{\wtU}^{-1}=g_m$. 
The map $f^k: X \to X$ has finitely many branch values, so there are only finitely many branch values of the model map $g_m$.   
Since the set $\MMM$ of model maps is finite, for fixed $U$, as $\wtU$ and $k$ vary, there are only finitely many branch values of 
$f^k$ in $U$.  Hence $U \intersect P_f$ is finite.  Since $\UUU_0$ is finite, $P_f$ is finite.
\qed
%

%The model map 
%\[ g_{\wtU, U} = \psi_U \circ f^k \circ \psi_{\wtU}^{-1}: \wtV \to V\]
%is an fbc.   The map $f^k: X \to X$ has finitely many branch values.  Hence, for fixed $\wtU$, there are only finitely many branch values of $f^k: \wtU \to U$ in $U$, and so there are only finitely many branch  values of $g_{\wtU, U}$ in $V$.  The definition of topologically finite type implies that for fixed $U$, as $\wtU$ varies, the set of maps $g_{\wtU, U}$ arising as above is finite.  Hence 
%\[ \# \{ v \in V: \exists m \in \MMM \mbox{ such that } V = V_m, v \in B_{g_m}\}< \infty.\]
%Hence, for fixed $U$ and variable $f^k: \wtU \to U$, there are only finitely many possibilities for the location in $U$ of a critical value of $f^k$.  This implies that $U \intersect P_f$ is finite.  Since $\UUU_0$ is finite, $P_f$ is finite.  
%\qed

\subsection{Finite subdivision rules}
\label{secn:fsr}

In this subsection, we show that another natural source of examples of finite type dynamics comes from the {\em finite subdivision rules} considered by Cannon, Floyd, and Parry.  We first briefly summarize their definition, focusing on the case when the underlying dynamics takes place on the two-sphere; cf. \cite{cfp:fsr}.  
\gap

\noindent{\bf Finite subdivision rules on the two-sphere.}  A finite subdivision rule (f. s. r.) $\RRR$ consists of a finite 2-dimensional CW complex  $S_\RRR$, a subdivision $\RRR(S_\RRR)$ of $S_\RRR$,  and a continuous cellular map $\phi_\RRR: \RRR(S_\RRR) \to S_\RRR$ whose restriction to each open cell is a homeomorphism.   When the underlying space of $S_\RRR$ is homeomorphic to the two-sphere $S^2$ (for concreteness, we consider only this case) and $\phi_\RRR$ is orientation-preserving, $\phi_\RRR$ is a postcritically finite branched covering of the sphere with the property that pulling back the tiles effects a recursive subdivision of the sphere; below, we denote such a map by $f$.   That is, for each $n \in \N$, there is a subdivision $\RRR^n(S_\RRR)$ of the sphere such that $f$ is a cellular map from the $n$th to the $(n-1)$st subdivisions.   Thus, we may speak of {\em tiles} (which are closed 2-cells), {\em faces} (which are the interiors of tiles), {\em edges}, {\em vertices}, etc. at {\em level } $n$.  It is important to note that formally, an f. s. r. is {\em not} a combinatorial object, since the map $f$, which is part of the data, is assumed given.   In other words:  as a dynamical system on the sphere, the topological conjugacy class of $f$ is well-defined.  
A subdivision rule $\RRR$ has {\em mesh going to zero} if for every open cover of $S_\RRR$, there is some integer $n$ for which each tile at level $n$ is contained in an element of the cover.  It has {\em bounded valence} if there is a uniform upper bound on the valence of any vertex at any level.    In this case, $f$ is a Thurston map without periodic critical points.  

In \cite{kmp:ph:cxci} it is shown that if $\RRR$ is a bounded valence finite subdivision rule on the 
sphere with mesh going to zero, then there are integers $n_0, n_1$ with the following property.  
Let $t$ be a tile of $\RRR^{n_0}(S^2)$, let $D_t$ be the star of $t$ in $\RRR^{n_0+n_1}(S^2)$ 
(that is, a ``one-tile neighborhood'' of $t$), and let $U_t$ be the interior of $D_t$.  
Let $\UUU_0$ be the finite open covering of $S^2$ defined by sets of the form $U_t$.  
Then each $U$ belonging to $\UUU_0$ and all of its iterated preimages are Jordan domains, and with respect to $\UUU_0$, 
the dynamical system $f: S^2 \to S^2$ is topologically cxc.  

The fact that $\RRR$ has bounded valence implies that up to cellular isomorphism, as $t$ varies through tiles at all levels, there are only finitely many possibilities for the cell structure of $D_t$ in $\RRR^{n_0+n_1}(S^2)$.   Consequently, if $t$ is a tile at level $n+k$, up to pre- and post-composition by cellular isomorphisms in domain and range, there are only finitely many possibilities for the cellular map $f^k: D_t \to f^k(D_t)$.   This implies that with respect to $\UUU_0$, such a map $f$ is of topologically finite type.  We conclude 

\begin{cor}
\label{cor:fsr_are_ft}
Suppose $\RRR$ is a bounded valence finite subdivision rule on the two-sphere with mesh going to zero.  Let $\UUU_0$ be the level zero good open sets as constructed in \cite[\S 4.3]{kmp:ph:cxci}.  
Then in the natural metric $d_\varepsilon$, the subdivision map $f: S^2 \to S^2$ is of metric finite type with constant $C=1$.  
\end{cor}

Compare \cite[Lemma 4.3]{cfp:xcI} and \cite[\S\,8]{bonk:meyer}.  

The preceding corollary implies that, up to similarity, there are only finitely many possible shapes of tiles.  We now explain this precisely.

Recall that by construction, each $U \in \UUU_0$ is a union of tiles at level $N_0 = n_0 + n_1$.  Choose a representative collection  of iterates comprising a dynatlas $(\VVV, \MMM)$.  Suppose $V_m \in \VVV$  is a model open set.  By construction, there is some $k_m$ such that $f^{k_m}: V_m \to U \in \UUU_0$.   Then $\cl{V}_m$ is a union of finitely many cells at level $N_0+k_m$.  Since $\MMM$ is finite, there exists $M \in \N$ such that each model open set $V_m$ is a union of finitely many tiles $s$ at a uniform level $M$ independent of $m$.  Then the set of such tiles $s$ arising in this way is finite.  

Now suppose $U \in \mathbf{U}$ is arbitrary and $|U|=n$.  By construction, there is a model open set $V_m \in \VVV$ and a similarity $\psi_U: U \to V_m$.  Let us say that a {\em tile at level $n+M$} is a set of the form $\psi_U^{-1}(s)$ where $s \subset V_m$ is a tile as in the previous paragraph.    We conclude that up to similarity, there are only finitely many tiles.

\gap

%\input sec3.tex

%\chapter{Summary of results due to V. Nekrashevych}

\section{Selfsimilar groups} 

In this chapter, we summarize some results of V. Nekrashevych from \cite{nekrashevych:book:selfsimilar}.    

\subsection{Group actions on rooted trees}
\label{secn:ssg}

Let $X$ be an alphabet consisting of $d \geq 2$ symbols.  For $n \geq 1$
denote by
$X^n$ the set of words of length $n$ in the alphabet $X$ and let $X^0 =
\{\emptyset\}$ consist of the empty word.   
Let $X^* = \union_n X^n$ and $X^{-\omega} = \{...x_3x_2x_1 : x_i \in X\}$.  The
length, or level, of a word $w$ will be denoted $|w|$.
For $x\in X$ and $w\in X^*$, the left shift map $\tau: X^* \to X^*$ is defined by $\tau(xw)=w$
and  the right shift map $\sigma: X^* \to X^*$ is defined by $\sigma(wx)=w$.  
Obviously, these two shifts commute: $\tau\circ \sigma = \sigma \circ \tau$. 

Let $G$ be a finitely generated group acting faithfully on $X^*$ in a
manner which preserves the lengths of words and which is transitive on
each $X^n$.   We write the action as a right action, so that $w^g$ is the
image of $w \in X^*$ under the action of $g$. The action is called {\em
selfsimilar} if for each $g \in G$ and $x \in X$ there exists $h \in G$
such that for all $w \in X^*$, 
\[ (xw)^g = x^g w^h.\] The element $h=h(x,g)$ is uniquely determined, is
called the {\em restriction} of $g$ to $x$, and is denoted $g|_x$.  More
generally, given $u \in X^*$ and $g \in G$ one finds that for all $v \in
X^*$, 
\[ (uv)^g = u^g v^{g|_u}\] for a uniquely determined element $g|_u$ called
the restriction of $g$ to $u$.    One finds readily the identities 
\[ g|_{uv} = (g|_u)|_v\] for all $u, v \in X^*$ and all $g \in G$, and
(remembering that the action is on the right)
\[ (gh)|_v = (g|_v)(h|_{v^g}).\]

\noindent{\bf Example.}  Let $T(X)$ denote the infinite rooted tree defined as follows.  The 
vertex set
is $X^*$.  For all $w \in X^*$ and $x \in X$, an edge joins $w$ and $wx$.  The root $o$ is the empty word $\emptyset$.  Let $G = \Aut(T(X))$ denote the
group of automorphisms of $T(X)$.  Note that  $G$ acts transitively on $X^n$ for each $n$.   For $x \in
X$ let $T_x$ denote the rooted subtree spanned by the set of vertices $xX^*$ with the root $x$.  Then the map 
\[ l_x: T(X) \to T_x, \;\; w \mapsto xw\] is an isomorphism of rooted
trees.   Given
$g \in G$ and $x \in X$, one has 
\[ g|_x = l_{x^g}^{-1} \circ g|_{T_x} \circ l_x\] on vertices.  Indeed,
the map 
\[ \rho: G \to G^d \rtimes \Sym(X), \;\; \rho(g) = ((g|_{x_1}, \ldots,
g|_{x_d}),
\sigma(g))\] is an isomorphism, where $\Sym(X)$ is the symmetric group on
$X$ and
$\sigma(g)$ is the action of $g$ on $X=X^1$.  In particular, the action of $G$ on $X^*$ is selfsimilar.  
\gap

\subsection{Contracting actions}  
In principle, the selfsimilarity of the action means that the image of a
word under the action of a group element can be recursively computed.  However, there exists the  danger that the word lengths of restrictions may blow up as the recursive
algorithm progresses.  The following definition is designed to capture a
robust condition which ensures that this does not occur. 

The action of $G$
on $X^*$ is called {\em contracting} if there exists a finite subset
$\KKK$ of $G$ such that for each $g \in G$, there is a ``magic level''
$m(g)$ with the following property:  for all words $v$ with $|v| \geq m(g)$,
the restrictions
$g|_v$ lie in $\KKK$.   The smallest such subset $\KKK$ is called the {\em
nucleus} of the action, denoted $\NNN$. It follows easily that there is always a ``good
generating set'' $S$ which contains the nucleus $\NNN$ and which is closed
under restrictions.  The proof is the greedy algorithm.  Start with an
arbitrary generating set $S$.  Replace $S$ with $S\union \{s|_x : x \in X, s \in S\}$, and repeat.  The process stabilizes, by the identities of restrictions given
above and the definition of contracting. 

The action is called {\em recurrent} if it acts transitively on the first level $X^1$ and for some (equivalently, for any) $x \in X$, the homomorphism $\Stab_G(x) \to G$ given by $g \mapsto g|_x$ is surjective.  If the action is recurrent, then it is transitive at each level.  

\subsection{Selfsimilarity complexes}

Suppose $G$ is a finitely generated group equipped with a faithful, selfsimilar, level-preserving, transitive-on-each-level, contracting group action on $X^*$ as in the previous section, and let $S$ be a
finite generating set for
$G$. Denote by $||g||$ the minimum length of a word in the generators representing $g$. 

Associated to this data is an infinite cellular 1-complex
$\Sigma = \Sigma(G,S)$ with labelled vertices, oriented labelled edges, and a basepoint, defined as
follows.  The 0-cells are the set of words $X^*$, labelled by the word. 
The 1-cells come in two types:  {\em vertical} edges, running ``up'' from $w$ to $xw$
and labelled $x$ (yes, this is backwards from the construction of the tree
$T(X)$)  and {\em horizontal} edges, running ``over'' from $w$ to
$w^s$ and labelled $s$, where $s \in S$.   The basepoint $o$ is the vertex corresponding to the empty word.   We denote by $\Sigma'$ the subcomplex spanned by words of length at least $1$, so that $\Sigma'$ is obtained by removing the root $o$ and all edges incident to $o$.

We think of the root $o$ at the bottom.  
The subcomplex spanned by the vertical edges is isomorphic to the tree $T(X)$, which branches upwards.
For each $n \in \N$, the subcomplex spanned by $X^n$ is the so-called {\em Schreier graph} of the action of $G$ on $X^n$ with respect to the generating set $S$.  Within such a Schreier graph, there may be multiple edges joining a pair of vertices,  and edges may join a vertex to itself.  The hypothesis that the action is transitive on each level implies that for each level $n$, this subgraph is connected.

We endow $\Sigma$ with the unique length metric so that 1-cells are locally isometric to the unit interval; 
the distance between vertices is thus the usual combinatorial distance.  
We denote the distance between points $p,q \in \Sigma$ by $|p-q|$.  
With this additional structure we may speak of {\em geodesics}, i.e. 
paths $\gamma: [0,L] \to \Sigma$ satisfying $|\gamma(0)-\gamma(t)| = t$ for all $0 \leq t \leq L$, 
where $L=+\infty$ is permitted.  The space $\Sigma$ is {\em geodesic}: any two points of $\Sigma$ 
are joined by at least one, and possibly several, geodesics.  
It is also {\em proper}: closed balls (of finite radius) are compact.   

The definitions easily imply that given any vertex $w=x_n\ldots x_2x_1 \in \Sigma$, there is exactly one geodesic joining $o$ to $w$, namely, the one whose image traverses the sequence of labelled vertical edges starting with $x_1$, then $x_2$, etc. and ending in $x_n$.  
Similarly, an infinite geodesic starting at the root $o$, which we call a {\em ray} for short, necessarily consists entirely of vertical geodesic segments.  The set of rays is thus identified with the set $X^{-\omega}$.  
Two rays (equivalently, two infinite sequences) $\ldots x_3x_2x_1$ and $\ldots y_3y_2y_1$ are {\em asymptotically equivalent} if there exists some constant $C>0$ such that each is contained in some $C$-neighborhood of the other.   Equivalently: there exists a sequence $g_1, g_2, g_3, \ldots$ of elements of $G$, ranging over only a finite set, such that for each $n \in \N$, $(x_n \ldots x_1)^{g_n} = y_n \ldots y_1$. 
The following theorem summarizes the principal features of these objects (see Definition \,\ref{def:bdrySigma} for $\bdry\Sigma$).

\gap

\begin{thm}
\label{thm:props_of_Sigma}
\be
\item The quasi-isometry type of $\Sigma$ is independent of the chosen
generating set $S$.
\item The space $\Sigma$ is Gromov hyperbolic.
\item The $\sim_{asymp}$-equivalence classes are finite, of cardinality at
most that of the nucleus, $\mathcal{N}$.  The Gromov boundary at infinity
$\bdry \Sigma$ is homeomorphic to the quotient $X^{-\omega}/\sim_{asymp}$.
It is compact, connected, metrizable, has topological dimension at most
$\#\NNN-1$, and, if
the action is recurrent, is locally connected.
\item The right shift $\sigma$ preserves the asymptotic equivalence
relation and defines a natural cellular $d$-to-one covering map $F:
\Sigma' \to \Sigma$ whose boundary values $\bdry F: \bdry \Sigma \to \bdry
\Sigma$ coincide with the map on the quotient space
$X^{-\omega}/\sim_{asymp}$ induced by $\sigma$.  The map $\bdry F$ is
continuous, surjective, and at most finite-to-$1$.
\item If the action is recurrent, then $\bdry F$ is a
branched covering map of degree $d$.
\eb
\end{thm}

The boundary value map $\bdry F: \bdry \Sigma \to \bdry \Sigma$ is defined as follows.  Given $\xi \in \bdry \Sigma$, to compute $\bdry F(\xi)$, choose
a representing geodesic ray $R$;
then $\bdry F(\xi)$ is represented by the geodesic ray $F(R\intersect \Sigma')$.  
\gap

\pf (1), (2), and (3)  are implied by Lemma 3.8.4, Theorem 3.8.6, and
Theorem 3.6.3 of
\cite{nekrashevych:book:selfsimilar}, respectively; see also
\cite{nekrashevych:hyperbolic}.   (4) is a general consequence of
hyperbolicity and the fact that $F$ is cellular.
It is clear that $\bdry F$ is continuous and surjective.
It is finite-to-one since a point of the boundary is represented by at
most finitely many rays,
and every ray has finitely many preimages. 
(5) seems to be
implied by the discussion in \cite[\S
4.6.1]{nekrashevych:book:selfsimilar} but for completeness we
include a proof in \S 6.4.  \qed

\noindent{\bf Augmented trees.} If in addition $S$ is a good generating set for a contracting action
i.e., one which contains the nucleus and which is closed under restrictions, then  
the selfsimilarity complex inherits additional structure making 
it an  {\em augmented rooted tree} in the sense of V. Kaimanovich 
\cite{kaimanovich:sierpinski}:  
if vertices $u_1=x_1v_1$ and $u_2=x_2v_2$ at the same level are joined by a 
horizontal edge, then so are $v_1$ and $v_2$ (we allow $v_1=v_2$).   That is, the complex $\Sigma$ is built out of  ``squares'' of the form 
\begin{equation}
\label{eqn:augtree}
\begin{CD}
xw @>s >>(xw)^s = x^s w^{s|_x} \\
@AxAA @AAx^sA\\
w    @>>s|_x>w^{s|_x}
\end{CD}
\end{equation}
where the vertices $w$ and $w^{s|_x}$ possibly coincide.     It follows easily that any geodesic can be inductively modified so it 
is in {\em normal form}.  That is, it consists of first a (possibly empty) vertical 
segment traversed downward, then a horizontal segment, followed finally by a (possibly 
empty) vertical segment traversed upwards.  
\gap

\noindent{\bf Shift maps.}  This has the following geometric/dynamic implication. 
The right shift $\sigma: X^* \to X^*$ induces a cellular, label-preserving, degree $d$ covering map $F: \Sigma' \to \Sigma$.  On vertices, the left shift map $\tau$ simply pushes each vertex to the vertex immediately below it; for $k \leq |v|$ we denote $\tau^k(v)=v^{[-k]}$.  
Since the left and right shifts commute, $F\circ \tau = \tau \circ F$ on vertices.   Note that the map $F$ is  $1$-Lipshitz, hence so is any iterate of $F$.   
\gap

\noindent{\bf Geometric interpretation of contraction.} 
In the next chapter, we will make essential use of an interpretation of the contracting
condition into geometric language.   First, note that by concatenating diagrams as in (\ref{eqn:augtree}) 
along vertical edges, we see that for all $g \in G$ and for all $u \in X^*$, we have 
$|| \; g|_u \; || \leq ||g||$.
Now let $l \geq 1$ be a positive integer.  Suppose $u_1, u_2$ are vertices at the same level $n$, 
and suppose $l=\min_{g \in G}\{ ||g|| : u_1^g = u_2\}$, so that $l$ is the length of the shortest horizontal edge-path at level $n$ 
from $u_1$ to $u_2$.  
The contracting property implies that there is ``magic level'' $m(l)$ such that whenever $|v| \geq m(l)$, 
we have $g|_v \in \NNN$, and so $g|_v \in S$ since $S$ is assumed good.  In particular, $||g|_v|| \leq 1$.  
This means that in the selfsimilarity complex $\Sigma$, we have the following ``quadrilateral'' whose side lengths are indicated 
(the top horizontal segment need not be geodesic in $\Sigma$):
%commutative diagram 
\begin{equation}
\label{eqn:contracting}
\begin{CD}
u_1=vw @>l=||g || >>(vw)^g=v^gw^{g|_v}=u_2 \\
@A{|v|}AA @AA{|v^g|}A\\
w    @> >\leq 1>w^{g|_v}
\end{CD}
\end{equation}
Fix now $l \geq 0$ and $n \geq m(l)$, and consider the Schreier subgraph spanned by $X^n$, with the restricted graph 
``horizontal'' metric $d_{hor}$.  Then for arbitrary $u_1, u_2 \in X^n$, the restriction $\tau^{m(l)}: X^n \to X^{n-m(l)}$ satisfies
\[ d_{hor}(\tau^{m(l)}(u_1), \tau^{m(l)}(u_2)) \leq \lceil d_{hor}(u_1, u_2)/l \rceil\]  
where $\lceil x \rceil$ is the smallest integer greater than or equal to $x$.

We also note the following.  Let $l \geq 1$ be a positive integer, and suppose there exists a horizontal path in 
$\Sigma$ between vertices which is a geodesic of length $m(l)\cdot  l$ in the ambient metric on $\Sigma$ 
(here, we do not mean the restriction of the metric to a horizontal subset).  
By concatening $m(l)$ copies of the diagram in (\ref{eqn:contracting}) along vertical edges, it is clear that then 
$m(l)\cdot l \leq m(l) + m(l) + m(l) = 3\cdot m(l)$.  It follows that $l \leq 3$ and 
so the maximum length of a horizontal geodesic is bounded above by $H_\Sigma := m(3)\cdot 3$.  
From this observation the hyperbolicity of $\Sigma$ follows easily; cf. 
\cite[\S 3.8]{nekrashevych:book:selfsimilar}, \cite{kaimanovich:sierpinski} and \cite{kmp:gromov}. 
In addition, for any geodesic in normal form,  its horizontal segment has  length at most
$H_{\Sigma}$.

%\input sec4.tex

%\chapter{Finiteness theorems for selfimilarity complexes}
\section{Finiteness principles for selfsimilarity complexes}
\label{ch:finiteness}

In this chapter, we assume we are given a faithful level-transitive selfsimilar contracting recurrent action of a finitely generated group $G$ on the set of words $X^*$ in a finite alphabet $X$, and are given a symmetric generating set $S$ which is closed under restrictions and which contains the nucleus $\NNN$ of the action.  Let $\Sigma$ denote the selfsimilarity complex associated to $G$ and $S$; it is a Gromov hyperbolic augmented tree.   

The main result of this chapter is a fundamental finiteness result (Theorem \ref{thm:finiteness_principles}) concerning the induced dynamics $F: \Sigma'  \to \Sigma$ that is reminiscent of the finiteness of cone types for a Gromov hyperbolic group, which we now describe.

Let $G$ be a finitely generated group and $S$ a generating set such that if $s\in S$, then $s^{-1}\in S$.
Any element $g$ in $G$ can thus be written as as a product of elements, i.e. a  word, in $S$. 
The minimum number of elements of $S$ required to write $g$ as  product of elements of $S$ is 
by definition the word length $||g||$ of $g$.
Let $\Sigma$ denote the {\em Cayley graph} of $G$ with respect to $S$.
Its vertex set  is the set of elements of $G$, and $(g,g')\in G\times G$ defines an edge
if $g^{-1}g'\in S$. The edges of $\Sigma$ are naturally labelled by elements of $S$.
The length metric on $\Sigma$ in which each edge is isfometric to a unit interval turns $\Sigma$ 
into a proper, geodesic metric space on which $G$ acts isometrically by left-translation.  
The distance   from $g$ to the neutral element $e$ is the word length of $g$, 
and to any other vertex $g'$ is $||g^{-1}g'||$.   

The  {\em cone} $C_g$ of an element $g\in G$ is the set of vertices $w\in G$ such that
$g$ lies in a geodesic segment joining $e$ to $w$ i.e., $||w||=||g|| + ||g^{-1}w||$.

The {\em cone type} $T_g$  of $g$ is $T_g = g^{-1}C_g$ i.e., 
$T_g=\{w \in G  :  gw \in C_g\}.$   
By definition, if $T_g = T_h$, then the restriction 
\[ \phi: C_g \to C_h, \;\;\; x \mapsto hg^{-1}x\]
of the left-translation map $L_{hg^{-1}}: \Sigma \to \Sigma$ is a well-defined isometry.  

Given a positive integer $k$ and $g \in G$, the {\em $k$-tail} of $g$ is the set of those 
$h \in G$ for which $||h||\leq k$ and $||gh|| < ||g||$.  The following theorem was proved 
by Cannon in the case of cocompact Kleinian groups; the proof of the general case is very similar and 
may be found in \cite[Theorem III.$\Gamma$.2.18]{bridson:haefliger:book}. 

\begin{thm}[Cone types finite]
\label{thm:cone_types_finite} 
Let $G$ be a $\delta$-hyperbolic group.  Then the cone type of $g$ is determined by the $(2\delta+3)$-tail of $g$.   
Hence, there are only finitely many cone types.    
\end{thm}
That is, the geometry of the cone $C_g$ of $\Sigma$ is determined by a uniformly small amount of combinatorial data 
associated to $g$.  Another result with a similar flavor (and proof) is an analogous result 
for so-called {\em half-spaces} defined by geodesic segments; cf. \cite[Theorem 3.21]{cannon:swenson:characterization}.

\subsection{Subcomplexes, cones, shadows}
\label{secn:cones_etc}

We emphasize that we regard the selfsimilarity complex $\Sigma$ as a CW
1-complex with edges locally isometric to Euclidean unit intervals.  
Each 1-cell is equipped with a distinguished orientation and a
distinguished label drawn from the finite set $S$.  The local picture near each vertex at
level $>0$ is the same.   We denote by $|a-b|$ the distance between points $a,b$.  

The following observation will be of central importance in the following chapter.  
From the hyperbolicity of $\Sigma$, it follows that there exists a constant $C$ with the following property.  
Suppose $R_1, R_2: [0,\infty) \to \Sigma $ are two infinite vertical rays from the basepoint $o$.  
Then $R_1 \sim R_2$ if and only if $|R_1(t)-R_2(t)| \leq C$ for all $t \geq 0$.   
In our setting we have $C=1$; see \cite[Theorem 3.6.3]{nekrashevych:book:selfsimilar} and \S 6.1 below.  
Thus in addition to the scale of $H_\Sigma$, the maximum length
of a horizontal geodesic segment, the scale of 1 arises as an important
quantity as well.  

In the remainder of this section, we work exclusively in the complex $\Sigma$, equipped with the path-metric of \S 4.
In the next section, we transfer the results obtained here to the boundary $\bdry \Sigma$.  

\subsubsection{Horizontal, vertical}
The {\em level} of a subset of $\Sigma$ is its distance from the origin.  
A subset of $\Sigma$ is {\em horizontal} if all its vertices lie at the same level;
it is {\em vertical} if no two vertices have the same level.  
Given two vertices $a, b$ at the same level, the {\em horizontal distance} 
$|a-b|_{hor} = \min\{ l(\gamma) : \mbox{$\gamma$ is horizontal and joins $a$ to $b$}\}$.  
If $V \subset \Sigma$ is horizontal we set $\diam_{hor}(V)=\sup\{ |a-b|_{hor} : a, b \in V\}$.
\gap

\noindent{\bf Notation.}  If $w=uv$ with $|u|=k$ we set, following Kaimanovich, $v=w^{[-k]}$.   
The notation $[ov]$ stands for the vertical geodesic ray joining the basepoint $o$ and a vertex $v$; it is unique.  
More generally, $[a,b]$ denotes a geodesic segment joining $a$ and $b$; it is unique if it is vertical.  
For a horizontal subset $V \subset \Sigma$ we will write $|V|$ for the level of $V$.

\subsubsection{Balls and neighborhoods}
Given a vertex $v$ of $\Sigma$ and an integer $r \geq 1$ we denote by
$B(v,r)$ the set of vertices $w$ such that $|v-w| \leq r$ and by
$B_{hor}(v,r)$ the intersection of $B(v,r)$ with the set of vertices at
level $|v|$.  {\em We emphasize that $B(v,r)$ is a set of vertices and not a
union of vertices and edges.}  Given a subset $A$ of the vertices
of $\Sigma$ we denote by  $B(A,r) = \union_{a \in A}B(a,r)$ and set 
$B_{hor}(A,r)=\union_{a \in A}B_{hor}(a,r)$ similarly.  

\begin{lemma}[Balls map to balls]
\label{lemma:balls_map_to_balls}
Fix an integer $k>0$.  Then for all $v \in X^*$, all
$\tilde{v} \in F^{-k}(v)$, and all $r \leq \diam_{hor} X^{|v|}$ we have $F^k(B_{hor}(\tilde{v},
r)) = B_{hor}(v,r)$. 
\end{lemma}

\pf The inclusion $\subset$ follows immediately
since $F^k$ is cellular, and the inclusion $\supset$ follows by
path-lifting, using the fact that $F^k$ is a covering map.  
\qed

\subsubsection{Distances between subsets; neighborhoods}
Let $A,B$ be subsets of $\Sigma$.  We set 
\[ |A-B| = \inf\{|a-b| \; : \; a \in A, b \in B\}.\]
If $A$ and $B$ are horizontal at the same level, we analogously set 
\[ |A-B|_{hor} = \inf\{|a-b|_{hor} \; : \; a \in A, b \in B\}.\]
Thus, the level of a set $A$ is  $|A| = |A - o|$.
For $r>0$ denote by $N_r(A)=\{x \in\Sigma : \exists a \in A, \; |x-a| < r\}$ the open $r$-neighborhood of $A$.  
The {\em Hausdorff distance} between $A$ and $B$ is  
\[ d_{haus}(A,B)=\inf\{r : A \subset N_r(B), B \subset N_r(A)\}.  \]

\subsubsection{Cones}
\begin{figure}
\label{fig:cone_umb}
\includegraphics[width=5in]{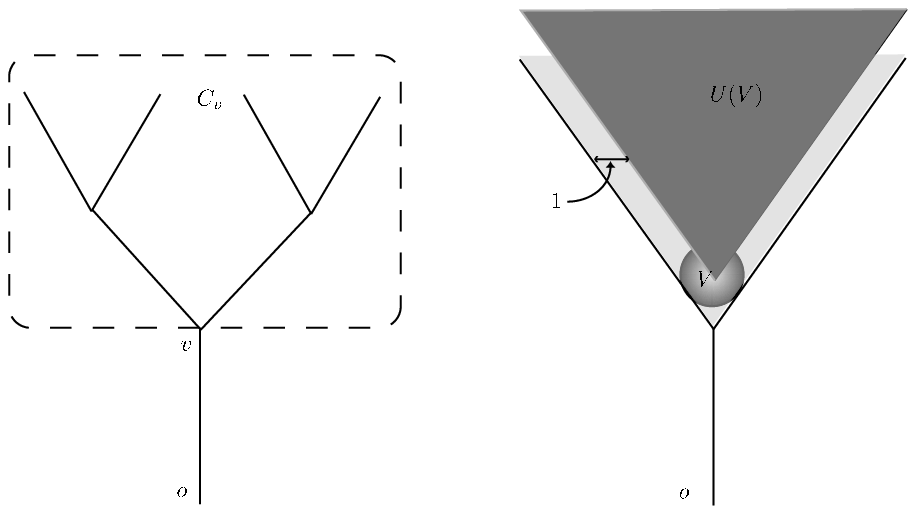}
%\psfrag{3}{$o$}
%\psfrag{1}{$v$}
%\psfrag{2}{$C_v$}
%\psfrag{4}{$1$}
%\psfrag{5}{$U(V)$}
%\psfrag{V}{$V$}
%\includegraphics[width=5in]{umb}
\begin{center}
%\caption
{\sf Figure \ref{fig:cone_umb}: 
At left, the cone $C_v$. At right, the shadow $S(V)$ (light) and the umbra $U(V)$ (dark).}
\end{center}
\end{figure}

Given a vertex $v \in \Sigma$ we define the {\em cone} $C_v$ as the
subset of vertices of $\Sigma$ given by 
\[ C_v = \{wv | w \in X^*\}.\]
See Figure \ref{fig:cone_umb}.

Note that given vertices $v_1, v_2$ we have that exactly one of the
following holds: 
\[ C_{v_1} \intersect C_{v_2} = \emptyset \mbox{ or } C_{v_1}
\subset C_{v_2} \mbox{ or } C_{v_2} \subset C_{v_1}.\]

If $v=wx$ then 
\[ F(C_v) = C_w\]
since $F$ acts as the right shift on vertices.  

\subsubsection{Induced subcomplexes}
Given a set $V$ of vertices of $\Sigma$, we denote by $\Sigma(V)$ the
{\em induced subcomplex} on $V$, i.e. the subcomplex of $\Sigma$ whose 0-cells
are the elements of $V$ and whose 1-cells consist of those 1-cells of
$\Sigma$ whose endpoints lie in $V$. Often we will denote such
subcomplexes by the symbol $\Gamma$.  

The {\em isomorphism type} of a subcomplex is its equivalence class under 
the equivalence relation generated by label- and orientation-preserving cellular isomorphisms.  

A connected component of an induced subcomplex is again an induced
subcomplex.  In particular, if $\Gamma$ is a connected induced subcomplex and
$\widetilde{\Gamma}$ is a component of $F^{-k}(\Gamma)$ then $\widetilde{\Gamma}$ is
again an induced subcomplex.  The following lemma will be important. 

\begin{lemma}[At least two apart]
\label{lemma:at_least_two_apart}
Let $V \subset X^*$ be a set of vertices, $\Gamma = \Sigma(V)$, and let 
$\Gamma_i$, $i=1,2$ be two connected components of $\Gamma$.
Then $\Gamma_1 \neq \Gamma_2$ iff
$|\Gamma_1 - \Gamma_2|
\geq 2$.  
\end{lemma}

\pf If $v_i \in \Gamma_i$ satisfy $|v_1-v_2|=1$, then $v_1$ and $v_2$ are
incident to the same edge in $\Sigma$.  This edge lies in $\Gamma$ since $\Gamma$
is an induced subcomplex.  Therefore $\Gamma_1 =\Gamma_2$.  
\qed

\subsubsection{Shadows}
Given a set of vertices $V$ we define the {\em shadow} $S(V)$ by 
\[ S(V) =\Sigma(\union_{v \in V}C_v).\]
See Figure \ref{fig:cone_umb}.
The shadow $S(B_{hor}(v,r))$ of a horizontal ball of radius $r$ about $v$ will be denoted more succinctly by $S(v,r)$.

If $V$ is a point or a horizontal ball, a shadow is somewhat analogous to a half-space in the Cayley graph of a Gromov hyperbolic group.  Many of the geometric results below have analogies in the group setting; compare \cite[\S 3.3]{cannon:swenson:characterization}.

The following properties are easily verified.  Below, $^cA$ denotes the complement of $A$; a set $A$ is $Q$-{\em quasiconvex} if every geodesic between points in $A$ lies in a $Q$-neighborhood of $A$.  

\begin{lemma}[Properties of shadows]
\label{lemma:properties_of_shadows}
Let $V \subset X^*$.  
\be

\item $u \in {^cS(V)}\intersect X^* \implies [ou] \subset {^cS(V)}$.  If
in addition $|u| >\max\{|v| : v \in V\}$ then $C_u \subset {^cS(V)}$.  

\item If $V$ is horizontal, $u_i \in S(V)$, $|u_1|=|u_2| > |V|$, and
$|u_1-u_2|=1$ (i.e. they are joined by a horizontal edge) then
$u_i^{[-1]}$, $i=1,2$ are in $S(V)$ and either coincide or are joined by a horizontal edge.
Hence, any horizontal path in $S(V)$ can be "pushed down" to a
horizontal path in $\Sigma(V)$ of the same or shorter length.  

\item $S(V)$ is connected iff $\Sigma(V)$ is connected. 

\item $S(V)$ is $Q$-quasiconvex where $Q=\max\{\roundup{\diam V/2},\roundup{H_\Sigma/2}\}+1$. 

\item If $V_1, V_2$ are horizontal subsets at the same level, and $|V_1-V_2| \geq 2$ then
$|S(V_1)-S(V_2)| \geq 2$.  

\item 
Let $\wtV, V$ be horizontal vertex sets and $k >0$ an integer.
Then  $ F^{k}(\Sigma(\wtV))=\Sigma(V) \iff 
F^k(S(\wtV)) = S(V)$.

\item Suppose $\Gamma = \Sigma(V)$ is connected.  Let $\wtV =
F^{-k}(V)$, $\widetilde{\Gamma} = \Sigma(\wtV)$.  Then
\be
\item 
$F^k(\widetilde{\Gamma}) = \Gamma$;
\item if $\widetilde{\Gamma}_i$ is a connected component of $\widetilde{\Gamma}$ then 
$F^k(\widetilde{\Gamma}_i) = \Gamma$;
\item if
$\widetilde{\Gamma}_i$, $i=1,2$ are distinct components of
$\widetilde{\Gamma}$ then $|\widetilde{\Gamma}_1 - \widetilde{\Gamma}_2|
\geq 2$ and hence $|S(\wtV_1)-S(\wtV_2)| \geq 2$.
\eb
\eb
\end{lemma}

We give proofs of the not-so-straightforward assertions. 

\pf  Assertion {\bf 2.}  follows immediately from the definitions, since we have assumed $\Sigma$ 
is an augmented tree.

{\bf 3.} Suppose $\Sigma(V)$ is connected.   It is enough to show that any pair of
vertices $w_1, w_2 \in S(V)
\intersect X^*$ are joined by a path in $S(V)$.  By definition, there are $v_i \in V$ such
that 
$w_i \in C_{v_i}\subset S(V)$.   Since $S(V)$ is an induced subcomplex this implies 
that the vertical geodesic segment $[v_iw_i] \subset S(V)$.  Since $v_i \in
\Sigma(V)$ and $\Sigma(V)$ is connected, there is a horizontal path joining 
$v_1$ to $v_2$ in $\Sigma(V) \subset S(V)$.  Hence $w_1, w_2$ are joined by a path in
$S(V)$.  

Conversely, suppose $\Sigma(V)$ is disconnected.  Write $\Sigma(V)=\Gamma = \Gamma_1 \sqcup
\Gamma_2$ where $\Gamma_1, \Gamma_2$ are disjoint nonempty induced subcomplexes.  Then $|\Gamma_1 - \Gamma_2| \geq 2$ by Lemma \ref{lemma:at_least_two_apart}. 
If $w_i \in C_{v_i}$ for some $v_i\in\Gamma_i$, then $|w_1-w_2| \geq 2$  since otherwise pushing down to the level of
$V$ will yield points at distance at most $1$, a contradiction.  
Hence $S(V)$ is disconnected.   

{\bf 4.} Let $u_1, u_2$ be two vertices in $S(V)$ and let $\gamma$ be a normal 
form geodesic joining them. 
Suppose the horizontal portion $\gamma'$ of this geodesic lies at a level larger than $|V|$.  
Then we may write $\gamma = \gamma_1 * \gamma' * \gamma_2$ where $\gamma_i$ are vertical and lie in $S(V)$.  
Since the length of $\gamma'$ is at most $H_\Sigma$ we have that $\gamma$ lies in a 
$\roundup{H_\Sigma/2}+1$-neighborhood of $S(V)$.   
Otherwise, the definition of shadows implies that $\gamma$ meets $V$ in two points $v_1, v_2$, 
so that $\gamma = \gamma_1 * \gamma' * \gamma_2$ where the $\gamma_i$ are vertical and in $S(V)$ and $\gamma'$ is a normal form geodesic joining $v_1$ and $v_2$.  Again, the length of $\gamma'$ is at most $\diam V$, hence every point on it lies in a $\roundup{\diam V/2}+1$-neighborhood of $V$, hence of $S(V)$ as well.

{\bf 5.} This follows from the previous paragraph on induced subcomplexes.

{\bf 6.} The right and left shifts commute. 
\qed

This subsection concludes with some technical
results regarding preimages of balls and of their corresponding shadows.  The former will later be used
in the proof of Theorem \ref{thm:bounded_degree}; the latter will be used in Chapter \ref{chap:bdry}.  

\begin{lemma}[Pullbacks of balls]  
\label{lemma:pullbacks_of_balls}
Suppose that $k, r >0$ are integers and $v\in X^*$.    
Let $\Gamma=\Sigma(B_{hor}(v,r))$.   
Then 
\be
\item $F^{-k}(B_{hor}(v,r)) = \bigcup_{F^k(\tilde{v})=v}B_{hor}(\tilde{v},r)$.

\item If $\widetilde{\Gamma}$ is a connected component of $F^{-k}(\Gamma)$ and $\wtV =
F^{-k}(v) \intersect \widetilde{\Gamma}$ then 
\[ \widetilde{\Gamma} = \Sigma\left(\union_{\tilde{v}\in \wtV} B_{hor}(\tilde{v},r)\right).\]
Moreover:  given any pair $\tilde{v}', \tilde{v}'' \in \wtV$ there are
$\tilde{v}_i
\in
\wtV$,
$i=0, \ldots, p$ such that $\tilde{v}_0 = \tilde{v}', \tilde{v}_p = \tilde{v}''$, and for each
$1 \leq i \leq p$ we have
$|\tilde{v}_{i-1}-\tilde{v}_i|_{hor} \leq 2r+1$.  

\item If $\widetilde{\Gamma}_i$, $i=1,2$ are two distinct components of $F^{-k}(\Gamma)$
then $|\widetilde{\Gamma}_1-\widetilde{\Gamma_2}| \geq 2$.  
\eb
\end{lemma}

\noindent{\bf Remark:}  We have always $\diam_{hor} V \leq \diam_{hor}\Sigma(V) \leq 
\diam_{hor} V+1$.  Equality in the upper bound can occur, e.g. when there are loops at each
vertex of a pair of vertices whose horizontal distance realizes $\diam_{hor} V$, and one or both of these loops fails to lift under $F$.  
\gap

\pf {\bf 1.}  This follows easily from [Lemma \ref{lemma:balls_map_to_balls}, Balls map to
balls].  {\bf 3.}  is a special case of Lemma \ref{lemma:properties_of_shadows}(7).  

We now prove the inclusion $\subset$ in {\bf 2.}.  First, 
the set of vertices $\tilde{V}$ is clearly contained in the induced
complex on the right-hand side.  Let $\tilde{w} \in \tilde{\tau} \subset \wtG$ map under 
$F^k$ to $w$, where $\tilde{\tau}$ is an edge.  Let $\tau =F^k(\tilde{\tau}) \subset
\Gamma$ be a closed one-cell of $\Sigma$ containing
$w$ and let
$v_i', i=1,2$ denote the (possibly indistinct) vertices comprising $\bdry \tau$.  
Since $\Gamma = \Sigma(B_{hor}(v,r))$ is connected, there are horizontal paths $\gamma_i
\subset
\Gamma$ which are geodesics with respect to the intrinsic path metric in $\Gamma$ joining $v$
to $v_i'$ of length $\leq r$.  Let $e \subset \Gamma$ be the edge-path given by 
$\gamma_1 *\tau * \gamma_2^{-1}$  where $*$ denotes concatenation.  Since $\wtG$ is a component of the
preimage of $\Gamma$, the map $\wtG \to \Gamma$ is a covering map and thus has the
path-lifting property.  Let $\tilde{e}$ be the unique lift of this edge-path containing
$\tilde{w}$.  Then $\tilde{e}$ contains as a sub-edge-path an edge-path given by
$\tilde{\gamma}_1 * \tilde{\tau} *
\tilde{\gamma_2}$ where $\tilde{\gamma_i} \mapsto \gamma_i$  under
$F^k$.  Thus $\tilde{\gamma}_i$ joins some $\tilde{v}_i \in \wtV$ to some $\tilde{v}_i' \in
F^{-k}(v_i')$.  Hence $\tilde{v}_i' \in B_{hor}(\tilde{v}_i,r)$ and so $\tilde{\tau}
\subset \Sigma(\cup_{\tilde{v} \in \wtV}B_{hor}(\tilde{v},r))$.

To prove the other inclusion, suppose $\tilde{w} \in \Sigma(\cup_{\tilde{v} \in
\wtV}B_{hor}(\tilde{v},r))$.  Then $\tilde{w} \in \tilde{\tau}$, a 1-cell with $\bdry \tilde{\tau} =
\{\tilde{v}_1', \tilde{v}_2'\}$.  Let $\tilde{\gamma}_i$ be a geodesic in
$\Sigma(\cup_{\tilde{v} \in \wtV}B_{hor}(\tilde{v},r))$ joining $\tilde{v}_i'$ to some $\tilde{v}_i
\in \wtV$, so that $|\tilde{v}_i'-\tilde{v}_i|=|\tilde{v}_i'-\wtV|$, $i=1,2$.  Since $\wtV \subset \wtG$ 
and $\wtG$ is connected, there exists a path
$\tilde{\gamma}$ in $\wtG$ joining $\tilde{v}_1$ to $\tilde{v}_2$.  Thus, the edge-path
$\tilde{\gamma}_1 * \tilde{\gamma} * \tilde{\gamma}_2 ^{-1}$ lies in $\wtG$ and joins
$\tilde{v}_1'$ to $\tilde{v}_2'$.  Therefore $\tilde{v}_i' \in \wtG$ and hence $\tilde{\tau}
\subset \wtG$ since $\wtG$ is an induced subcomplex.  Hence $\tilde{w}$, which lies in
$\tilde{\tau}$, is in $\wtG$. 

To prove the second assertion, suppose there is a nontrivial partition of $\wtV$ into
disjoint subsets $\wtV_1 \union \wtV_2$ such that $|\wtV_1 - \wtV_2| > 2r+1$.  This implies
that 
\[ |\union_{\tilde{v}_1 \in \wtV_1}B_{hor}(\tilde{v}_1,r) - \union_{\tilde{v}_2 \in
\wtV_2}B_{hor}(\tilde{v}_2, r)| > 1.\]
Hence the induced subcomplex $\wtG$ is not connected.   
\qed

From the properties of shadows (Lemma \ref{lemma:properties_of_shadows}) we get a
result similar to the previous one for shadows.  

\begin{lemma}[Pullbacks of shadows]
\label{lemma:pullbacks_of_shadows}
Fix a positive integer $r>0$, and let $v$ be a vertex at level $>r$.  
Denote by $S$ the shadow $S(v,r)$.  Then for all $k>0$, 
\be 

\item $ F^{-k}(S)\intersect X^* =
\bigcup_{F^k(\tilde{v})=v}S(\tilde{v},r)\intersect X^*$.  

\item If $\widetilde{S}$ is a component of $F^{-k}(S)$ and $\wtV = F^{-k}(v)
\intersect \widetilde{S}$ then 
\[ \widetilde{S} = S(\wtV) = \Sigma\left(\union_{\tilde{v} \in
\wtV}S(\tilde{v},r)\intersect X^*\right).\]

\item If $\widetilde{S}_i$, $i=1,2$ are two distinct components of $F^{-k}(S)$
then $|\widetilde{S}_1-\widetilde{S}_2| \geq 2$. 
\eb 
\end{lemma}

\subsubsection{Umbrae}
\label{subsecn:umbrae} 

The darkest part of a shadow is called its {\em umbra}.

Given a vertex set $V$ we define its {\em umbra} $U(V)$ as the subcomplex
induced by the set of vertices $u \in S(V) \intersect X^*$ such that
$|u-{^cS(V)}| > 1$; see Figure \ref{fig:cone_umb}.  Then 
\[ X^*\intersect(S(V)-U(V)) = V \sqcup \{ z \in X^* : |z-{^cS(V)}|_{hor}=1\}.\]

The lemma below gives a sufficient condition for an umbra to be nonempty.  
Easy examples (e.g. \cite[Figure 3.1]{nekrashevych:book:selfsimilar}) show that it is not, however, a necessary condition.

\begin{lemma}[Unit ball inside implies umbra is nonempty]
\label{lemma:umbra_is_nonempty}
If $B_{hor}(v,1) \subset V$ then for all $x \in X$ we have $C_{xv} \subset U(V)$.
\end{lemma}

\pf Briefly: the map $G$ induced by the left-shift is 1-Lipschitz.  In more detail: let $y \in {^cS(V)}\intersect X^*$ and suppose $|wv-y|_{hor}=1$ for some word $w\in X^*\setminus\{\emptyset\}$.  
Let $y' = y^{[-|w|]}$.  
Then $|v-y'| \leq 1$ by Lemma \ref{lemma:properties_of_shadows},(2), and so $y' \in B_{hor}(v,1)$. But 
since $y \in {^cS(V)}$, we have $y' 
\in {^cS(V)}$ by Lemma \ref{lemma:properties_of_shadows}(1) and so $B_{hor}(v,1)
\not\subset V$.  
\qed 

Since the maps $F, G$ induced by the right and left shifts commute, Lemma \ref{lemma:properties_of_shadows} implies 

\begin{lemma}[Naturality of umbrae]
\label{lemma:naturality_of_umbrae}
When restricted to components, umbrae are natural with respect to the dynamics.
That is: suppose $V$ is a vertex set, $\Gamma = \Sigma(V)$ is connected, $\widetilde{\Gamma}$ is a
connected component of $F^{-k}(\Gamma)$, $\wtV = F^{-k}(V) \intersect
\widetilde{\Gamma}$, and let $U=U(V), \wtU=U(\wtV)$.  Then $F^k(\wtU)=U$.   
\end{lemma}

The following observation will be used later (in the proof of Lemma
\ref{lemma:umbrae_contain_shadows})  to show that boundaries at infinity 
of umbrae are open.   

\begin{lemma}[Unit speed penetration]
\label{lemma:unit_speed_penetration}
Suppose $u \in U(V)$ and $v \in C_u$.  Then 
\[ |v-{^cS(V)}| \geq |v|-(|u|+m(H_\Sigma)).\]
\end{lemma}

\pf Let $y \in {^cS(V)}$.   If $|y| \leq |u|$ then $|v-y| \geq |v|-|u|$ and we are done.  
So now suppose $|y| > |u|$.    We have $[oy] \subset {^cS(V)}$ by Lemma
\ref{lemma:properties_of_shadows}(1).  Let
$[v'y']$ be the horizontal segment of a normal form geodesic joining $v$
and $y$.  The level $l$ of this segment is at most  $|u|+m(H_\Sigma)$,
for otherwise, pushing this segment down to the level of $u$ will
yield $|u-[oy]| \leq 1 \implies |u-{^cS(V)}|\leq 1$, contradicting $u \in
U(V)$.   Thus $l \leq |u|+m(H_\Sigma)$.  But $|v-y| \geq |v|-l \geq |v|-(|u| + m(H_\Sigma)$),
yielding the result. 
\qed

\subsection{Isometry types of shadows}
In this section, we prove a static finiteness result:  the isometry
type of a shadow $S(V)$ is determined by the isomorphism type of the induced subcomplex
$\Sigma(\widehat{V})$ of a certain neighborhood
$\widehat{V}$ of the defining subset $V$ called the {\em hull} of $V$.  
In the next section we add dynamics and prove a related finiteness result.  

We begin with an easier result which has the same flavor.  

\begin{thm}[Cone types]
\label{thm:cones_finite}
Suppose $v_i$ are vertices of $\Sigma$ and let $\Gamma_i =
\Sigma(B_{hor}(v_i,H_\Sigma))$, $i=1,2$.  If $\phi: \Gamma_1 \to \Gamma_2$ is an
isomorphism of 1-complexes with labelled edges such that $v_2=\phi(v_1)$, then the map 
\[ \phi: C_{v_1} \to C_{v_2} \]
given by 
\[ wv_1 \mapsto  wv_2, \;\;\;w \in X^* \]
is an isometry with respect to the distance function of $\Sigma$.  
\end{thm}

\begin{cor}[Finitely many cone types]
\label{cor:cones_finite}
There are only finitely many isometry classes of cone
types.  
\end{cor}

{\noindent\bf Proof of Corollary \ref{cor:cones_finite}:}  There are only finitely many possible isomorphism types
of labelled 1-complexes of the form $B_{hor}(v, H_\Sigma)$ since the horizontal valence of any vertex
is bounded by $\# S$.  
\qed

We view this as an analog of Cannon's observation of the finiteness of cone types for (Gromov)  
hyperbolic groups; cf. Theorem \ref{thm:cone_types_finite}.
Compare with \cite[\S 3.3]{nekrashevych:book:selfsimilar}.
%, where it is shown that there are at most $2^{\#\NNN}$ homeomorphism types of boundaries of cones $C_v$.
\gap

{\noindent\bf Proof of Theorem \ref{thm:cones_finite}:} For convenience we write $v=v_1$, $v^\phi=v_2$, $\Gamma=\Gamma_1$, $\Gamma^\phi =
\Gamma_2$, and $\phi(wv) = wv^\phi$.    

Let $w_1, w_2 \in C_v$.  Let $\gamma$ be a normal form geodesic joining
$w_1, w_2$.  Then $q \geq |v|$ where $q$ is the level of
the horizontal segment of $\gamma$.  

\begin{figure}
\label{fig:conefig1}
\includegraphics[width=5in]{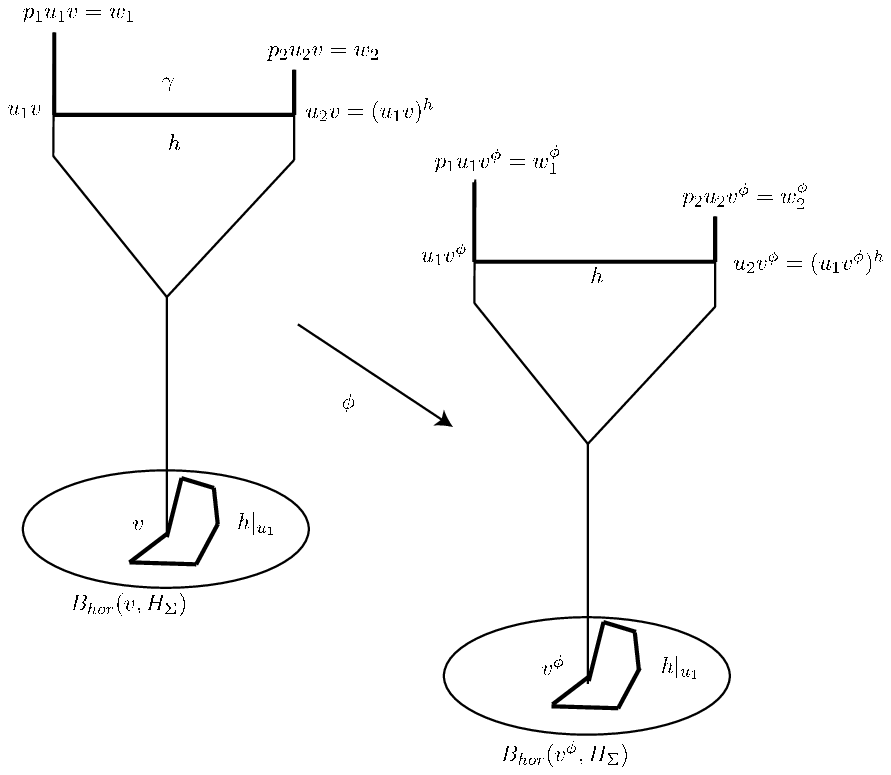}
\begin{center}
{\sf Figure \ref{fig:conefig1}: Finitely many cone types.  The path $\gamma$ is shown in bold.  The edge-path based at $v$ representing $h|_{u_1}$ is contained in the horizontal ball of radius $H_\Sigma$ based at $v$ and is shown in bold. }
\end{center}
\end{figure}

Write $w_i = p_i u_i v$ where $|u_i| = q-|v|$, $i=1, 2$ (see Figure \ref{fig:conefig1}).  The horizontal
segment of $\gamma$ joins $u_1v$ to $u_2v$ and is a horizontal edge-path
representing an element $h \in G$ where $||h||  \leq H_\Sigma$.  
Then 
\[ |w_1-w_2| = |p_1|  + ||h||  + |p_2|.\]

The theorem follows once we show 
%\begin{equation}
%\label{eqn:h_and_phi}
$(u_1v^\phi)^h = u_2v^\phi$.
%\end{equation}
Assuming this, it follows that the edge-path from $w_1^\phi = p_1u_1v^\phi$ to
$w_2^\phi=p_2u_2v^\phi$ given by going first down towards the root along
$p_1$ to $u_1v^\phi$ , over via the horizontal edge-path representing $h$ to $(u_1v^\phi)^h = u_2v^\phi$, and
then up along $p_2$ has length equal to $|w_1-w_2|$ and thus $|w_1^\phi -
w_2^\phi| \leq |w_1 - w_2|$.  The inequality $|w_1-w_2| \leq |w_1^\phi - w_2^\phi|$ then
follows by considering $\phi^{-1}$.  

We now show $(u_1v^\phi)^h = u_2v^\phi$.
By the definition of a selfsimilar action, 
\[ u_2v=(u_1v)^h  =  u_1^h v^{h|_{u_1}} \implies u_2=u_1^h,\;\;\;
v=v^{h|_{u_1}}.\]  We have $||h|_{u_1}|| \leq ||h|| \leq H_\Sigma$.  Recalling the definition of $\Gamma$, it follows that  the restriction $h|_{u_1}$
is represented by an edge-path in $\Gamma$ starting at $v$.  It is a loop since
$v^{h|_{u_1}}=v$.  Since $\phi: \Gamma \to \Gamma^\phi$
is an isomorphism of labelled graphs, the same edge-path representing
$h|_{u_1}$, when starting at $v^\phi$, is also a loop.   Hence 
$(v^\phi)^{h|_{u_1}} = v^\phi$.
So $ (u_1v^\phi)^h = u_1^h(v^\phi)^{h|_{u_1}} = u_2v^\phi$
and the proof is complete.
\qed

\noindent{\bf Hulls.}  
Fix a positive integer $D$.  Given $V \subset X^n$ a horizontal set of vertices of diameter
$\leq D$, the {\em $D$-hull} of $V$ is defined as
\[ \widehat{V} = \bigcup_{i=0}^{\roundup{D/2}} \bigcup_{v' \in
V^{[-i]} } B_{hor}(v', H_\Sigma) \] 
where $V^{[-i]} = \{v' | v' = v^{[-i]}, v \in V\}$.

If $D=\diam V$, the hull of $V$ is approximately a horizontal $H_\Sigma$-neighborhood of
the convex  hull of $V$.  The hull of $V$ contains all vertices lying on normal form geodesic
segments joining elements of $V$.  To see this, let $\gamma$ be any such geodesic; it has length
at most  $D$ and so its vertical portions each have length bounded by
$\roundup{D/2}$.  Then 
$|V|-\roundup{D/2} \leq |\gamma| \leq |V|$.  Since the horizontal portion has length at most
$H_\Sigma$ it follows that the vertices of $\gamma$ lie in the hull.  Hence, 
the induced subgraph  $\Sigma(\widehat{V})$ of the hull of $V$ is always connected.  

\begin{lemma}[Naturality of hulls]
\label{lemma:naturality_of_hulls}
If $F^k(\wtV) = V$ then $F^k(\widehat{\wtV})=\widehat{V}$. 
\end{lemma}

\pf The assumption implies $F^k(\wtV^{[-i]}) = V^{[-i]}$, $i=1, \ldots,
\roundup{D/2}$.  The lemma then follows by Lemma \ref{lemma:balls_map_to_balls}, [Balls map
to balls].  
\qed

We say that two induced subcomplexes of $\Sigma$ are {\em isometric}
if there is a cellular isomorphism between them which is distance-preserving.
The relation generated by this property is clearly an equivalence relation.  
The {\em isometry type} of a subcomplex is defined as its corresponding equivalence class.   

The theorem below enunciates the following principle:  the isometry type of a shadow of a
horizontal subset is determined by the isomorphism type of a (suitably large, depending on 
its diameter) associated hull.  

\begin{thm}[Finitely many shadow types]
\label{thm:shadow_type_determined_by_hull}  
For $i=1,2$, suppose $V_i$ are horizontal sets
of vertices of $\Sigma$ of diameter $\leq D$, let $\widehat{V}_i$ be the corresponding
$D$-hulls, and let
$\Gamma_i =
\Sigma(\widehat{V}_i)$.   If $\phi: \Gamma_1 \to \Gamma_2$ is an isomorphism (of 1-complexes
with labelled edges) then
the map 
\[ \phi: S(V_1) \to S(V_2) \]
given by 
\[ \phi(wv) = w\phi(v), v \in V_1, w \in X^*\]
is an isometry with respect to the distance function of $\Sigma$ and an
isomorphism of complexes.  

\end{thm} 
\gap

The proof is more or less exactly the same as the proof of Theorem
\ref{thm:cones_finite} above.

\pf Again write $V = V_1, V^\phi = V_2, \Gamma=\Gamma_1, \Gamma^\phi =
\Gamma_2$.  Suppose $w_1, w_2
\in S(V)$.  Let $v_i \in V \intersect [ow_i], i=1,2$.  

Conceptually it is easiest to consider now two cases.
\gap

{\bf Case 1:} There exists a normal form geodesic $\gamma$ from $w_1$ to
$w_2$ whose level (that is, the level of its horizontal part), $l$, is at least $|V|$.  

Write 
\[ w_i = p_i u_i v_i , \;\;\; i=1, 2, \;\; |u_1|=|u_2|, |v_1|=|v_2|, \;\;
|u_1v_1|=l.
\] Then by definition, there is some $h\in G$ such that $(u_1v_1)^h = u_2v_2 = u_1^h v_1^{h|_{u_1}}$ where $||h|| \leq H_\Sigma$. 
The restriction $h|_{u_1}$ again satisfies $|| \; h|_{u_1} \; || \leq H_\Sigma$ and
so viewed as an edge-path based at $v_1$ is a path of length $\leq H_\Sigma$.  
So $|v_1 - v_2| \leq H_\Sigma$.  Thus this edge-path lies in $\Gamma$.  Since $\phi: \Gamma \to
\Gamma^\phi$ is an isomorphism, we find that
$(v_1^\phi)^{h|_{u_1}} = v_2^\phi$ and consequently $(u_1v_1^\phi)^h  =
u_2v_2^\phi$.  With the same reasoning as in the previous proof we
conclude $|w_1^\phi - w_2^\phi|
\leq |w_1 - w_2|$.  By symmetry, equality holds.

{\bf Case 2:}  Otherwise, the level of the horizontal part
of $\gamma$ is strictly less than $|V|$.   The subsegment $\gamma'$ of $\gamma$ joining
$v_1$ to
$v_2$ is a geodesic, hence has length bounded by $\diam V$.  
The vertical distance between the horizontal segment of $\gamma$ and $V$
is at most $\roundup{D/2}$.  Hence $\gamma' \subset \Gamma$.  
If we write $w_i=p_i v_i$, then the concatenation of the geodesic segments $[p_1v_1^\phi,v_1^\phi]$, 
$\phi(\gamma')$ and
$[v_2^\phi,p_2v_2^\phi]$ is a curve joining $\phi(w_1)$ to $\phi(w_2)$ of length 
$|w_1-w_2|$ so
$|w_1^\phi-w_2^\phi|\leq |w_1-w_2|$.   By symmetry, equality holds.  
\qed

\subsection{Isometry types of maps between shadows}
In this section, we prove a dynamical version of Theorem
\ref{thm:shadow_type_determined_by_hull}.  

Let $F_i: \widetilde{Z}_i \to Z_i$, $i=1,2$  be maps between metric spaces.  We say that $F_1,
F_2$ are {\em isometrically isomorphic as maps} if there are isometries $\tilde{\phi}: 
\widetilde{Z}_1 \to \widetilde{Z}_2$ and $\phi: Z_1 \to Z_2$ such that the diagram 
%commutative diagram 
\[ 
\begin{array}{ccc}
\widetilde{Z}_1 & \stackrel{\tilde{\phi}}{\longrightarrow} & \widetilde{Z}_2
\\
F_1 \downarrow & \; & \downarrow F_2 \\
Z_1& \stackrel{\phi}{\longrightarrow} &Z_2\\
\end{array}
\] 
commutes.  
The {\em isometry type} of a map $F: \widetilde{Z} \to Z$ between metric spaces is its equivalence
class under this equivalence relation.  

We define similarly the {\em isomorphism type} of a cellular label-preserving map between CW
1-complexes with labelled edges.  

The following theorem says that the isometry type of a map between shadows is determined by the isomorphism type of the map
between induced subcomplexes of hulls.  

\begin{thm}[Shadow map determined by hull map]
\label{thm:shadow_map_determined_by_hull_map}
Suppose $V_i, \wtV_i$ are horizontal vertex sets of diameter $\leq D$, let
$\widehat{V}_i, \widehat{\wtV}_i$ denote the corresponding $D$-hulls, and suppose 
$F^{k_i}:
\Sigma(\wtV_i) \to \Sigma(V_i)$.  Suppose $\tilde{\phi}: \Sigma(\widehat{\wtV}_1) \to
\Sigma(\widehat{\wtV}_2)$ and
$\phi:
\Sigma(\widehat{V}_1) \to \Sigma(\widehat{V}_2)$ are isomorphisms of labelled complexes which
satisfy
$F^{k_2}\circ \widetilde{\phi} =
\phi \circ F^{k_1}
$ on $\Sigma(\widehat{\wtV}_1)$.  Then $F^{k_i}: S(\wtV_i) \to S(V_i)$, $i=1,2$ 
and the isometries $\tilde{\phi}: S(\wtV_1) \to S(\wtV_2)$, $\phi:S(V_1)
\to S(V_2)$ given by Theorem \ref{thm:shadow_type_determined_by_hull} also satisfy $F^{k_2}\circ \widetilde{\phi} =
\phi \circ F^{k_1}
$.  
\end{thm}

\pf Lemma \ref{lemma:properties_of_shadows}(6) implies that $F^{k_i}: S(\wtV_i) \to
S(V_i)$, $i=1,2$.  It remains only to prove the claimed commutativity
property.  

Let $\tilde{u}_1 \in S(\wtV_1)$ and write (as in the
definition of $\tilde{\phi}$) $\tilde{u}_1 = w\tilde{u}_1'$ where $\tilde{u}_1'
\in \wtV_1$ and $w \in X^*$.    Let $u_1 =
F^{k_1}(\tilde{u}_1)$ and
$u_1' = F^{k_1}(\tilde{u}_1')$, so that $F^{k_1}(\tilde{u}_1) = wu_1'$.  
Let $l=|w|$.  

The desired commutativity follows immediately since the two shifts commute: 

\begin{center}
\includegraphics[width=3in]{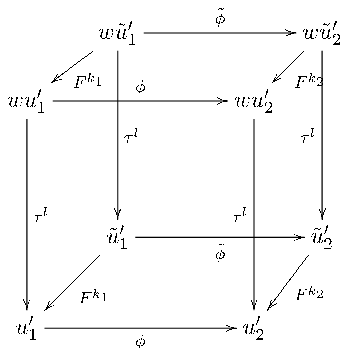}
\end{center}

\qed

\begin{cor}[Diameter bound implies finite]
\label{cor:diameter_bound_implies_finite}
For a fixed diameter bound $D$, there are only
finitely many isometry types of maps of shadows $F^k: S(\wtV) \to S(V)$ where  $\wtV$ is
a horizontal subset of diameter $< D$ and where 
$F^k:\Sigma(\wtV) \to \Sigma(V)$.  
\end{cor}

\pf Up to isomorphism, there are only finitely many isomorphism types 
of labelled graphs of the form $\Sigma(\widehat{\wtV}), \Sigma(\widehat{V})$, 
where $\wtV, V$ are horizontal subsets of diameter $\leq D$.   Hence, up to isomorphism, 
there are only finitely many  cellular, label-preserving maps $\Sigma(\widehat{\wtV})\to \Sigma(\widehat{V})$.
\qed

\subsection{Finiteness principles}

Let $v \in X^*$ and $r>0$ be an integer.  
If $B_{hor}(v,r)$ is a ball and $\Gamma = \Sigma(B_{hor}(v,r))$ is the induced subcomplex,
then  a component of $F^{-k}(\Gamma)$ need not be the induced subcomplex of a ball.  Instead, 
it may be the induced subcomplex of a union of overlapping balls.  The theorem below says that these components cannot be
too large, provided the level $|v|$ is large enough.  More precisely, it says that  at sufficiently high levels, for all iterates, components of preimages of 
induced subcomplexes of horizontal balls map
by uniformly bounded degree.  

\begin{thm}[Bounded degree]
\label{thm:bounded_degree}
Let $r>0$ be an integer.   Then there exists  a positive integer bound $C=C(r)$
and a level $n=n(r)$ such that $(\forall v, |v|>n)(\forall k>0)(\forall \tilde{v}
\in F^{-k}(v))$ if
$\Gamma = \Sigma(B_{hor}(v,r))$, $\widetilde{\Gamma}$ is a component of $F^{-k}(\Gamma)$,
and
$\wtV=F^{-k}(v)\intersect \widetilde{\Gamma}$, then 
\be
\item $\# \wtV < C$
\item $\diam_{hor}\wtV < D:=(2r+1)(C+1)$.
\eb
\end{thm} 

Using Theorem \ref{thm:bounded_degree}, we first derive the main result of this section.

\begin{thm}[Finiteness Principles]
\label{thm:finiteness_principles}
Fix an integer $r>0$.  Then there are only finitely many isometry types of maps of the form 
\[ F^k: \widetilde{S} \to S \]
and of the form 
\[ F^k: \widetilde{U} \to U\]
where $S=S(v,r)$, 
$|v|>n(r)$, $\widetilde{S}$ is a connected component of $F^{-k}(S)$, and
$\wtU, U$ are their corresponding umbrae.  
\end{thm}

\pf By Corollary \ref{cor:diameter_bound_implies_finite} and Lemma
\ref{lemma:naturality_of_umbrae}, it suffices to show that the diameter of a component 
$\widetilde{\Gamma}$ of $F^{-k}(\Gamma)$ is uniformly bounded whenever $\Gamma =
\Sigma(B(v,r))$ and $|v|$ is sufficiently large.  The bound is implied by 
[Theorem \ref{thm:bounded_degree}, Bounded degree].  
\qed

We now establish Theorem \ref{thm:bounded_degree}.  
The proof   depends on a purely algebraic result about
contracting selfsimilar group actions.  Recall that $m(L) = \max\{ m(g) : ||g|| \leq L\}$ where $m(g)$ is as in \S 4.2.   

In the theorem below, $\Stab_G(v)$ is the stabilizer of $v$ in $G$, $B_G(1,L)$ is the ball of
radius
$L$ about the identity in $G$ with respect to the word metric $||\cdot||$, the group $\genby{\Stab(v)\intersect B_G(1,L)}$ is the subgroup of $G$
generated by their intersection, and $\genby{\Stab(v)\intersect B_G(1,L)}.vw$ is the orbit of the
vertex $vw$ under the action of this subgroup.

\begin{thm}[Orbits under vertex stabilizers]
\label{thm:orbits_under_vertex_stabilizers}
Let $G$ be a finitely generated group, equipped with a selfsimilar, contracting group action on $T(X)$.
Fix a finite ``good''
generating set $S$ for $G$ so that (i) $S$ is closed under restriction, and (ii)
$\mathcal{N} \subset S$.  Let $N=\# \NNN$.
Then for all $w \in X^*$, and 
for all $L \geq 1$, if $|v| \geq m((N+1)L)$, we have 
\[ \#\genby{\Stab_G(v)\intersect B_G(1,L)}.vw \leq 1+q+q^2 + \ldots + q^{N+1}\]
where $q=\#\Stab_G(v)\intersect B_G(1,L)$.
\end{thm}

\gap

{\noindent\bf Proof of Thm. \ref{thm:orbits_under_vertex_stabilizers}:} 
Let $\mathcal{H}_{v,L}=\Stab(v)
\intersect B_G(1,L)$, $w \in X^*$, and consider the one-complex $\mathcal{K}$ with vertex set equal to the set of elements $\tilde{v}$ in the orbit $\genby{\mathcal{H}_{v,L}}.vw$ and edges  from $\tilde{v}$ to $\tilde{v}^h$, $h \in \mathcal{H}_{v,L}$. 
(This is just the Schreier graph of the action of $\genby{\mathcal{H}_{v,L}}$ on the orbit of $vw$.)  
Note that by construction, $\mathcal{K}$ is connected.  

Suppose  $|v| \geq m((N+1)L)$ holds and that in the intrinsic graph metric of $\mathcal{K}$ 
there exists a geodesic edge-path $W$ with $N+1$ vertices.   
On the one hand,  the restriction to $W$ of the iterated left shift $\tau^{m(L(N+1))}|_W$ is injective, 
since the right-hand substrings of length $|w|$ are all distinct.  
On the other hand, by contraction and the definition of $m(\cdot)$, 
the vertices in the image $\tau^{m(L(N+1))}(W)$ are all joined to one another by edges labelled by elements 
of the nucleus.  We conclude that the graph-diameter of $\mathcal{K}$ is at most $N$.  

In the $1$-skeleton of the complex $\mathcal{K}$, the valence of each vertex is 
at most $q$. Hence the number of vertices in a sphere of radius $n$ is bounded by $q^n$.
Since $\diam\, \mathcal{K}$ is bounded by $N$, it follows that the total number of vertices in $\mathcal{K}$ is bounded by  
$$1 + q +  q^2 + \ldots + q^{N+1}$$ and the proof is complete.

\qed

{\noindent\bf Proof of Thm. \ref{thm:bounded_degree}:}  Conclusion 2 follows immediately from 1,
Lemma \ref{lemma:pullbacks_of_balls}, and the triangle inequality.  To prove 1, let $\tilde{v}_0
= vw_0$, $\tilde{v}_p =vw_p \in \wtV$ be arbitrary.  
By Lemma \ref{lemma:pullbacks_of_balls}(2) there exist $\tilde{v}_i = vw_i \in \wtV$, $i=1,
\ldots, p-1$ such that for $1\leq i \leq p$ we have $|\tilde{v}_{i-1}-\tilde{v}_i|_{hor} \leq
2r+1=L$.  This implies $\tilde{v}_p \in \genby{\Stab(v)\intersect B_G(1, L)}.vw_0$.  Since
$\tilde{v}_p\in \wtV$ is arbitrary we conclude that $\wtV \subset \genby{\Stab(v)\intersect
B_G(1, 2r+1)}.vw_0$ and so the cardinality of this set is bounded by $C(L)$, the
constant in Theorem \ref{thm:orbits_under_vertex_stabilizers}.   
\qed

%\input sec5.tex

%\chapter{Boundary dynamics}

\section{$\bdry F: \bdry \Sigma \to \bdry \Sigma$ is finite type}\label{chap:bdry}

As in the previous chapter, we assume we are given a faithful, selfsimilar, level-transitive, contracting, recurrent action of a finitely generated group $G$ on the set of words $X^*$ in a finite alphabet $X$, and are given a symmetric generating set $S$ which is closed under restrictions and which contains the nucleus $\NNN$ of the action.  Let $\Sigma$ denote the selfsimilarity complex associated to $G$ and $S$; it is a Gromov hyperbolic augmented tree.   

In this chapter, we apply the finiteness principle, Theorem \ref{thm:finiteness_principles}, to prove that the induced dynamics $\bdry F: \bdry \Sigma \to \bdry \Sigma$ is metrically of finite type with respect to a certain covering $\UUU_0$ when equipped with the visual metric on the boundary (Theorem \ref{thm:boundary_dynamics_finite_type}).  The covering $\UUU_0$ will consist of connected components of the boundary of umbrae associated to the shadows of horizontal balls of radius one at a sufficiently deep level.   

Recall that an infinite geodesic ray based at $o$ is necessarily vertical.

\subsection{Metrics on the boundary}
\label{secn:Gromov_metric}

\begin{defn}Let $(u_1, u_2)$ be an ordered pair of vertices of $\Sigma$, and let $R_1, R_2: [0,\infty) \to \Sigma$ be infinite geodesic rays starting from $o$.
\be
\item The {\em Gromov inner product} with respect to the basepoint $o$ is 
\[ \gromprod{u_1}{u_2}=\frac{1}{2}\left( |u_1|+|u_2|-|u_1-u_2|\right).\]

\item The {\em normal form geodesic inner product} is 
\[ \level{u_1}{u_2} = \min_\gamma \{|\gamma| \}\]
where the minimum is over the set of normal form geodesics joining $u_1$
to $u_2$. 

\item The
{\em divergence inner product} is 
\[ \divlevel{R_1}{R_2} = \max\{t : |R_1(t)-R_2(t)| \leq 1\}.\]
\eb
\end{defn}

The lemma below says that these three products are all comparable.

\begin{lemma}
\label{lemma:products_comparable}
\be
\item For all pairs $(u_1, u_2)$ of vertices, we have 
\[ \level{u_1}{u_2} - \frac{H_\Sigma}{2} \leq \gromprod{u_1}{u_2} \leq
\level{u_1}{u_2}.\]
\item
\be
\item $d_{haus}(R_1, R_2) < \infty \iff \divlevel{R_1}{R_2}= \infty$
\item if $\divlevel{R_1}{R_2}<\infty$, then $\forall
s,t>\divlevel{R_1}{R_2}$, 
\[ \divlevel{R_1}{R_2} \leq \level{R_1(s)}{R_2(t)} \leq
\divlevel{R_1}{R_2}+m(H_\Sigma).\]
\eb
\eb
\end{lemma}

\pf {\bf 1.}  If $u_1, u_2$ lie on the same geodesic ray through $o$ then
$\gromprod{u_1}{u_2} = \min\{|u_1|, |u_2|\} = \level{u_1}{u_2}$.  
Otherwise, let 
\[ \gamma = [u_1u_1'] \union [u_1'u_2'] \union [u_2'u_2]\]
be a normal form geodesic, where the middle segment is the horizontal
portion.  Then 
\[ \begin{array}{ccl}
\gromprod{u_1}{u_2} & = & \frac{1}{2} \left(|u_1|+|u_2| - ( |u_1|-|u_1'|
+ |u_1' - u_2'| + |u_2|-|u_2'|) \right) \\
\; & = & \level{u_1}{u_2} - \frac{1}{2}|u_1' - u_2'| 
\end{array}
\]
so
\[ \level{u_1}{u_2} - \frac{H_\Sigma}{2} \leq \gromprod{u_1}{u_2} \leq \level{u_1}{u_2}.\]

{\bf 2.}  (a) The implication $\Leftarrow$ is obvious.  
Note that for any rays $R_1$ and $R_2$ and for any $s,t \ge 0$, $|R_1(s)-R_2(t)|\ge |s-t|$ so that
if  $d_{haus}(R_1, R_2) \leq D$, then  for all $t$ we have $|R_1(t)-R_2(t)|\leq 2D$.
Hence for all $t$, $|R_1(t-m(2D))-R_2(t-m(2D))| \leq 1$.  
(b) The lower bound is obvious, since any geodesic joining two points at
horizontal distance one is necessarily horizontal.  The upper bound
follows easily.  
\qed

\begin{defn}[$\bdry \Sigma$]\label{def:bdrySigma}
The {\em limit space} $\bdry \Sigma$ of the action is the quotient space 
\[ \bdry \Sigma = \{\mbox{\rm geodesic rays }R:[0,\infty) \to \Sigma |
R(0)=o\}/\sim\]
where the {\em asymptotic equivalence relation} $\sim$ is defined by 
\[ R_1 \sim R_2 \iff d_{haus}(R_1, R_2) < \infty.\]
\end{defn}
Cf. \S 4.1. 
We will often denote by $\xi$ or by $R(\infty)$ the equivalence class of a geodesic ray 
$R$.
\gap

One may define another topology on $\partial\Sigma$ as follows. Given a geodesic ray $R_0$, a fundamental
system of neighborhoods of $R_0(\infty)$ is given by the sets $V(R_0,t)$ of points $R(\infty)$
which admit a representative with $\divlevel{R_0}{R}> t$.   This topology coincides with the quotient topology defined above.  

\gap

\noindent{\bf Remarks:}  In fact, in our case the following conditions are
equivalent:  
\be
\item $d_{haus}(R_1,R_2) < \infty $
\item $d_{haus}(R_1, R_2) \leq 8\delta$ where $\delta$ is the constant of
hyperbolicity (in the Gromov inner product definition).
\item $d_{haus}(R_1, R_2) \leq 1$.   
\eb
The equivalence of (1) and (2) holds for general $\delta$-hyperbolic
spaces; see \cite{ghys:delaharpe:groupes}.  The equivalence of (3) and
(1) is implied by Lemma \ref{lemma:products_comparable}.  By (3) and the definition of $\Sigma$, there are at most $\#S+1$ elements in any asymptotic equivalence class; with a bit more work one can improve this bound to $\#\NNN$, the cardinality of the nucleus.   

The Gromov, normal form geodesic, and divergence inner products of a
pair of points on $\bdry \Sigma$ are defined as follows:
\begin{defn}
Let $\xi_1, \xi_2 \in \bdry \Sigma$.  Then 
\[ \gromprod{\xi_1}{\xi_2} = \sup_{R_1\in \xi_1, R_2\in \xi_2}
\liminf_{s_i, t_j \to \infty} \gromprod{R_1(t_i)}{R_2(s_j)}\]
with analogous formulae used for $\level{\cdot}{\cdot}$ and
$\divlevel{\cdot}{\cdot}$.  
\end{defn}
The following fact is known; see \cite{ghys:delaharpe:groupes}.  
For all $R_1 \in \xi_1, R_2 \in \xi_2$, we have 
\[ \gromprod{\xi_1}{\xi_2} - 2\delta \leq \liminf_{s_i, t_j \to \infty}
\gromprod{R_1(t_i)}{R_2(s_j)} \leq \gromprod{\xi_1}{\xi_2}\,.\]

Thus, up to a universal additive constant, the choice of
representative ray and the choice of points (at sufficiently
high level) on these rays used in the computation is irrelevant. Moreover, by
Lemma
\ref{lemma:products_comparable}, there exists a constant $C_0$ such that for all $\xi_1, \xi_2$, 
\[ \left|  \max\{\gromprod{\xi_1}{\xi_2}, \level{\xi_1}{\xi_2},
\divlevel{\xi_1}{\xi_2}\} - \min \{\gromprod{\xi_1}{\xi_2}, \level{\xi_1}{\xi_2},
\divlevel{\xi_1}{\xi_2}\} \right|  < C_0\]
Therefore, any one of these quantities may
be approximated, to within $C_0$, by computing $\level{R_1(t)}{R_2(t)}$
for any $t$ which is sufficiently large.  

Tracing through the dependencies shows that one may in fact take this
additive constant to be $C_0 = 100(\delta + m(H_\Sigma))$.    

It is known that 

\be
\item $\gromprod{\xi_1}{\xi_2} = \gromprod{\xi_2}{\xi_1}$
\item $\gromprod{\xi_1}{\xi_2} = \infty \iff \xi_1 = \xi_2$
\item $\gromprod{\xi_1}{\xi_3} \geq \min\{\gromprod{\xi_1}{\xi_2},
\gromprod{\xi_2}{\xi_3}\} - \delta$.
\eb 
and we have similar statements for the other two products.  

For $\epsilon>0$ let $\varrho_\epsilon(\xi_1, \xi_2) =
\exp(-\epsilon\level{\xi_1}{\xi_2})$.  The function $\varrho_\epsilon$ does
not quite define a distance function, as the triangle inequality fails
slightly due in particular to the constant present in (3) above. But, by choosing
$\epsilon$ small enough and by using {\em
chains}, however, this does define a metric, a so-called {\em visual 
metric}, and the resulting distance function turns out to be comparable
to $\varrho_\epsilon$ to within universal mutliplicative constants which
tend to one as $\epsilon$ decreases to $0$.

By using the normal form geodesic inner product, we then have 

\begin{prop}[existence of visual metrics on $\Sigma$]
\label{prop:gromov_metric}
There exists a constant $\epsilon_0>0$ depending only on $\delta$
(equivalently, on $H_\Sigma$) such that for all
$\epsilon<\epsilon_0$, there exist a metric $d_\epsilon$ compatible with
the topology of $\partial\Sigma$ and a constant $C_\epsilon\ge 1$ (with
$C_\epsilon \to 1$ as $\epsilon \to 0$) with the following property.  For
all $\xi_1, \xi_2 \in \bdry \Sigma$, 
\[ \frac{1}{C_\epsilon} \leq
\frac{d_\epsilon(\xi_1,\xi_2)}{\varrho_\epsilon(\xi_1, \xi_2)}
\leq C_\epsilon.\]
\end{prop}

\subsection{Invariance properties of boundaries}
\label{secn:invariance_properties}

\begin{defn}[Starlike]
Let $W \subset \Sigma$.  We say that $W$ is {\em starlike} if, for any geodesic ray $R$, $R(t) \in W \implies R([t,\infty)) \subset W$.  Equivalently: for any vertex $v$, $v \in W$ implies $C_v \subset W$.  If $W$ is starlike, its {\em boundary}
$\bdry W$ is 
\[ \bdry W = \{\xi \in\bdry\Sigma: \exists R \in \xi, t >0\; \mbox{ with } \; R([t,\infty))
\subset W\}.\]   
\end{defn}

Cones and shadows are starlike by definition and umbrae are starlike by Lemma
\ref{lemma:properties_of_shadows}(1).   We will often use the following
facts without reference:

\begin{lemma}
\label{lemma:when_boundaries_intersect}
For $i=1,2$ let $V_i$ be horizontal sets of vertices at the same level, 
and let $U_i = U(V_i)$ be their corresponding umbrae.

If $\cl{\bdry U_1} \intersect \cl{\bdry U_2} \neq \emptyset$, then $|V_1 - V_2|_{hor} \leq 1$.
\end{lemma}

\pf
Let $\xi$ denote an element of the common intersection and suppose 
$(\xi_i^n)$ are sequences in $\partial U_i$ such that  
$\xi_i^n \to\xi$ as $n \to \infty$, $i=1,2$.  Then for all $n$, there exist rays $R_i^n$
representing $\xi_i^n$ such that $R_i^n \intersect V_i \neq \emptyset$ and 
$ \divlevel{R_1^n}{R_2^n} \to \infty.$
Hence there exists $n$ with $\divlevel{R_1^n}{R_2^n}$ greater than the common level
of $V_1, V_2$.  It follows that $|V_1-V_2|_{hor} \leq 1$.  
\qed

\begin{lemma}[Naturality of boundaries]
\label{lemma:naturality_of_boundaries}
Suppose $\widetilde{W}, W$ are starlike.
If $F^k(\widetilde{W}) = W$, then 
$\bdry F^k(\bdry \widetilde{W}) = \bdry W$.
\end{lemma}

\pf The containment $\subset$ is trivial.  To prove the other inclusion,
suppose $R(\infty) \in \bdry W$ and $w=R(t) \in W$.  Let $\gamma=[w, R(\infty)]$; it is a vertical geodesic ray.  
Let $\widetilde{R}$ denote the lift of this segment based at $\tilde{w} \in\widetilde{W}$.  
Since $\widetilde{W}$ is starlike, $\widetilde{R}(\infty) \in \widetilde{W}$.  
Clearly $\bdry F^k(\widetilde{R}(\infty)) = R(\infty)$.
\qed

Recall from \S \ref{secn:finite_type} the definition of quasisimilarity.  The proof of the lemma below 
is a direct consequence of the definitions. 

\begin{lemma}[Induced quasisimilarities]
\label{lemma:boundary:quasisimilarity} 
Let $W_i$, be shadows or umbrae of horizontal sets $V_i$ at levels $n_i$, $i=1,2$.    
Let $\phi: W_1 \to W_2$ be an isometry which preserves the relation
of being at the same level.   
Then $\phi$ induces a homeomorphism $\bdry \phi:\bdry W_1 \to
\bdry W_2$ which is a
$(C^2_\epsilon,\exp(-\epsilon(n_2-n_1))$-quasisimilarity, where
$C_\epsilon$ is the constant as in the definition of the visual metric; see \S
\ref{secn:Gromov_metric}.
\end{lemma}

\subsection{Metric estimates for boundaries of umbrae}

\noindent{\bf Notation.}  Let $R$ be a geodesic ray and $t\in\N$.  We denote by 
\[ S(R,t) = S(B_{hor}(R(t),1)) \]
the shadow corresponding to the closed horizontal ball of radius one centered at
$R(t)$, and by $U(R,t)$ the corresponding umbra.  
Clearly, the shadow $S(R,t)$ contains its umbra.  The following lemma provides a partial converse.

\begin{lemma}[Umbrae contain shadows]
\label{lemma:umbrae_contain_shadows}
Let $R$ be a geodesic ray.  
\be
\item For all $s>t$, $R(s) \in U(R,t)$. 
\item Let $V$ be any horizontal vertex set.  If $R(t) \in U(V)$ then
$S(R,t+c) \subset U(V)$ where $c=m(H_\Sigma)+3$.
\item In particular:  $U(R,t) \supset S(R,t+c)$.
\eb
\end{lemma}

\pf {\bf 1.}  Suppose $v$ is any vertex with $|R(s)-v|_{hor}=1$.  
Then $| R(t)-v^{[s-t]}|_{hor} \leq 1$, or in other words  
$v^{[s-t]} \in B_{hor}(R(t), 1)$; hence $v \in S(B_{hor}(R(t), 1)) = S(R,t)$.  
Since $v$ is arbitrary, we have then that $|R(s)-S(R(t))^c| > 1$.  
By the definition of umbra, $R(s) \in U(R,t)$.

{\bf 2.} Let $u=R(t)$, $v=R(t+c)$, $w \in S(R,t+c)$ and $w' = [ow] \intersect B_{hor}(v,1)$. 
Note that $|v-w'|\leq 1$.  By Lemma
[\ref{lemma:unit_speed_penetration}, Unit speed penetration], 
\[ |v-{^cS(V)}| \geq t+c-(t+m(H_\Sigma)) = 3.\]
Hence $|w'-{^cS(V)}| \geq \left| |w'-v|-|v-{^cS(V)}| \right| \geq 2$ and so $B(w',1)
\subset S(V)$.  Hence $w' \in U(V)$ and so $w \in U(V)$ since umbrae are
starlike.  

{\bf 3.}  This is an immediate consequence of (1) and (2).  
\qed

In the remainder of this section, 
\bi
\item $C_0$ denotes the universal additive constant as in \S \ref{secn:Gromov_metric};
\item $\epsilon_0$ is the positive constant as in Proposition \ref{prop:gromov_metric};
\item $\epsilon < \epsilon_0$ is fixed;
\item $C_\epsilon, d_\epsilon$ denote the constant and metric in Proposition
\ref{prop:gromov_metric}. 
\ib

\begin{lemma}[Diameter estimates]
\label{lemma:diameter_estimates}
For all $t>C_0+4m(H_\Sigma)$, and for all geodesic rays $R$, we have the
following diameter estimates:
\be
\item  
\[ B_\epsilon\left(R(\infty), \frac{1}{C_1}\exp(-\epsilon t)\right) \subset
\bdry S(R,t)
\subset B_\epsilon(R(\infty), C_1 \exp(-\epsilon t)) \]
where $C_1 = 2C_\epsilon \exp(\epsilon (C_0+m(H_\Sigma)))$  

\item  \[ B_\epsilon\left(R(\infty), \frac{1}{C_2}\exp(-\epsilon t)\right) \subset
\bdry U(R,t)
\subset B_\epsilon(R(\infty), C_2 \exp(-\epsilon t)) \]
where $C_2 = 2C_\epsilon \exp(\epsilon (C_0+4m(H_\Sigma)))$ 

\item  
\[ \frac{1}{C_3}\exp(-\epsilon t) <  \diam_\epsilon \bdry C_{R(t)} < C_3 \exp(-\epsilon t)\]
where $C_3 = C_\epsilon\exp(\epsilon(C_0 + \lceil \log_d B \rceil + 1))$, $d=\#X$, and $B \leq \#S+1$ is the maximum cardinality
of a horizontal unit ball (here, $C_{R(t)}$ is the cone above the vertex $R(t)$).  

\item For all $k \in \N$, for all $r <(1/C_1)\exp(-\epsilon k)$, and for all $\tilde{\xi}$, if $\xi=\bdry
F^k(\tilde{\xi})$ then 
\[ B_\epsilon\left(\xi, \frac{1}{C_4}\exp(k\epsilon)r\right) \subset \bdry F^k
(B_\epsilon(\tilde{\xi},r))\subset B_\epsilon(\xi, C_4\exp(k\epsilon)r).\]
where $C_4=C_1^2$.  
\eb

\end{lemma}

\pf {\bf 1.} We prove the upper bound first.  Let $u_1, u_2 \in
S(R,t)$ and suppose $u_i' \in [ou_i] \intersect B_{hor}(R(t), 1)$,
$i=1,2$.  Then $|u_1'-u_2'| \leq 2$ which implies that
$\level{u_1}{u_2} \geq t-1$.  Since $u_1, u_2$ are arbitrary elements of
$S(R,t)$ it follows that $d_\epsilon(\xi_1, \xi_2) \leq C_\epsilon 
\exp(-\epsilon(t-1))< C_1\exp(-\epsilon t)$ for all $\xi_1, \xi_2 \in
\bdry S(R,t)$.  

To prove the lower bound, let $\xi = R(\infty)$, and suppose $\xi' \in
\bdry\Sigma$ satisfies $\level{\xi}{\xi'} > t+m(H_\Sigma)+C_0$. 
Let $R'$ represent $\xi'$.  Then there exists infinitely many integers $s$ such
that $\level{R(s)}{R'(s)}>t+m(H_\Sigma)$.  This implies that $R'\intersect
B(R(t),1) \neq \emptyset$, i.e. that $\xi' \in
\bdry S(R,t)$ and the estimate follows.  

{\bf 2.}  By Lemma
[\ref{lemma:umbrae_contain_shadows}, Umbrae contain shadows, (2)]
we have $U(R,t) \supset S(R,t+c)$.  This together with part (1) above
yields the estimates. 

{\bf 3.}  The upper bound is clear.  To find the lower bound, let $k= 
\lceil \log_d(B) \rceil+1$ where $B$ is the maximum cardinality of a horizontal unit ball. 
Note that there are $d^k$ elements of the cone $C_{R(t)}$ at level $t+k$.  Thus some point
$u$ of $C_{R(t)}$ at level $t+k$ lies at horizontal distance greater than $1$ from $R(t+k)$.  Let $R'$ be a ray
through $u$.  Then  $\divlevel{R}{R'} \in [t, t +k]$ and the claim follows easily.  

{\bf 4.}  This follows immediately from the fact that shadows map to
shadows:  $F^k(S(R,t)) = S(F^k (R), t-k)$,  and estimate (1) above.  
The condition on $r$ ensures that we may choose $t>k$.
\qed

\begin{cor}[Umbrae are open]
\label{cor:umbras_open}
Let $V$ be any vertex set and $U$ its umbra.  Then $\bdry U$ is an
open subset.
\end{cor}

\pf $U$ is nonempty if and only if $\bdry U$  is nonempty, since
umbrae are starlike.  Suppose $\xi \in \bdry U$.  By definition (see \S\S \ref{subsecn:umbrae} and \ref{secn:invariance_properties}), this means there exists a ray $R$ and a level $t$ with $R(\infty)=\xi$ and $R(t) \in U$.  Hence by Lemma [\ref{lemma:umbrae_contain_shadows},
Umbrae contain shadows] $S(R,t+c) \subset U$.  Hence
$B_\epsilon(R(\infty),
\frac{1}{C_1}
\exp(-\epsilon(t+c))) \subset \bdry S(R,t+c) \subset \bdry
U$ by Lemma \ref{lemma:diameter_estimates}(1).  

\qed

\subsection{$\mathbf{\bdry F}$  is a branched covering}
\label{secn:is_br_cov}

Recall from Section \ref{secn:fbc} the definition of a finite branched covering.  
The main result of this section is 

\begin{thm}[$\mathbf{\bdry F}$ is a branched covering]
\label{thm:is_br_cov} The map $\bdry F: \bdry \Sigma \to
\bdry \Sigma$  is a finite branched covering map of degree $d=\#X$.
\end{thm}

\pf 
Theorem \ref{thm:props_of_Sigma}(4) implies that $\bdry F$ is continuous, open and finite-to-one.  Lemma \ref{lemma:naturality_of_boundaries} and Corollary \ref{cor:umbras_open} imply that $\bdry F$ is open.
The map $\bdry F$ is also closed since $\bdry\Sigma$ is compact and Hausdorff. 
Therefore, since $\bdry\Sigma$ is locally connected (Theorem \ref{thm:props_of_Sigma}),
it is enough to show
\[ \forall \xi \in \bdry \Sigma, \ \   \sum_{\txi \in \bdry F^{-1}(\xi)}\deg(\bdry F, \txi) = d.\]

In the sequel, we will denote by $\xi, \txi, \zeta$, etc. elements of $\bdry \Sigma$.  Recall the definition of local degree:
\[ \deg(\bdry F, \txi) = \inf_{U \ni \txi}\{\#\bdry F^{-1}(\bdry F(\zeta)): \zeta \in U\}\]
where the infimum is over all open $U$ containing $\txi$.

The proof will proceed by first interpreting the degree at a point on $\bdry \Sigma$
in terms of the map $F$ on $\Sigma$ itself, and then using the fact that $F$ is a covering map of degree $d$.

To this end, recall that a point $\xi \in \bdry \Sigma$ is an equivalence class of rays emanating from $o$.  Also, recall that balls are by definition subsets of vertices, and that for any subset $V$ of vertices, $\Sigma(V)$ denotes the induced subcomplex of $\Sigma$ containing $V$.  

Let us denote by 
\[ \Sigma(\xi)_n = \Sigma\left(\union_{R \in \xi} B_{hor}(R(n), 1)\right).\]
Observe that for all $n \in \N$ and all $\xi \in \bdry \Sigma$, 
\begin{equation}
\label{eqn:brcover:dropdown}
\Sigma(\xi)_n^{[-1]} \subset \Sigma(\xi)_{n-1}
\end{equation}
since $\Sigma$ is an augmented tree.
Recall that $F$ maps equivalence classes of rays surjectively onto equivalence classes of rays, and maps horizontal balls of radius one surjectively onto horizontal balls of radius one (cf. Lemma \ref{lemma:balls_map_to_balls}).  It follows that whenever $\bdry F(\txi) = \xi$, the restriction 
\[ F: \Sigma(\txi)_{n+1} \to  \Sigma(\xi)_n\]
is surjective for all $n \geq 0$.   
%
%Suppose  $v \in \Sigma(\xi)_n$ is any vertex, and consider the set 
%\[ F^{-1}(v) \intersect \Sigma(\txi)_{n+1}. \]
%Denote by 
%\[ \delta(\txi, v)=\#F^{-1}(v) \intersect \Sigma(\txi)_{n+1}\]
%and by 
%\[ \delta(\txi)_{n+1}=\max_{v \in \Sigma(\xi)_n}\delta(\txi, v).\]
%Fix $n$, for convenience set $k=\delta(\txi)_n$, and choose a vertex $v$ for which $\delta(\txi, v)=k$.  
%Since $F$ acts as the right shift, we may write 
%\[ F^{-1}(v) \intersect \Sigma(\txi)_{n+1}=\{vx_1, vx_2, \ldots, vx_k\}.\]
% By (\ref{eqn:brcover:dropdown}), $v^{[-1]}x_i \in \Sigma(\txi)_n$, $i=1, \ldots, k$, so we have that   
%\[ v^{[-1]}x_i \in F^{-1}(v^{[-1]}) \intersect \Sigma(\txi)_n , \, \;\; i=1, \ldots k.\]
%Hence  
%\[ \delta(\txi)_n \geq \delta(\txi, v^{[-1]}) \geq k=\delta(\txi, v) = \delta(\txi)_{n+1}.\]
%Hence as $n \to \infty$, the quantity $\delta(\txi)$ defined by 
%\[ \delta(\txi):=\lim_{n \to \infty} \delta(\txi)_n \]
%exists.  
%Let $\txi \in \bdry \Sigma$, set $\xi = \bdry F(\txi)$,  and let $n \in \N$.  
Then  the quantity 
\[ \delta(\txi)_{n+1} = \max_{v \in \Sigma(\xi)_n} \# \{ x\ : \ vx \in \Sigma(\txi)_{n+1}\}
=\max \{ \# F^{-1}(F(\tilde{v})) \intersect \Sigma(\txi)_{n+1} : \tilde{v} \in  \Sigma(\txi)_{n+1}\} \]
is the maximum cardinality of a fiber of this restriction. 
It follows that
\[ \delta(\txi)_{n+1} \leq \delta(\txi)_n\]
holds so that the limit 
\[ \delta(\txi)=\lim_{n \to \infty} \delta(\txi)_n\]
exists.  We are going to show that in fact $\delta(\txi)=\deg(\bdry F, \txi)$.
\gap

Although $F(\Sigma(\txi)_{n+1}) = \Sigma(\xi)_n$, the restriction $F|_{\Sigma(\txi)_{n+1}}$ need not be proper.  Let $\wtG(\txi)_{n+1}$ denote the unique component of $F^{-1}(\Sigma(\xi)_n)$ which intersects $\Sigma(\txi)_{n+1}$; a priori, it is larger than $\Sigma(\txi)_{n+1}$.  Its diameter, however, is uniformly bounded by a constant $D$ depending only on $d$ and the number of generators in $S$.  Let $m(D)$ denote the "magic level" such that for any pair of vertices $v_1, v_2$ at the same level, $|v_1-v_2|_{hor} \leq D$ implies $\left|v_1^{[-m(D)]} - v_2^{[-m(D)]}\right|_{hor} \leq 1$.  
\gap

\noindent{\bf Claim:}  {\em For all $n$ sufficiently large, the map 
\[ F: \wtG(\txi)_{n+1} \to \Sigma(\xi)_n\]
is a covering map of degree $\delta(\txi)$.}
\gap

\noindent{\bf Proof of Claim:} 
The definition of $\delta(\txi)$ shows that for all $n$ sufficiently large, the degree of the covering map $F: \wtG(\txi)_{n+1} \to \Sigma(\xi)_n$ is at least $\delta(\txi)$.  
We now establish the upper bound.  

Suppose  $n$ is large and $v \in \Sigma(\xi)_n$ is any vertex.  Write $v=uw$ where $w=v^{[-m(D)]}$, and consider now the set 
\[ F^{-1}(v) \intersect \wtG(\txi)_{n+1}=\{uwx_1, \ldots, uwx_k\}. \]
By definition, there exists a ray $R \in \xi$ such that $|R(n)-v| \leq 1$.
Choose a lift $\wtR$ of $R$ under $F$ representing $\txi$ such that $F(\wtR) = R$.  Then by definition 
\[ |\wtR(n+1)-uwx_i|_{hor} \leq D, \ i=1, \ldots, k.\]
By the definition of $m(D)$ we have 
\[ |\wtR(n+1-m(D))-wx_i|_{hor} \leq 1, \ i=1, \ldots, k\]
which implies that $wx_i \in \Sigma(\txi)_{n+1-m(D)}$ and hence that 
\[ k \leq \delta(\txi, v^{[-m(D)]})_{n+1-m(D)} \leq \delta(\txi)_{n+1-m(D)}=\delta(\txi)\]
provided that $n$ is sufficiently large.    
\qedspecial{Claim.}

Since $F$ is a covering map of degree $d$, the previous claim implies immediately that for all $\xi \in \bdry \Sigma$, 
\[
\sum_{\txi \in \bdry F^{-1}(\xi)} \delta(\txi) = d.
\]
The following lemma completes the proof of Theorem \ref{thm:is_br_cov}.  

\begin{lemma}[Interpretation of local degree]
\label{lemma:interpretation_of_degree}
For all $\txi \in \bdry \Sigma$, 
\[ \deg(\bdry F, \txi) = \delta(\txi).\]
\end{lemma}

\pf We first establish $\leq$.  Suppose $\bdry F (\txi) = \xi$ 
and $k$ is the local degree at $\txi$.  For $n \in \N$ let 
\[ \wtU_n = \{\tzeta \in\bdry\Sigma: \ \gromprod{\txi}{\tzeta} > n+C_0\}\]
and
\[ U_n = \bdry F(\wtU_n)\]
where $C_0$ is as in \S 6.1. 
Then $\{\wtU_n\}$ is a basis of neighborhoods of $\txi$.  

Then by the definition of local degree, for all large $n$, we have  
\[ k = \deg(\bdry F, \txi) = \max_{\zeta \in U_n}
\{\# \bdry F^{-1}(\zeta) \intersect \wtU_n\}.\]
Let us fix such an $n$ and consider a point $\zeta \in U_n$ which realizes the above supremum.
Let $R \in \xi, S \in \zeta, \wtR \in \txi$ be representing rays.  
Then by the definition of $\wtU_n$, there exist preimages $\wtS_i$, $i=1, \ldots, k$, of $S$ under $F$ 
such that for each such $i$, the divergence level satisfies $\divlevel{\wtR}{\wtS_i} > n$.    
Note that the rays $\wtS_i$ are all distinct, since they are preimages 
of a ray.  By the definition of divergence level, for all such $i$, 
\[ |\wtR(n)-\wtS_i(n)|_{hor} \leq 1.\]
Hence 
\[ \wtS_i(n) \in \Sigma(\txi)_n\]
and so since each $\wtS_i(n)$ maps to the same point $F(\wtS_i(n)) = S(n-1)$ we have 
\[ k \leq \delta(\txi)_{n} = \delta(\txi)\]
provided that $n$ is sufficiently large.
\gap

To establish the other bound, we make use of the following claim.  
The idea is that if $f: X \to Y$ is an fbc, then the image of the branch locus in $Y$ is nowhere dense, 
so that any $y \in Y$ is a limit of points each having the maximal number $d=\deg(f)$ of distinct preimages.    This is what we are going to show by means of the claim below.  Roughly, here is the idea of the proof.    Suppose instead of $\bdry F$ we consider a subhyperbolic rational map.  Take a very small ball $B$ which intersects the Julia set.  Take a sequence of inverse branches of this ball which {\em realizes the maximum  of the degree over all inverse branches} to obtain an iterated preimage $\widetilde{B}$ of $B$.   Then all further iterated preimages of $\widetilde{B}$ must be unramified, and the union of all of these further preimages is dense in the Julia set.  
\gap

\noindent{\bf Claim:}  {\em There are a universal level $n_0$ and a constant $M$ such that for any
vertex $\tilde{v}$ with $|\tilde{v}|\geq n_0$,
there exists a vertex $\tilde{w}$ with $|\tilde{v}|=|\tilde{w}|$
such that (i) $|\tilde{v}-\tilde{w}|_{hor} \leq M$, and (ii) for all $k \in \N$,
any pair of distinct preimages of
$\tilde{w}$ under $F^{-k}$
are at least two (horizontal) units apart.}
\gap

\noindent{\bf Proof of Claim:}  Let $r=1$ and let $n=n(r)$ be the constant given by Theorem \ref{thm:bounded_degree}.
Choose a vertex $u \in X^n$ at level $n$  arbitrarily, and 
set $\Gamma_0 = \Sigma(B_{hor}(u,1))$.

By Theorem \ref{thm:bounded_degree},
\[ p_0 = \sup\{ \deg(F^k: \wtG_0 \to \Gamma_0) : k \in \N, \wtG_0 \mbox{ is a component of }\; F^{-k}(\Gamma_0)\}\]
is finite.  Suppose $F^{k_0}: \wtG_0 \to \Gamma_0$ realizes this supremum.  For convenience set $\Gamma = \wtG_0$.  Then for all $k \in \N$ and all components $\wtG$ of $F^{-k}(\Gamma)$,
\[ \deg(F^k: \wtG \to \Gamma) = 1.\]
Let $n_0=n+k_0$ be the level of $\Gamma$.
Since by assumption the action is level-transitive, each horizontal subcomplex is connected.  
Hence the quantity $M$, defined as the maximum horizontal distance of a vertex at level $n_0$ to a vertex in $\Gamma$, exists.

Given now any vertex $\tilde{v}$ at level $n>n_0$, put $l=n-n_0$, and let $v=F^l(\tilde{v})$.
By the definition of $M$, there exists a vertex $w \in \Gamma$ and a horizontal 
edge-path $\gamma$ of length at most $M$ joining $v$ to $w$.  Let $\tilde{w}$ be the endpoint of the lift of $\gamma$ under $F^l$ based at $\tilde{v}$.  Then $|\tilde{v}-\tilde{w}|_{hor} \leq M$.

By construction, for all  $j > 0$, each connected component $\wtG$ of the preimage of $\Gamma$
under $F^j$ maps to $\Gamma$ by degree one.  Hence there are $d^j$ such preimages.
By  Lemma \ref{lemma:properties_of_shadows}, for a fixed $j$, any two such preimages $\wtG^i$, $i=1,2$ are at least two units apart.  In particular, this holds for $j=k+l$.  We conclude
that each of the $d^k$ inverse images of $\tilde{w}$ under $F^{-k}$
are at least two apart, and the Claim follows.
\qedspecial{Claim}

Now suppose that $k=\delta(\txi)$.  Then for all large $n$ we have 
\[ k = \delta(\txi)_{n+1}.\]
Fix $n$ large, let $v \in \Sigma(\xi)_n$, and suppose $\delta(\txi, v)=\delta(\txi)_{n+1}$, so that 
\[ F^{-1}(v) \intersect \Sigma(\txi)_{n+1} = \{\tv_1, \ldots, \tv_k\}.\]
Apply the Claim (with $\tv=v$ in the hypothesis) to obtain a vertex $w$ (called $\tw$ in the conclusion) for which $|v-w|_{hor} \leq M$, and for which all of the iterated preimages of $w$ at a given level are at least two horizontal units apart.  Let $\gamma$ be a 
horizontal edge-path joining $v$ and $w$ and let $\tilde{w}_i$ be the unique preimage 
of $w$ obtained by lifting $\gamma$ under $F$ based at $vx_i = \tilde{v}_i$.  
Then the $\tilde{w}_i$'s are all distinct.  
Let $S$ be a ray through $w$ and let $\wtS_i$ be the lifts of $S$ through $\tilde{w}_i$.  
Since the $\tilde{w}_i$'s are at least two apart, the rays $\wtS_i$ are in distinct 
equivalence classes $\tzeta_i$, (Lemma \ref{lemma:at_least_two_apart}), each of which maps to 
$\zeta = [S]$.  Since this occurs for all $n$ large, we have 
\[ \deg(\bdry F, \txi) \geq k=\delta(\txi).\]

\qedspecial{Lemma}

\qedspecial{Theorem \ref{thm:is_br_cov}}

\subsection{Dynamics on $\mathbf{\bdry \Sigma}$}
\label{secn:UUU_n}

In the next two sections, we prove that the dynamics on $\bdry \Sigma$ is
of finite type with respect to the visual metric.  
In this section, we define the family $\UUU_n$, $n=0,1,2,\ldots$ of open
covers of $\bdry \Sigma$, and collect the necessary finiteness results.  
In the next section, we show that the
dynamics on the boundary is of finite type.  

The construction of the family $\UUU_n, n \in \N$, is somewhat technical, so we first give the general idea.  Recall the definition of the {\em cone} $C_v$ associated to a vertex $v \in X$ (\S \ref{secn:cones_etc}).   Let us pretend for the moment that the boundary at infinity of any cone is both open and connected.  We would like to take $\UUU_0$ to be the set $\{\bdry C_v : |v|=n_0\}$ of boundaries at infinity of  the set of cones at some fixed level $n_0$.  Suppose $C$ is such a cone, $\widetilde{C}$ a preimage of $C$ under some iterate $F^k$, and $p=\deg(F^k: \widetilde{C} \to C)$.  Since probably $\widetilde{C}$ is not a cone but rather a union of cones, we may lose control of the degree $p$ .  However, as the level $n_0$ of  the vertex defining $C$ increases, the cone $C$ gets smaller.  Eventually, it is small enough so that the degree $p$ is uniformly bounded independent of $k$.    To make this rigorous, we use the finiteness principles, which assert that up to isometry, there are  only finitely many local models for the map $F^k: \widetilde{C} \to C$.

Unfortunately, a priori we do not know if the boundaries of cones are open and connected.  So, 
to make the above heuristic argument precise:
\bi
\item we work  with shadows and  umbrae, to get open sets of the boundary; 
\item we take a finite covering by {\em connected components} of umbrae at some fixed level $n_0$ sufficiently large so that the finiteness principles apply; 
\item we add basepoints to these connected components to aid in indexing preimages.
\ib
\gap

\noindent{\bf Remark:}  There are known conditions which imply that the boundaries of cones are connected and equal to the closure of their interiors; see \cite[\S 3.3.3]{nekrashevych:book:selfsimilar}.  

\gap

\noindent{\bf Construction of $\UUU_n$.} 
Take $r=1$ and let $n_0 = n(r)$ as in the statement of Theorem
[\ref{thm:finiteness_principles}, Finiteness Principles]. Let 
\[ \umb_0 \subset \{(U,R) | U = U(R,n_0)\} \]
be a {\em finite} set such that 
\[ \bigcup_{(U,R) \in \umb_0} \mho(U,R) \]
covers $\bdry \Sigma$, where $\mho(U,R)$ denotes the connected
component of $\bdry U$ containing $R(\infty)$.  
This set exists since umbrae are open (Corollary \ref{cor:umbras_open}) and $\bdry \Sigma$ is locally connected 
(Theorem \ref{thm:props_of_Sigma}).

Let 
\[ \umb_n = \{(\wtU, \wtR) | F^n:(\wtU,\wtR) \to (U,R) \in
\umb_0\}\]
where $\wtU \subset \Sigma$ is a connected component of
$F^{-n}(U)$, $F^n(\wtR) = R$, and $\wtR(t) \in \wtU$ for
all $t >> 0$.

Given $(U,R) \in \umb_n$ we let $\mho(U,R)$ denote the connected
component of $\bdry U$ which contains $R(\infty)$, equipped with
the basepoint $R(\infty)$.  Finally, we set 
\[ \UUU_n = \{\mho(U,R) | (U,R) \in \umb_n\}.\]
We emphasize that the elements of $\UUU_n$ are sets equipped with basepoints.  While this leads to redundancy which could be avoided, it is convenient since it aids in the bookkeeping of preimages that follows.   By abuse of notation, the index $n$ will here be called the {\em level} of
$(U,R) \in \umb_n$, even though $|o - U| = n+n_0$.  
Since $\bdry F$ is a finite branched covering (Theorem \ref{secn:is_br_cov}), we  have 
\begin{equation}
\label{eqn:boundary:alternative}
\UUU_{n+1} = \bdry F^{-1}\UUU_n
\end{equation}
in the sense that $\UUU_{n+1}$ is the set of all pairs
$(\widetilde{\mho}, \tilde{\xi})$ which map under $\bdry F$
to a pair $(\mho, \xi) \in \UUU_n$; 
note that the same underlying set $\widetilde{\mho}$ may arise more than once when equipped with different basepoints.

This completes the definition of the family $\UUU_n, n \in \N$.  

\gap

\noindent{\bf Application of finiteness principles.}
By Theorem [\ref{thm:finiteness_principles}, Finiteness
Principles] and the choice of $n_0$, there are only finitely many isometry classes
of maps of the form 
\[ F^k: \wtU \to U\]
where $(U,R) \in \umb_n$.  The elements of $\umb_n$ are sets with basepoints.   The
supremum of the degrees     
\[ p = \sup_k \sup_{(\wtU, \wtR) \in \umb_{n+k}}\deg(F^k|\wtU \to U) <
\infty\] 
is finite, and $\umb_0$ is finite, so we conclude:\\
\gap

{\em There are only finitely many isometry
classes of maps of pairs 
\[ F^k: (\wtU, \wtR) \to (U,R)\]
where $(\wtU, \wtR) \in \umb_{n+k}$ and $(U,R) \in \umb_n$.}
\gap

\subsection{Boundary dynamics is metrically finite type}

In this section, we prove that the boundary dynamics on $\Sigma$ is of finite type, 
hence is cxc by Theorem \ref{thm:finite_type_implies_cxc}. 
Actually, we will show the following slightly stronger statement in which the control
of diameters is more precise:  

\begin{thm}[Boundary dynamics is of finite type]
\label{thm:boundary_dynamics_finite_type}
The dynamics on $\bdry\Sigma$ is metrically finite type, hence is cxc.  

More precisely, we have the following.   There exists a family of open covers $\UUU_n$,
$n=0,1,2,\ldots$, constants $C,
\lambda > 1$ and a finite set of
$\MMM$ of {\em pointed model maps} 
\[ g_m: (\wtV_m,v_m) \to (V_m,v_m), m \in \MMM\]
where $\wtV_m, V_m$ are connected metric spaces, $\tilde{v}_m \in
\wtV_m, v_m \in V_m$, with the following property.  If 
\[ \UUU_{n+k}\ni (\wtmho, \tilde{\xi}) \mapsto (\mho, \xi) \in \UUU_n\]
under $\bdry F^k$, then there exist homeomorphisms 
$\widetilde{\psi}: \wtmho \to \wtV_m$ and $\psi: \mho \to V_m$ depending only on $\wtmho,
\mho$ respectively such that 
%commutative diagram 
\[ 
\begin{array}{ccc}
(\wtmho,\tilde{\xi}) & \stackrel{\widetilde{\psi}}{\longrightarrow} &
(\wtV_m,\tilde{v}_m)
\\
\bdry F^k \downarrow & \; & \downarrow  g_m
\\ (\mho,\xi) & \stackrel{\psi}{\longrightarrow}
&(V_m,v_m) \\
\end{array}
\] 
commutes, and such that $\widetilde{\psi}, \psi$ are respectively
$(C,\lambda^{n+k})$ and $(C,\lambda^n)$-quasisimilarities.
\end{thm}

Like the metric constructed in Section \ref{secn:top_finite_type} for topologically finite type maps, the metric $d_\varepsilon$ above has the property that $\diam U \asymp \exp(-\varepsilon n)$ where the constants are independent of $U$ and $n$.

\gap

\pf First, axiom [Expans] holds immediately from Lemma \ref{lemma:diameter_estimates}, 
since any open set on the boundary contains the boundary of a cone.  
Since any cone eventually maps onto all of $\Sigma$, its boundary eventually maps onto all of $\bdry \Sigma$, 
so axiom [Irred] holds as well.  

We now establish the existence of a dynatlas.  
Let $\MMM$ denote the set of isometry classes of maps of pairs
$F^k: (\wtU, \wtR) \to  (U,R)$ where $(\wtU, \wtR) \in \umb_{n+k}$
and $(U,R) \in \umb_n$.  For each $m \in \MMM$ choose a representative 
\[ F^{k_m}: (\wtU_m, \wtR_m) \to (U_m,R_m);\]
one may choose this representative to have minimal level if desired (it is more convenient
here not to normalize the sets to have diameter $1$).  
Let  $(\wtV_m,\tilde{v}_m) = (\mho(\wtU_m,\wtR_m),\wtR_m(\infty))$, 
$(V_m,v_m) =
(\mho(U_m,R_m),R_m(\infty))$, and
$g_m =
\bdry F^{k_m}|_{\wtV_m}$.  

Suppose now $\UUU_{n+k} \ni (\wtmho, \tilde{\xi}) \mapsto (\mho,\xi) \in
\UUU_n$ is induced by $\bdry F^k: (\wtU,\wtR) \to (U,R)$
where $\wtmho = \mho(\wtU,\wtR)$ and $\mho = \mho(U,R)$.  Then there
exists $m \in \MMM$ and isometries $\wtphi: \wtU \to \wtU_m$, $\phi: U
\to U_m$ such that the diagram below commutes:
\[ 
\begin{array}{ccc}
(\wtU,\wtR) & \stackrel{\widetilde{\phi}}{\longrightarrow} &
(\wtU_m,\wtR_m)
\\
F^k \downarrow & \; & \downarrow  F^{k_m}
\\ (U,R) & \stackrel{\phi}{\longrightarrow}
&(U_m,R_m) \\
\end{array}
\] 
Suppose $U_m$ has level $n_m$.  Then $\wtU$ has level $n_m + k_m$ while
the levels of $\wtU, U$ are $n+k$ and $n$, respectively.  By Lemma
\ref{lemma:boundary:quasisimilarity}, the maps 
\[ \wtpsi = \bdry \wtphi|_{\wtmho}: \wtmho \to \wtV_m \]
and
\[ \psi = \bdry \phi|_\mho: \mho \to V_m \]
are respectively $(C_\epsilon^2, \exp(-\epsilon(n_m+k_m-(n+k))))$- and
$(C_\epsilon^2,
\exp(-\epsilon(n_m-n)))$-quasisimilarities.  So $\wtpsi, \psi$ are $C$-quasisimilarities, where $C = C_\epsilon ^2$ 
and the proof is complete.  
\qed

%\input sec6.tex

%Bibliography

%\bibliographystyle{math}
%\bibliography{refs}
%\input finitetype2.bbl
     \def\cprime{$'$}

\end{document}